\documentclass[10pt,draft]{book}
\usepackage{amssymb}
\usepackage{latexsym}

\newtheorem{theorem}{Theorem}[section]
\newtheorem{lemma}[theorem]{Lemma}
\newtheorem{definition}[theorem]{Definition}
\newtheorem{proposition}[theorem]{Proposition}
\newtheorem{corollary}[theorem]{Corollary}
\newtheorem{notation}[theorem]{Notation}
\newtheorem{conjecture}{Conjecture}
\newtheorem{problem}{Problem}
\newtheorem{condition}{Condition}

\newenvironment{proof}{\paragraph{\it Proof.}}{$\square$\vskip0.4cm}
\newenvironment{remark}{\paragraph{\it Remark.}}{\vskip0.4cm}
\newenvironment{example}{\paragraph{\it Example.}}{\vskip0.4cm}

\newcommand{\nc}{\newcommand}
\nc{\oper}[1]{\mathop{\mathchoice{\mbox{\rm #1}}{\mbox{\rm #1}}
{\mbox{\rm \scriptsize #1}}{\mbox{\rm \tiny #1}}}\nolimits}
\nc{\operlimits}[1]{\mathop{\mathchoice{\mbox{\rm #1}}{\mbox{\rm #1}}
{\mbox{\rm \scriptsize #1}}{\mbox{\rm \tiny #1}}}}
\nc{\M}{{\cal M}}
\nc{\N}{{\cal N}}
\nc{\Q}{{\cal Q}}
\nc{\J}{{\cal J}}
\nc{\higgsinfty}{${{\rm Higgs}_\infty}\ $}
\nc{\Qu}{{\mathbb Q}}
\nc{\calC}{{\cal C}}
\nc{\calD}{{\cal D}}
\nc{\calV}{{\cal V}}
\nc{\calB}{{\cal B}}
\nc{\C}{{\mathbb C}}
\nc{\R}{{\mathbb R}}
\nc{\calK}{{\cal K}}
\nc{\calR}{{\cal R}}
\nc{\Proj}{{\mathbb P}}
\nc{\Z}{{\mathbb Z}}
\nc{\Nat}{{\mathbb N}}
\nc{\E}{{\mathbb E}}
\nc{\U}{{\mathbb U}}
\nc{\dbar}{{\overline{\partial}}}
\nc{\K}{K_\Sigma}
\nc{\Dirac}{{D\!\!\!\! \slash}} 
\nc{\Ha}{{\cal H}}
\nc{\Hy}{{\mathbb H}}
\nc{\uPhi}{{\mathbf \Phi}}
\nc{\cM}{{{\overline\M}}}
\nc{\bM}{{\mathbf M}}
\nc{\calO}{{\cal O}}
\nc{\A}{{\cal A}}
\nc{\bL}{{\mathbb L}}
\nc{\bP}{{\mathbb P}}
\nc{\bV}{{\mathbb V}}
\nc{\calL}{{\cal L}}
\nc{\cD}{{\cal D}}
\nc{\bd}{{\bar{d}}}
\nc{\bg}{{\bar{g}}}
\nc{\G}{{\cal G}}
\nc{\bG}{{\overline{\G}}}
\nc{\Pic}{\mathop{\rm Pic}\nolimits}
\nc{\coker}{\mathop{\rm coker}\nolimits}
\nc{\im}{\mathop{\rm im}\nolimits}
\nc{\rank}{\mathop{\rm rank}\nolimits}
\nc{\ch}{\mathop{\rm ch}\nolimits}
\nc{\td}{\mathop{\rm td}\nolimits}
\nc{\tr}{\mathop{\rm tr}\nolimits}
\nc{\pr}{\mathop{\rm pr}\nolimits}
\nc{\ad}{\mathop{\rm ad}\nolimits}
\nc{\Hom}{\mathop{\rm Hom}\nolimits}
\nc{\End}{\mathop{\rm End}\nolimits}
\nc{\Div}{\mathop{\rm Div}\nolimits}
\nc{\trace}{\mathop{\rm tr}\nolimits}
\nc{\sing}{\mathop{\rm Sing}\nolimits}
\nc{\const}{\oper{const.}}
\nc{\Coeff}{\operlimits{Coeff}}
\nc{\Res}{\operlimits{Res}}
\nc{\cE}{{\cal E}}
\nc{\compE}{{E\stackrel{\Phi}{\rightarrow}{E\otimes K}}}
\nc{\ucompE}{{\E_{\tM}\stackrel{\uPhi}{\rightarrow}
{\E_{\tM}\otimes\K}}}
\nc{\tD}{{\widetilde{D}}}
\nc{\tF}{{\widetilde{F}}}
\nc{\tM}{{\widetilde{\M}}}
\nc{\tN}{{\widetilde{\N}}}
\nc{\D}{{\mathbf D}}
\nc{\cL}{{\cal L}}
\nc{\eqed}{e^\circ_d}
\nc{\teqed}{\widetilde{e}^\circ_d}
\nc{\al}{\alpha}
\nc{\be}{\beta}
\nc{\ze}{\zeta}
\nc{\ga}{\gamma}
\nc{\la}{\lambda}
\nc{\La}{\Lambda}
\nc{\si}{\sigma}
\nc{\xisi}{\xi^\Sigma}
\nc{\sisi}{\si^\Sigma}
\nc{\eqal}{\alpha^\circ}
\nc{\eqbe}{\beta^\circ}
\nc{\eqga}{\gamma^\circ}
\nc{\eqpsi}{\psi^\circ}
\nc{\eqze}{\zeta^\circ}
\nc{\eqeta}{\eta^\circ}

\nc{\binomial}[2]{\mbox{\Large $#1 \choose #2$}}

\nc{\stack}[2]{
{\begin{array}{c}
\scriptstyle #1 \\ \scriptstyle #2 \end{array}} }

\begin{document} 
\pagestyle{empty}

\title{Geometry of the moduli space of Higgs bundles}
\author{Tam\'as Hausel \\ Trinity College, Cambridge}
\date{August, 1998 \vskip3cm 
{\it Thesis submitted for the degree of Doctor of Philosophy \\ in the 
University of Cambridge}}
\maketitle
\newpage
\thispagestyle{empty}
$ $
\newpage
\thispagestyle{empty}
\vspace*{22cm}

{\em\hspace*{13cm} to Kati}
\newpage
\thispagestyle{empty}
\newpage 
\thispagestyle{empty}
$ $
\newpage

\vspace*{5cm} 

\noindent {\huge \bf Preface}

\vspace*{1cm}

\noindent The present dissertation is my own
work, except where attributed to others. It is not the
outcome of work done in collaboration,  except
Chapters~\ref{cohomology} and \ref{minfty}. 

Chapters~\ref{cohomology} and \ref{minfty} describe joint work
with Michael Thaddeus. We started a
correspondence in early 1997 about \cite{thaddeus1} and related
problems. The results appearing in these chapters were mostly achieved
when we participated in the Research in Pairs program in the
Mathematisches Forschungsinstitut Oberwolfach for three weeks in June,
1998. 

I was advised by my supervisor that Chapters 1-5 are sufficient for
a Ph.D. thesis. I added  Chapters \ref{cohomology}
and \ref{minfty} to the thesis 
because I believe that they make the
presented work more compact.

\tableofcontents

\frontmatter

\pagestyle{headings}

\chapter{Acknowledgements}

My principal debt of gratitude goes to my supervisor Nigel Hitchin,
who first suggested the problem of this thesis and then helped in
numerous ways during the course of the research. Many of the ideas
appearing here were found in or after a supervision with him.  

Collaboration with Michael Thaddeus has resulted in
Chapters~\ref{cohomology} and \ref{minfty} of the present thesis. 
His influence is noticeable in
the rest of the thesis as well. 
Working with him proved to be very
stimulating.  

A great amount of inspiration was provided by my college contact Michael
Atiyah through exciting conversations and his own mathematical work.

I learned a lot from Frances Kirwan and Graeme Segal from
their works, lectures and discussions. 

I am grateful to Philip Boalch, Richard Earl, G\'abor Etesi, Manfred Lehn,
Ambrus P\'al, Christian Pauly,
Rich\'ard Rim\'anyi, Bal\'azs Szendr{\H o}i and Andr\'as Szenes, 
who explained to me some parts of
Mathematics appearing in the present thesis. I especially thank
Bal\'azs Szendr{\H o}i for his thorough reading of the draft
versions of this dissertation.  

I was able to write this thesis because I was taught Mathematics by
many excellent teachers throughout the years including K\'aroly
Bezdek, Bal\'azs
Csik\'os, Judit K\"oves,
Endre Makai Jr., G\'abor Moussong, J\'ozsef
Pelik\'an, Istv\'an Reimann,
J\'anos Szenthe, Andr\'as Sz{\H u}cs, Imr\'en\'e Thiry and
J\'anosn\'e Tillman. 

The support of my family in making the thesis and 
helping me to become a mathematician have been indispensable.

I gratefully acknowledge financial support from Trinity College,
Cambridge in the form of an Eastern European Bursary and Research
Scholarship and from the British Government in the form of an Overseas
Research Scheme award. 
I especially thank B\'ela Bollob\'as for introducing me
to Cambridge and Trinity College. 

I spent the academic year 1997-98
in Oxford, as a visiting student, where the hospitality
of the Mathematical Institute and St. Catherine's College, Oxford
helped a lot in finishing the thesis.  

The joint results with Michael
Thaddeus were mostly achieved in the Mathematisches Forschungsinstitut 
Oberwolfach, where we participated in the Research in Pairs program for
three weeks in June, 1998. 
The RiP program is supported by Volkswagen-Stiftung. 
I thank the Institute for its
hospitality which made our research there so fruitful. 

For difficult calculations I used the mathematical software packages
Mathematica, Maple and Macaulay in this reverse lexicographical order. 
In particular, Macaulay 2 proved to be essential in the calculation of the
cohomology ring in Chapter~\ref{cohomology}.

\chapter{Introduction}
\pagestyle{headings}

\section{Motivation: Interaction with Physics}
\label{motivation}

Traditional Mathematical Physics, the subject of mathematically 
rigorous results in Theoretical Physics, has always been a bridge between
Mathematics and Physics. 
Since the late 70's we have been witnessing a profoundly new
interaction between Geometry and Particle Physics, implying an
intrinsic connection between them. 
Let us
mention a few examples: Donaldson's theory of
4-manifolds and Seiberg-Witten theory, enumerative geometry and sigma
models, the Jones polynomial and Chern-Simons theory, the Verlinde
formula and Conformal Field Theory. 
A common feature in these pairs is that the theories are motivated by
both Geometry and Particle Physics. 

In the late 70's physicists started to use sophisticated geometrical
and topological methods for the description of the non-perturbative 
aspects of Quantum Field Theories. It turned out that the
independently developed geometrical methods could provide new insights
into Physics. Indeed, these were used for testing physical theories
in the absence of the available technology for experiments.

Geometers welcomed the renewed interest in their work -- and tried
to solve the new problems suggested by physicists. It resulted in an
exchange of ideas -- physicists revealed new directions in Geometry, and
geometers delivered solutions in Physics.

By now links of this new kind are so numerous that one is tempted to
hope that a new subject is emerging with mixed mathematical and
physical motivations, that may be called Quantum Geometry. 

While the present thesis is concerned only with Mathematics, the motivation
of the original question comes from Physics. In the next section, by explaining
the motivation for the problem of this thesis, we show a few examples  
of the different types of interaction between Geometry and Particle Physics.

\newpage
\section{Statement of the problem}
\label{statement}

Analyzing the conjectured
S-duality in N=2 supersymmetric Yang-Mills theory, 
which is a proposed $SL(2,\Z)$ symmetry of the theory, Sen in \cite{sen} could
predict the dimension  of the space of $L^2$ harmonic forms $\Ha_k$ on 
the universal cover of the 
moduli space of $SU(2)$ magnetic monopoles of charge $k$,
by speculating that there must be an $SL(2,\Z)$-action on the space 
$\bigoplus\Ha_k$,
which represents bound electron states of the theory.

The moduli space of monopoles $M_k$ of charge $k$ is the parameter space of
finite energy and charge $k$ solutions to Bogomolny equations, 
which can be interpreted as the one dimensional reduction of the
self-dual $SU(2)$ Yang-Mills equations on $\R^4$. 

The parameter or moduli space $M_k$ of magnetic monopoles of charge
$k$ is a
non-compact manifold, with $\pi_1(M_k)=\Z_k$, and
has a natural hyperk\"ahler and complete metric on it, which comes from 
an abstract construction 
(the so-called hyperk\"ahler quotient construction\footnote{cf. \cite{hitchin-et-al}}) 
and known
explicitly only in the case $k=2$, when $M_2$ 
is called the Atiyah-Hitchin
manifold\footnote{For further details see Subsection~\ref{gaugebogomolny}.}. 

When $k=2$
Sen's conjecture says that $\dim(\Ha_2)=1$. By knowing the metric of
$M_2$ explicitly, Sen was able to
find a non-trivial $L^2$ harmonic form on the universal cover
$\tilde{M}_2$, giving some support
for his conjecture and in turn for S-duality.   

For higher $k$ however 
Sen's conjecture says something about a metric which is not known 
explicitly. Nevertheless the statement is quite interesting from a 
mathematical point of view as the space of $L^2$ harmonic forms on a 
non-compact complete Riemannian manifold is not well understood. 

Hodge theory tells us that in the compact case the space of $L^2$ harmonic
forms is naturally isomorphic to the De-Rham cohomology of the manifold.
However in the non-compact case 
there is no such theory, and indeed 
the harmonic space depends crucially
on the metric.

Nevertheless some part
of Hodge theory survives for complete Riemannian
manifolds\footnote{cf. 
\cite{derham} Sect. 32 Theorem 24 and Sect. 35 Theorem 26},
such as the Hodge decomposition theorem which states that 
for a complete Riemannian manifold $M$ the space 
$\Omega^*_{L^2}$ of $L_2$ forms on $M$ has an orthogonal decomposition
$$\Omega^*_{L^2}=\overline{d(\Omega^*_{cpt})}\oplus \Ha^* \oplus 
\overline{\delta(\Omega_{cpt}^*)},$$ also that $\Ha^*=\ker(d)\cap \ker(d^*)$.    
An easy corollary\footnote{Cf. \cite{segal-selby}} of these results says that the
composition  
$$H^*_{cpt}(M)\rightarrow \Ha^*\rightarrow
H^*(M)$$ is the forgetful map.

By calculating the image of $H^*_{cpt}(\tilde{M}_k)$ in $H^*(\tilde{M}_k)$ 
Segal and Selby could give a lower bound for the
harmonic forms on the moduli space of magnetic monopoles which coincides
with the dimension given by Sen's conjecture (see \cite{segal-selby}).
This purely mathematical result is thus a supporting evidence for 
the conjectured S-duality in $N=2$ SYM of theoretical Physics. 

In this thesis we will investigate the analogue of Sen's conjecture for
Hitchin's moduli space $\M$ of Higgs bundles of fixed determinant of 
degree $1$ over a Riemann surface
$\Sigma$ of genus $g>1$. 
Hitchin introduced $\M$ in \cite{hitchin1} by considering the $2$-dimensional
reduction of the $4$-dimensional self-dual Yang-Mills equation\footnote{For
details see Subsection~\ref{gaugehitchin}.}. 
The moduli space $\M$ is a simply connected non-compact manifold of
dimension
$12g-12$ with a
complete hyperk\"ahler metric on it.  
Led by the similarities
between the spaces $M_k$ and $\M$ and their origin, we 
address the following problem: 

\begin{problem} What are the $L^2$ harmonic forms on $\M$?
\label{l2}
\end{problem}

In Chapter~\ref{intersection} we solve the topological part of this problem by
considering the so-called {\em virtual Dirac bundle}\footnote{For its gauge
theoretical construction see Section~\ref{gaugedirac}.} $\D$. 
By examining the degeneration locus of $\D$  we prove:

\begin{theorem} The forgetful map 
$$j_\M:H_{cpt}^*(\M)\rightarrow H^*(\M)$$
is $0$. 
\label{main}
\end{theorem}

This says that unlike the case of $\tilde{M}_k$ the topology of $\M$ does not give
the existence of $L^2$ harmonic forms. We can state this fact
informally as: ``There are no topological $L^2$ harmonic forms on 
Hitchin's moduli space of Higgs bundles''. 

Segal and Selby's result together with 
Sen's conjecture suggest that for $\tilde{M}_k$ the topology gives all the
harmonic space. Led by this and supported by the discussion in
Subsection~\ref{toy}
we can formulate the following
conjecture: 

\begin{conjecture} There are no non-trivial $L^2$ harmonic forms on
Hitchin's moduli space of Higgs bundles. 
\label{conjecture}
\end{conjecture}

It would be interesting to see whether a physical argument could back
this conjecture\footnote{Note that the conjecture does not hold for {\em parabolic} Higgs
bundles, as the toy example after Theorem~\ref{fura} shows. Note also that
Dodziuk's vanishing theorem \cite{dodziuk} shows that there are no
non-trivial 
$L^2$ {\em holomorphic} forms on $\M$, 
since the Ricci tensor of a hyperk\"ahler
metric is zero.}. We know of one serious appearance 
of Hitchin's moduli space
of Higgs bundles in the Physics literature. In \cite{bershadsky-et-al}
a topological 
$\sigma$-model with target space $\M$ arises as certain limit of $N=4$
supersymmetric Yang-Mills theory. However it is not clear whether
$L^2$ harmonic forms on $\M$ have any physical interpretation in this theory.

The mathematical counterpart of a topological sigma model with target
space $\M$ is
the enumerative geometry of curves in $\M$. In order to consider
enumerative problems in $\M$ one first has to understand
$H^*(\M)$,
the ordinary 
cohomology ring of $\M$. The purpose of Chapter~\ref{cohomology} is to
describe the cohomology ring of $\M$, establishing the mathematical
background of calculating the physical theory of 
$\cite{bershadsky-et-al}$. This is intended in a future
work. Chapter~\ref{minfty} 
places the understanding of the cohomology ring of $\M$ into a
more general framework, and finds intimate relations with the work of
Atiyah-Bott, and Kirwan. 

All in all the purpose of this thesis is to give a comprehensive
description of the cohomology of Hitchin's moduli space of Higgs
bundles with the aim of using it later for calculating physically
interesting theories.

\newpage
\section{Structure and Results}

The thesis has two parts. The first part is a collection of results
mainly from the literature which are needed in the second part. 
The second part has four chapters. The first two cover more or less the papers
\cite{hausel1} and \cite{hausel2}, respectively, while the last two are made up from
\cite{hausel-thaddeus}.

In Chapter~\ref{compact} the emphasis is on the $\C^*$-action, and the
approach is from Symplectic Geometry. In
Chapter~\ref{intersection}
we focus on the
hypercohomology of Higgs bundles and the point of view is in Algebraic
Geometry. In 
Chapter~\ref{cohomology} the central feature is equivariant cohomology
and the extra tool is Equivariant Topology.
Finally in Chapter~\ref{minfty} the main object is the resolution tower
and the new methods are from Mathematical Gauge theory and Homotopy theory. 

In Chapter~\ref{compact}, using Lerman's construction of symplectic
cutting, we consider a canonical compactification of  
$\M$, producing a projective variety
$\cM=\M\cup Z$, with orbifold singularities. For this we thoroughly examine the $\C^*$-action on $\M$, given
by scalar multiplication of  the Higgs field. We find that the downward flows
correspond to the components of the nilpotent cone while the upward
flows correspond to the Shatz stratification. The former result
will be exploited in Chapter~\ref{intersection}, while the latter in
Chapter~\ref{cohomology} and in Chapter~\ref{minfty}. 
Moreover we give a detailed study of the spaces Z and $\cM$. 
In doing so we reprove some assertions of Laumon and Thaddeus on
the nilpotent cone.

In Chapter~\ref{intersection} we prove the physically motivated 
Theorem~\ref{main}. 
For this we
consider the virtual Dirac bundle on $\M$ 
which is
the analogue of the virtual Mumford bundle on $\N$, the moduli space
of rank $2$ stable bundles of fixed determinant of odd degree over 
a projective curve $\Sigma$.
Using results from Chapter~\ref{compact} about the nilpotent cone
we apply Porteous' theorem to the downward flows to obtain a proof of Theorem~\ref{main}
that all 
intersection numbers in the compactly supported
cohomology of $\M$ vanish, i.e. ``there are no {\em topological} 
$L^2$ harmonic forms on $\M$''. 
Our result
generalizes two facts. One is the well known vanishing of the Euler characteristic of
$\N$, which gives the vanishing of one intersection number on
$\M$. The other is the vanishing of the ordinary cohomology class of
the Prym variety, the generic fibre of the Hitchin map, which gives
the vanishing of $g$ intersection numbers on $\M$. We prove here that the
rest of the $g^2$ intersection numbers also vanish.   
Our proof
shows that the vanishing of all intersection numbers of 
$H^*_{cpt}(\M)$ is given by
relations analogous to the Mumford relations in the cohomology ring of
$\N$. 

In Chapter~\ref{cohomology} we give a conjectured complete description 
of the cohomology ring of $\M$. We use an approach motivated by
Kirwan's work on the proof of the Mumford conjecture, to prove 
that the equivariant cohomology ring of $\tM$ is generated by
universal classes. However our proof is completely geometric in
nature and rests on the degeneracy locus description of the upward flows.
We conclude the chapter by explaining a -- computer supported -- conjectured complete description
of $H^*_I(\M)$, the subring of $H^*(\M)$ generated by the universal
classes $\alpha$, $\beta$ and $\gamma$. We support this conjecture by
proving that the conjectured ring and $H^*_I(\M)$ have the same  
Poincar\'e polynomial. Also we find the first
two relations. The second of which turns out to be $\be^g=0$,
showing that Newstead's conjecture is still true over $\M$. 

In Chapter~\ref{minfty} we construct a {\em resolution tower} of 
smooth $S^1$-manifolds: $$\tM\cong\tM_0\subset \tM_1\subset ... \subset \tM_k \subset
...,$$ from the moduli spaces $\tM_k$ of stable Higgs bundle with a 
pole of order at most $k$ at a fixed point, and of
degree $1$, and consider the direct limit  $$\tM_{\infty}=\lim_{k\to
\infty}\tM_k.$$
We prove that its cohomology is a free graded algebra on the universal
classes and that $$i_0^*:H^*(\tM_\infty)\to H^*(\tM)$$ is
surjective, therefore  a resolution of the cohomology ring of $\tM$. 
We show how $\tM_\infty$ can be used to provide a 'finite
dimensional' and purely geometric proof of the Mumford conjecture. 
To shed light on many striking features of $\tM_\infty$ we show that 
$\tM_{\infty}$ is homotopy
equivalent to $B{\bG}$, the classifying space of the gauge group
modulo constant scalars, and that they are also homotopy equivalent as
stratified spaces. We finish the chapter by proving that even the homotopy of the
resolution  tower stabilizes in the spirit of the Atiyah-Jones conjecture.

We conclude the thesis by summarizing the work in the previous
chapters from the point of view of the compactification $\cM$.

\mainmatter

\part{General facts}
\label{general}

\chapter{Moduli spaces}

The subject of the thesis is the investigation of a particular moduli
space, the moduli space of Higgs bundles. In this chapter we give an
introduction to moduli spaces in general and collect results from the
literature which
will be needed in the second part. The chapter has two sections. In
the first one we deal with gauge theory in the second with 
algebraic moduli spaces.

\section{Mathematical gauge theory}

The interests of mathematicians in gauge theory started in the late
70's 
with the appearance  of a few pioneer papers, such as  
\cite{atiyah-etal-1} and
\cite{atiyah-etal-2}. The papers were concerned with the Yang-Mills
equations in gauge theory, which for mathematicians meant a branch
of differential geometry, namely connections on fibre bundles.

Yang and Mills introduced the Yang-Mills equations on Minkowski
$4$-space in 1954 as a 
non-Abelian generalization of Maxwell's equations. Later the same
equations over Euclidean $4$-space were also considered and these  
equations are  our concern in the next subsection\footnote{For a
more detailed introduction see \cite{atiyah}.}.

\subsection{Yang-Mills equations in $4$ dimensions}

Let $G$ be a Lie group (usually $U(2)$, $SU(2)$ or $SO(3)$) and $P$ be a
principal $G$-bundle over $\R^4$. If $A$ is a connection on $P$, then its curvature  
$F(A)\in \Omega^2(\R^4;\ad(P))$ is a two-form with values in
$\ad(P)$, where $\ad(P)=P\times_G{\mathfrak g}$ is the vector
bundle
associated to the adjoint representation. 

In principle a physical theory is given by its Lagrangian. 
The {\em Yang-Mills Lagrangian} (or {\em functional} or {\em energy function} or {\em action}) is:
\begin{eqnarray} S(A)=-\int_{\R^4}\tr(F(A)\wedge
\ast F(A)),\label{functional}\end{eqnarray} where $$\ast:\Omega^2(\R^4;\ad(P))\to 
\Omega^2(\R^4;\ad(P))$$ is the Euclidean Hodge star operator. 

The corresponding Euler-Lagrange equations, which describe the
critical points of the functional $S$, are called 
the {\em Yang-Mills equations}:  \begin{eqnarray} d_A\ast F(A)=0, \label{yangmills}\end{eqnarray}
Recall that the Bianchi identity says that \begin{eqnarray} d_A F(A)=0,
\label{bianchi} \end{eqnarray}
which is formally similar\footnote{As a matter of fact
 a possible duality between the two equations is the source of the conjectured
S-duality, which was mentioned in the Introduction. In Maxwell's
 theory, which is a Yang-Mills theory on Minkowski $4$ space with
 $G=U(1)$,  this is the well known duality between electricity and
 magnetism.\label{sduality}}  to (\ref{yangmills}).  

The absolute minima of the Yang-Mills Lagrangian are given
by the {\em self-dual} 
$$F(A)=\ast F(A)$$ 
and {\em anti self-dual} 
$$F(A)=-\ast\! F(A)$$ 
{\em Yang-Mills equations}. Physically a solution of finite energy 
to these equations
represents a {\em multi-instanton}. Note that because of the Bianchi
identity (\ref{bianchi}) it is immediate that 
solutions to the last two equations  
satisfy the original Yang-Mills equations (\ref{yangmills}).
Mathematically the last two equations are equivalent by reversing the
orientation of the underlying $\R^4$. We will only consider the
self-dual equations here. 

We will need later the form of the equations above in terms of a
trivialization of 
$P$ over $\R^4$. If the basic
coordinates of $\R^4$ are  $(x_1,x_2,x_3,x_4)$ , then  the connection $A$ is
described
by a Lie algebra-valued $1$-form: 
\begin{eqnarray} A=A_1dx_1+A_2dx_2+A_3dx_3+A_4dx_4.\label{connection}
\end{eqnarray}  
and then its curvature is $$F(A)=dA+A\wedge A.$$
Alternatively $$F(A)=\sum_{i<j} F_{ij}dx_i\wedge\ dx_j,$$
where $$F_{ij}=\left[\frac{\partial}{\partial x_i}+A_i,
\frac{\partial}{\partial x_j}+A_j\right],$$
or if we write $$\nabla_i=\frac{\partial}{\partial x_i}+A_i,$$
then $$F_{ij}=[\nabla_i,\nabla_j].$$

In this trivialization the self-dual Yang-Mills equations then take the form
\begin{eqnarray}
\begin{array}{c}
F_{12}=F_{34},\\ 
F_{13}=F_{42},\\ 
F_{23}=F_{14}. \end{array} {\Bigg \} } 
\label{selfdual}
\end{eqnarray} 

An important aspect of the Yang-Mills equations is that they are {\em gauge
invariant}. 
To explain this consider a $C^{\infty}$ section $g$ of the bundle of 
groups $P\times_{\rm Ad} G$.  It is an automorphism of the principal bundle which leaves
each fibre invariant. It is called a {\em gauge
transformation}, and the group of all gauge transformations is called the
{\em gauge group}. 
A gauge  transformation transforms any connection $A$ on $P$ in a
natural way to another connection $g(A)$. Moreover 
 if $A$ is a
solution to the self-dual Yang-Mills equations then $g(A)$ is also a
solution. 
Thus $\G$ acts naturally on the finite energy 
solution space of the self-dual
Yang-Mills equations. Because we do not consider two gauge equivalent
solutions different we define the {\em moduli space} (or {\em
parameter space}) to be the quotient.  

It can be shown that in nice cases the moduli space will inherit a
natural structure of a finite dimensional smooth manifold. 
The study of these moduli spaces is the central problem of the
mathematical
gauge theory.  The significance of the
subject became more apparent in the early 80's 
with the seminal work of Donaldson and Atiyah-Bott.

Donaldson, analyzing the moduli space of solutions to the anti-self dual
Yang-Mills equations over a closed  oriented Riemannian $4$-manifold $M$, could prove 
many fundamental theorems about the topology of the underlying
differentiable manifold\footnote{For more details consult \cite{donaldson-kronheimer}.}.    
     
Atiyah and Bott, in \cite{atiyah-bott}, 
extensively studied the $2$-dimensional Yang-Mills equations, relating
the subject to stable vector bundles on algebraic curves and Morse theory. 
In the next subsection we explain  how they related Yang-Mills theory to
the algebraic geometry of vector bundles on curves.

\subsection{Yang-Mills equations in $2$ dimensions}
\label{yangmills2}

For the sake of simplicity we restrict our attention to $G=U(n)$, and
let $P$ be a principal $U(n)$-bundle over a Riemann surface
$\Sigma$. We also fix a metric on $\Sigma$ compatible with the
complex structure. 
Now if $A$ is a connection on $P$ its curvature is: $$F(A)\in \Omega^2(\R^4;\ad(P)).$$
The space $\A$ of all connections on $P$ is naturally an affine space 
modelled on the infinite dimensional vector space 
$\Omega^{0,1}(\Sigma,\ad(P)\otimes\C)$. 
As we explained above,
the infinite dimensional 
gauge group $\G=\Gamma(\Sigma,P\times_{\rm Ad}U(n))$ 
acts naturally on the infinite dimensional affine space $\A$.  

The strategy of \cite{atiyah-bott} is to calculate
$H^*_\G(\A)$, the
$\G$-equivariant cohomology\footnote{For the definition of equivariant cohomology see
Subsection~\ref{equivariantcoh}.} 
of $\A$, in two different ways.  
First is a direct way by noting that $\A$, being
homeomorphic\footnote{To get the correct topology on $\A$ as a Banach
manifold, one has to
consider Sobolev connections. Because we want to give only an intuitive
picture we do not spell out the details here, but refer to $\S 14$ of
\cite{atiyah-bott}.
For related remarks see the
footnotes at the end of  
Subsection~\ref{gaugemk}.} to a
topological vector space, is contractible, thus $$H^*_\G(\A)\cong
H^*(B\G),$$ and $H^*(B\G)$ can be calculated
directly\footnote{For the result see Section~\ref{cohomologyofN}.}.

The other approach for calculating $H^*_\G(\A)$ is Morse
theory. Morse theory gives information about the cohomology of the 
whole space in terms of the cohomology of the non-degenerate critical submanifolds
of the Morse function. 
Recall that  we have the  Yang-Mills functional
$S:\A\rightarrow \R$, defined by (\ref{functional}), which gives a function
on $\A$. 

The key idea of \cite{atiyah-bott} is to think of $S$ as a
$\G$-equivariant Morse function on the 
infinite dimensional Banach manifold
$\A$, and calculate $H^*_\G(\A)$ from Morse theory. 
As in the previous subsection,
critical points of $S$ are given by solutions to the $2$-dimensional 
Yang-Mills equations (the Euler-Lagrange equations of the Yang-Mills
functional)
of the form
(\ref{yangmills}) with 
$$\ast:\Omega^2(\Sigma;\ad(P))\to 
\Omega^0(\Sigma;\ad(P))$$
now being the Hodge star operator corresponding to the fixed metric on
$\Sigma$. 

The solution space has infinitely many components: one component,
corresponding to the absolute minimum of the Yang-Mills functional,
contains  
irreducible connections,  whereas the others contain only reducible ones. A famous
result of Narasimhan and Seshadri identifies the moduli space of
irreducible solutions to the Yang-Mills equations with the 
moduli space of stable vector bundles on $\Sigma$, a space, which had
been investigated by algebraic geometers for decades. 

Thus,
provided that the Morse function is perfect, the
Morse theory approach gives a method to calculate the cohomology of
the moduli space of stable bundles,
by knowing the cohomology
of the moduli space of reducible connections, which can be inductively
calculated via this process. 
Though they couldn't proceed with the analysis\footnote{Later this
was done by Daskalopoulos in \cite{daskalopoulos}.}, they found an
alternative way. Namely there is a one-to-one correspondence between
unitary connections on $P$ and holomorphic structures on the
associated $C^\infty$ 
complex vector bundle $\calV=P\times_{\ad}{\mathfrak g}$. Thus
one can identify $\A$ with the complex affine
space $\calC$ of holomorphic structures 
on ${\calV}$. Moreover the latter space has the advantage of being independent
of
the metric on $\Sigma$, transferring the problem from differential
geometry to holomorphic, and indeed algebraic geometry.  
 
In Section~\ref{cohomologyofN} we will explain how Atiyah and Bott
applied Morse theory to $\calC$, motivated by the above heuristic argument of
Morse theory over $\A$.

\subsection{Bogomolny equations in $3$ dimensions}
\label{gaugebogomolny}

In the 80's, further parts of gauge theory became the subject of 
study by mathematicians. A good example is the problem of magnetic
monopoles. To explain the mathematical background, consider a connection $A$
on a principal $G$-bundle $P$ over $\R^4$ of the form (\ref{connection}). 
If we make the assumption that the Lie algebra-valued functions $A_i$
are independent of $x_4$, then $A_1$,$A_2$ and $A_3$ define a
connection: $$A_1dx_1+A_2dx_2+A_3dx_3$$ over $\R^3$, while $A_4$ becomes a
function on $\R^3$ traditionally denoted by $\phi$ and called the
{\em Higgs field}. 

Under this procedure the self-dual Yang-Mills equation reduces to
the so-called {\em Bogomolny equation}: $$F(A)=\ast d_A\phi.$$ 
Moreover the reduced energy function takes the form
$$S(A)=-\int_{\R^3}\tr(F(A)\wedge \ast F(A))+\tr(d_A\phi\wedge \ast d_A\phi).$$  

This procedure is called {\em dimensional reduction}. Thus 
the Bogomolny equation is the $1$-dimensional reduction of the
self-dual Yang-Mills equation.

Physically a solution of finite energy to the Bogomolny equation
represents a {\em magnetic monopole}. It can be shown that any
solution to the Bogomolny equation with finite energy has energy 
$8\pi k$, where $k$ is a positive integer, which is called the {\em charge}
of the monopole. 

As in the Introduction we 
denote by $M_k$ the moduli space of charge $k$ and $SU(2)$ magnetic monopoles,
which is the charge $k$ solution space to the $SU(2)$ Bogomolny
equation modulo gauge transformations. This is a non-compact
smooth manifold of dimension $4k-4$, with an inherited complete
hyperk\"ahler metric. 

Atiyah and Hitchin could find the metric explicitly on the
$4$-dimensional 
manifold $M_2$ and by examining its geodesics they
could provide\footnote{For further details see \cite{atiyah-hitchin}.}
a description of the 
low energy scattering 
of $SU(2)$ magnetic monopoles of charge $2$.

\subsection{Hitchin's self-duality equations in $2$ dimensions}
\label{gaugehitchin}

In 1987 Hitchin in \cite{hitchin1} considered the $2$ dimensional
reduction
of the self-dual
Yang-Mills equations, whose moduli space is the central object 
of the present thesis. We now describe this construction.

If we assume (further) that the Lie algebra valued functions $A_i$
in (\ref{connection}) are independent of both $x_3$ and $x_4$, then
$A_1$ and $A_2$ define a connection $$A=A_1dx_1+A_2dx_2$$ on $\R^2$,
while $A_3$ and $A_4$ become Lie algebra-valued functions on $\R^2$,
which we relabel as $\phi_1$ and $\phi_2$. Furthermore we introduce
$\phi=\phi_1-i\phi_2$, called the {\em complex Higgs field}.

From a coordinate independent point of view, we have a connection $A$
on a principal $G$-bundle $P$ over $\R^2$ together with an auxiliary field
$$\phi\in\Omega^0(\R^2;\ad(P)\otimes \C).$$ If we moreover write 
$z=x_1+ix_2$ and introduce $$\Phi=\frac{1}{2}\phi dz\in
\Omega^{1,0}(\R^2;
\ad(P)\otimes \C)$$ and $$\Phi^\ast=\frac{1}{2}\phi^\ast d\bar{z}\in \Omega^{0,1}(\R^2;\rm
ad(P)\otimes \C)$$ then the $2$-dimensional reduced self-dual Yang-Mills equation
becomes 

\begin{eqnarray}
\begin{array}{c} F(A)=-[\Phi,\Phi^*],\\
d_A^{\prime\prime}\Phi=0.\end{array}{\Bigg \} } \label{equation} 
\end{eqnarray}
These equations are called {\em Hitchin's self-duality equations}.

Unfortunately there are no finite energy solutions to 
Hitchin's
self-duality 
equations on $\R^2$. However by exploiting the conformal invariance of
(\ref{equation}) we can write down  Hitchin's self-duality 
equations
over a compact Riemann surface $\Sigma$ by demanding $A$ to be a
connection on a principal $G$-bundle $P$ over $\Sigma$ and 
the Higgs field to be $$\Phi\in \Omega^{1,0}(\Sigma;\ad(P)\otimes \C).$$

Fortunately solutions to Hitchin's self-duality equations over a
Riemann surface do exist and Hitchin made an extensive study of
their parameter or moduli space $\M$ in \cite{hitchin1} for the case
of $G=SU(2)$ and $G=SO(3)$. We list some of his results in Subsection~\ref{nagysrac}.

\subsection{Dirac equations in $2$ dimensions}
\label{gaugedirac}

To explain the gauge theoretic origin of the virtual Dirac bundle
$\D$ of Chapter~\ref{intersection}, following Hitchin
\cite{hitchin4}, 
we first consider the Dirac 
equation in $\R^4$.

Let $\psi_1$ and $\psi_2$ be scalar functions.
The ordinary {\em Dirac equation}
is of the form:

$$
\left(\frac{\partial}{\partial x_1}+\left(\begin{array}{cc}i&0\\
0&-i\end{array}\right)\frac{\partial}{\partial x_2}+
\left(\begin{array}{cc}0&1\\ -1&0\end{array}\right)
\frac{\partial}{\partial x_3}+
\left(\begin{array}{cc}0&i\\ i&0\end{array}\right)\frac{\partial}{\partial x_4}\right)
\left(\begin{array}{c}\psi_1\\ \psi_2\end{array}\right)=0.
$$

If $A$ is a self-dual Yang-Mills connection on $P$, and $\calV$ is the
vector bundle associated to $P$ in a vector representation, 
then the Dirac equation 
coupled to $A$ is:
$$\Dirac_A(\psi)=0,$$
where $\psi$ is a twisted spinor and 
the Dirac operator coupled to $A$ 
$$\Dirac_A: \Gamma(S_+\otimes \calV)\to \Gamma(S_-\otimes \calV)$$ is given by 
$$\Dirac_A=\nabla_1+\left(\begin{array}{cc}i&0\\
0&-i\end{array}\right)\nabla_2+\left(\begin{array}{cc}0&1\\ -1&0\end{array}\right)
\nabla_3+
\left(\begin{array}{cc}0&i\\ i&0\end{array}\right)\nabla_4,$$
where   $$\nabla_i=\frac{\partial}{\partial x_i}+A_i.$$ 

After the dimensional reduction we consider a connection on
a principal $G$-bundle $P$ over a
Riemann surface $\Sigma$, the rank $2$ vector bundle $\calV$ is
associated to $P$ in a suitable vector representation and $\phi$ a complex Higgs field. 
Then the above coupled Dirac equation takes the following
shape: 
\begin{eqnarray}\begin{array}{c}\bar{\partial}_A\psi_1+\phi\psi_2=0\\ 
\partial_A\psi_2+\phi^*\psi_1=0,\end{array}{\Bigg \} }\label{diracequ}\end{eqnarray} 
now with $$\psi_1\in  \Omega^{1,0}(\Sigma,\calV)$$ and  
$$\psi_2\in  \Omega^{0,1}(\Sigma,\calV).$$ 
Suppose that $(A,\phi)$ is a solution to Hitchin's self-duality
equations (\ref{equation}) and consider the vector space of solutions 
to the equation (\ref{diracequ}). Using Hodge theory of elliptic
complexes
it can be shown that this vector
space is canonically isomorphic to the Dolbeault definition of the  
hypercohomology\footnote{For the {\v C}ech definition see Subsection~\ref{hypercohomology}}
vector space 
$\Hy^1(\Sigma,E_A\stackrel{\Phi}{\rightarrow}E_A\otimes K)$, 
where $E_A$ is the holomorphic vector bundle corresponding to the
connection $A$ on $\calV$, and $\Phi\in
H^0(\Sigma,\End (E)\otimes K)$ is the corresponding Higgs field. 
Thus we can assign to any point $(A,\Phi)\in \M$ the complex
vector space 
$\Hy^1(\Sigma,E_A\stackrel{\Phi}{\rightarrow}E_A\otimes K)$. This will
not be a vector bundle in general but only a coherent sheaf $\D$, which we call
the {\em virtual Dirac bundle}. 
We will define $\D$ rigorously in
Chapter~\ref{intersection} and use it to prove Theorem~\ref{main}.

\newpage
\section{Algebraic Geometry of vector bundles on curves}

As we already mentioned \cite{atiyah-bott}  established a link
between Yang-Mills theory on Riemann surfaces and the algebraic geometry
of vector bundles over projective algebraic curves via the theorem of
Narasimhan and Seshadri. 
Recall that, in a suitable form, this theorem asserts
that there is a one-to-one correspondence between rank $n$ semi-stable
vector bundles and 
those connections on a principal $U(n)$-bundle, which give 
absolute minima of the Yang-Mills functional (\ref{yangmills}).    

In this section we recall some
definitions and results from the theory of
vector bundles over projective algebraic curves, which we will need
later.

\subsection{Moduli spaces of stable bundles}
\label{kissrac}

Algebraically a Riemann surface $\Sigma$ is the same as a non-singular
complex projective algebraic curve. Moreover holomorphic vector bundles over $\Sigma$
correspond to complex algebraic vector bundles. We will not
distinguish between them. 

The main problem of the subject is to classify all vector
bundles\footnote{A vector bundle always means an algebraic vector
bundle, unless otherwise stated.} over the curve
$\Sigma$. This was done for genus $0$ curves by Grothendieck in
\cite{grothendieck}, and
for genus $1$ curves by Atiyah in \cite{atiyah2}, both in
1957. However the
genus at least $2$ case has proved to be much more difficult. From now on we  restrict our
attention to the genus at least $2$ case.   

%
%

First recall that $C^\infty$ vector bundles over $\Sigma$ are
classified by their ranks and degrees.  Thus we will concentrate on
$\calC$, the complex affine space of holomorphic
structures
on a fixed $C^\infty$ complex vector bundle ${\calV}$ over
$\Sigma$ of rank $n$ and degree $d$. Certainly we do not want to distinguish between isomorphic 
vector bundles, therefore we consider the complexified gauge group 
$\G^c={\rm Aut}({\calV})$ of complex automorphisms of ${\calV}$, which acts
naturally on $\calC$. Since an orbit of this action is clearly the set of
vector bundles isomorphic to a given one, the main problem reduces to understand 
the orbits of this action. However the space $\calC/\G^c$, though
clearly parametrizes all vector bundles with underlying $C^\infty$
bundle $\calV$, is not even Hausdorff\footnote{This is because  of the
so-called {\em jumping phenomenon}.}, hence its description is 
rather hopeless. One attempt of getting around this problem was
provided by Mumford's Geometric Invariant
Theory\footnote{Cf. \cite{fogarty-etal}} in the 60's.      
The idea is to take a $\G^c$-invariant open subset $\calC_{s}$ 
of $\calC$ in order to
get a Hausdorff space $\calC_{s}/\G^c$, and indeed, as we will see
later, a smooth algebraic
variety, and in good cases a projective one. 

%
%

\begin{definition} The {\em slope} of a vector bundle $E$ 
on $\Sigma$ is defined by: $$\mu(E)={\rm deg}(E)/\rank(E).$$ Moreover 
$E$ is {\em semi-stable} (resp. {\em stable}) 
if it has at least as large (resp. strictly larger) slope than any of its proper
subbundles. Finally $\calC_{ss}\subset \calC$ denotes the open set of
semi-stable, $\calC_{s}\subset \calC$ the open set of stable bundles. 
\label{slope}
\end{definition}    

In some sense semi-stable bundles are the analogues of simple finite groups in
the classification of finite groups. For example one has the analogue
of the Jordan-H\"older theorem:

\begin{theorem}[Harder, Narasimhan] Every holomorphic bundle $E$ has
a {\em canonical filtration}: \begin{eqnarray}
0=E_0\subset E_1 \subset E_2 \subset \dots \subset
E_r=E,\label{filtration} \end{eqnarray}
with $D_i=E_i/E_{i-1}$ semi-stable and $$\mu(D_1)>\mu(D_2)>\dots >
\mu(D_r).$$
\label{harder-narasimhan}
\end{theorem}
Thus as in the classification of finite groups, one first starts with
the classification of semi-stable vector bundles.
As we mentioned above, Geometric Invariant Theory constructs a
projective algebraic variety for the parameter space of semi-stable
vector bundles:

\begin{theorem} When $(n,d)=1$ then $\calC_s=\calC_{ss}$ and 
the moduli space  $N(n,d)=\calC_{s}/\G^c$
of rank $n$ and degree $d$ 
stable bundles over $\Sigma$ is a smooth projective
algebraic variety of dimension $n^2(g-1)+1$.
\end{theorem}

\begin{example} 1. Clearly every line bundle is stable thus $N(1,d)$
contains all line bundles of degree $d$. 
This space is the well-known Jacobian, which
we relabel as $\J_d=N(1,d)$, and use $\J$ for $\J_1$. This is an Abelian variety of dimension
$g$. In this case the classification problem is clearly
settled. 

2. The next moduli spaces in the list are the rank $2$ moduli spaces
$N(2,d)$, with $d$ odd. These are all isomorphic, thus we restrict our
attention to $N(2,1)$, which we rename as $\widetilde{\N}$. \label{tN}
It is a smooth projective variety of dimension
$4g-3$. The determinant gives a map 
$det_\N:\widetilde{\N}\rightarrow \J$. For any 
$\Lambda\in \J$ the fibre $det_\N^{-1}(\Lambda)$ will be
denoted by $\N_{\Lambda}$, which is a smooth projective variety
of dimension $3g-3$. The map $f:\N_{\Lambda_1}\rightarrow
\N_{\Lambda_2}$ given by $f(E)= E\otimes (\Lambda_2\otimes
\Lambda_1^*)^{1/2}$, where $(\Lambda_2\otimes \Lambda_1^*)^{1/2}$ is a
fixed square root of $\Lambda_2\otimes \Lambda_1^*$, is an isomorphism
between $\N_{\Lambda_1}$ and $\N_{\Lambda_2}$. Hence we will write
$\N$ \label{N} for $\N_\Lambda$, when we do not want to emphasize the fixed line bundle
$\Lambda$. In words $\N$ is the moduli space of stable rank $2$ bundles of
fixed determinant of degree $1$ over the Riemann surface $\Sigma$.

The moduli spaces $\N$ and $\widetilde{\N}$ 
will appear all along in this thesis. They have been much
studied over the years. 
In particular their
cohomology rings have 
been described completely. We explain some of the   results in
this direction in Section~\ref{cohomologyofN}.
\end{example}

\subsection{Moduli spaces of Higgs bundles}
\label{nagysrac}

Hitchin in Theorem 4.3 of \cite{hitchin1} proved a generalization of
the above mentioned theorem of 
Narasimhan and Seshadri, which linked the solutions of his
self-duality equations (\ref{equation}) to algebro-geometric objects,
so-called stable Higgs  bundles\footnote{For a more thorough treatment see Subsection~\ref{vanthe}.}: 

\begin{definition} 
The complex $\compE$ with $E$ a
vector bundle on $\Sigma$, $K$ the canonical bundle of $\Sigma$, 
and $\Phi\in H^0(\Sigma,\Hom(E,E\otimes K))$, 
is called a {\em Higgs bundle}\footnote{The term {\em Higgs bundle} was first
used by Simpson in
\cite{simpson}.}, while $\Phi$ is called the 
{\em Higgs field}. 

The slope 
$\mu(\cal E)$ of a Higgs
bundle $\cE=\compE$ 
is defined as the slope\footnote{Cf. Definition~\ref{slope}} $\mu(E)$
of its vector bundle $E$. 
A Higgs
bundle is called {\em semi-stable}  (resp. {\em stable}) 
if it has at least as large (resp. strictly larger) slope 
than any of its proper $\Phi$-invariant subbundles.  
\label{stablehiggs}
\end{definition}

Thus Hitchin's gauge theoretic construction of $\M$, the moduli space
of solutions to Hitchin's self-duality equations, with fixed
determinant connection as explained in \cite{hitchin1}, in turn is
isomorphic to the moduli space of rank $2$ stable Higgs bundles with
trace-free Higgs field and 
fixed determinant of odd degree. 
This space is the central moduli
space of the present thesis: We fix a degree $1$ line bundle $\Lambda$ over
$\Sigma$, and denote by $\M$ \label{M} the moduli space of rank $2$ 
stable Higgs bundles with trace free Higgs field and determinant $\Lambda$. 

The moduli space $\M$ has many features which show its
importance. Probably the most important is that  the cotangent bundle
of $\N$, (defined above)  sits inside
$\M$ as an open dense subset. Namely, $(T^*_\N)_E$ is canonically isomorphic
to $H^0(\Sigma, \End_0(E)\otimes K_\Sigma)$ thus the points of $T^*_\N$
are stable Higgs bundles. 
 
After introducing the space $\M$, Hitchin gave its extensive description
in \cite{hitchin1}, \cite{hitchin2}. 
Here we restate those results, which we use later.

\begin{itemize}
\item 
$\M$ is a noncompact,
smooth complex manifold of complex dimension $6g-6$ containing $T^*_\N$ as a dense
open set.

\item
Furthermore $\M$ is canonically a Riemannian manifold with a complete
hyperk\"ahler metric. Thus $\M$ has complex structures parameterized
by $S^2$. One of the complex structures, for which $T^*_\N$ is a complex
submanifold, is distinguished, call it $I$. 
We will only be concerned with this
complex structure here. The others (apart from $-I$) are 
biholomorphic to each other and
give $\M$ the structure of a Stein manifold. From the corresponding 
K\"ahler forms one
can build a holomorphic symplectic form $\omega_h$ on ($\M,I$). 

\item

There is a map, called the Hitchin map $$\chi:\M\rightarrow
H^0(\Sigma,K^2)=\C^{3g-3}$$ defined by
$$(E,\Phi)\mapsto \det \Phi.$$ The Hitchin map is proper and an algebraically
completely integrable Hamiltonian system with respect to the holomorphic
symplectic form $\omega_h$, with generic fibre a Prym variety corresponding
to the spectral cover of $\Sigma$ at the image point.

\item

Let $\omega$ denote the K\"ahler form corresponding to the complex structure
$I$. There is a holomorphic $\C^*$-action on $\M$ defined by 
$(E,\Phi)\mapsto (E,z\cdot \Phi)$. The restricted action of $U(1)$ defined by
$(E,\Phi)\rightarrow (E,e^{i\theta}\Phi)$ is isometric and indeed
Hamiltonian with proper moment map $\mu$. The function 
$\mu$ is a perfect Morse function,
moreover:

$\mu$ has $g$ critical values: an absolute minimum $c_0=0$ and
$c_d=(d-\frac{1}{2})\pi$, where $d=1,...,g-1$.

$\mu^{-1}(c_0)=\mu^{-1}(0)=F_0=\N$
is a non-degenerate critical manifold of index $0$.

$\mu^{-1}(c_d)=F_d$ is a non-degenerate critical manifold of index
$2(g+2d-2)$ and is diffeomorphic to a $2^{2g}$-fold cover of 
the symmetric product $\Sigma_{\bd}$, where we used the notation $\bd=2g-2d-1$.

\item

The fixed point set $S$ of the involution $\sigma(E,\Phi)=(E,-\Phi)$ is
the union
of $g$ complex submanifolds of $\M$ namely, 
$$S=\N\cup \bigcup_{d=1}^{g-1}E^2_d,$$
where $E^2_d$ is the total space of a vector bundle $E^2_d$ over $Z_d$. Moreover,
$E^2_d$ is a complex submanifold of dimension $3g-3$.

\end{itemize}

\subsection{The moduli space of Higgs $k$-bundles}

Using Geometric Invariant Theory, Nitsure in \cite{nitsure} gave an
algebraic construction of $\M$ and many other related spaces. As we
will use some of them in Chapter~\ref{cohomology}, we define them here:

\begin{definition} Let $k\geq 0$ and $\tM_k$ \label{tMk} denote the moduli space of stable
rank $2$ Higgs bundles of degree $1$, with poles of order at most $k$ at a fixed
point $p\in \Sigma$. A Higgs bundle with pole is a complex
$E\stackrel{\Phi}{\rightarrow} E\otimes K\otimes L_p^k$ where $E$ is a rank
$2$ vector bundle over $\Sigma$, the line bundle $L_p$ corresponds to
the divisor $p\in \Sigma$ and the Higgs field with poles: $\Phi\in
H^0(\Sigma,\End(E)\otimes K \otimes L_p^k)$. For convenience we call
such a complex a {\em Higgs $k$-bundle} and $\Phi$ a {\em Higgs
$k$-field}. 
Moreover we call $\compE\otimes L_p^k$ stable if the slope of any 
$\Phi$-invariant line subbundle of $E$ is strictly smaller than $\mu(E)$. 
\label{poles}
\end{definition}

Proposition 7.4 of \cite{nitsure} then tells us:

\begin{theorem}[Nitsure] The space $\widetilde{\M}_k$ is a smooth
quasi-projective variety of dimension
$$8(g-1)+1+4k+\dim(H^1(\Sigma,K\otimes L_p^k)).$$ 
\end{theorem}

As in the case of $\widetilde{\N}$, the determinant map gives a map 
$$det_{\M_k}:\tM_k\rightarrow \J\times H^0(\Sigma,K\otimes L_p^k),$$ defined by 
$det_{\M_k}(E,\Phi)=(\Lambda^2E,\trace(\Phi))$. 
For any 
$\cL\in \J\times H^0(\Sigma,K\otimes L_p^k)$ the fibre $det_{\M_k}^{-1}(\cL)$ will be
denoted by $\M^k_{\cL}$. Just as in the stable vector bundle case any two fibres of
$det_{\M_k}$
are 
isomorphic. Usually we will write $\M_k$ \label{Mk} for
$\M^k_\cL$, when the Abelian Higgs bundle, with order $k$ pole,  
$\cL$ has zero Higgs field. 
The dimension of $\M_k$ is clearly
$$8(g-1)+1+4k+\dim(H^1(\Sigma,K\otimes L_p^k))-
(g+\dim(H^0(\Sigma,K\otimes L_p^k)))= 6(g-1)+3k.$$ By definition
$\M_0=\M$, hence $\dim(\M)=6(g-1)$, which checks up with the dimension
calculated by Hitchin. We will also use $\tM$ \label{tM} for $\tM_0$.

A consequence of the Geometric Invariant Theory construction of
Nitsure is the following\footnote{Cf. Remark 5.12 in \cite{nitsure}.}:

\begin{corollary} The spaces $\tM_k$ and consequently $\M_k$ are
quasi-projective varieties.
\end{corollary}

In particular Hitchin's moduli space $\M$ is also a quasi projective
variety.
We will examine in a certain sense the
canonical compactification of $\M$ in Chapter~\ref{compact}. 

Before we
proceed to the next section we insert here a gauge theoretic
construction of Nitsure's spaces $\tM_k$, which is more in the
spirit of \cite{atiyah-bott} and \cite{hitchin1} and shall be used
throughout the thesis.  

\subsubsection{Gauge theoretic construction of $\tM_k$}
\label{gaugemk}

We  denote by $\calC$ the complex affine
space of holomorphic structures\footnote{To make later
heuristic arguments about infinite dimensional manifolds 
precise we need to choose holomorphic structures of Sobolev class $L^2_1$.} 
on a fixed rank $2$ 
smooth, complex vector bundle $\calV$ of degree $1$. For any
integer $k$ 
consider the
infinite dimensional vector spaces\footnote{To be precise for $(q,r)=(0,1)$ and $(1,0)$ we
consider sections of Sobolev class $L^2_1$ for $(1,1)$ of class $L^2$.}  
$$\Omega^{q,r}_k=\Omega^{q,r}\left(\Sigma,\End(\calV)\otimes \calL_p^k\right),$$
where 
by $\calL_p^k$ we denoted the smooth line bundle
underlying the line bundle $L^k_p$, which is the line bundle of the
divisor $kp$. Now fix $k\geq 0$ and 
define a map \begin{eqnarray}\dbar_k:\calC\times \Omega^{1,0}_k\to
\Omega^{1,1}_k\label{dbark}\end{eqnarray} by sending the pair
$(E,\phi)$ to $\dbar_k^E\phi$, where $$\dbar_k^E:\Omega^{q,r}_k\to \Omega^{q,r+1}_k$$ is 
the $\dbar$ 
operator associated to the holomorphic structure $\End(E)\otimes L^k_p$
on $\End(\calV)\otimes\calL_p^k$. 
It is characterized by the property that for a local section $\phi$
the equation $\dbar_k^E\phi=0$ holds if and
only if $\phi$ is holomorphic with respect to the holomorphic
structure $\End(E)\otimes L_k^p$
on $\End(\calV)\otimes \calL_p^k$. 

Now we define $$\calB_k=\dbar_k^{-1}(0)\subset \calC\times \Omega_k$$
to be the subspace of pairs $(E,\phi)\in \calC\times \Omega_k$ with
$\phi$ being holomorphic. 
We denote by $\pr_k:\calB_k\to \calC$ the projection. Occasionally
$\calB$ will stand for $\calB_0$ and $\pr$ for $\pr_0$.

Now we have the complexified gauge group $\G^c=\Gamma({\rm
Aut}({\calV}))$, the group\footnote{To make later arguments precise we
have to consider gauge transformations of Sobolev class $L^2_2$} 
of complex automorphisms of $\calV$, acting on $\calC$ and on
$\Omega^{1,0}_k$, which induces an action on $\calB_k$. Let us denote
by $(\calB_k)^s\subset \calB_k$ the subspace of Higgs $k$-bundles which are stable. The subset
$(\calB_k)^s\subset \calB_k$ is clearly invariant under the complex
gauge group $\G^c$. 
Then we can form the quotient
$(\calB_k)^0/\G^c$, which is exactly the
moduli space of stable Higgs $k$-bundles $\tM_k$, we are after.

\newpage
\thispagestyle{empty}

\chapter{Geometry of manifolds}

This chapter deals with the symplectic geometry and topology of
manifolds in general. The results appearing here 
will be used in the second part for
the moduli space of Higgs bundles.

\section{$\C^*$-actions on K\"ahler manifolds}
\label{rizsa}

In this section we collect the results from the literature concerning
$\C^*$-actions on K\"ahler manifolds. At the same time we sketch the
structure of Chapter~\ref{compact}.

\subsection{Stratifications}
\label{rizsa1}

Suppose that we are given a K\"ahler manifold $(M,I,\omega)$ with
complex structure $I$ and K\"ahler form $\omega$. Suppose also that
$\C^*$ acts on $M$
biholomorphically with respect to $I$ and such that the K\"ahler form is 
invariant under the induced action of $U(1)\subset \C^*$. Suppose furthermore
that this latter action is Hamiltonian with proper moment map 
$\mu:M\rightarrow \R$, with finitely many critical points and 
$0$ being the absolute minimum of $\mu$.    
Let $\{F_\lambda\}_{\lambda\in A}$ be the set of the components of the fixed 
point set of the $\C^*$-action. 

We list some results of \cite{kirwan1} extended to our case. Namely, Kirwan's
results are stated for compact K\"ahler manifolds, but one can usually modify
the proofs for non-compact manifolds, with proper moment maps as above\footnote{cf. Chapter 9 in 
\cite{kirwan1}.}. 
  
There exist two stratifications in such a situation. The first one
is called the {\em Morse stratification} and can be defined as follows. 
The stratum $U^M_\lambda$, the so-called {\em upward Morse flow from $F_\lambda$}, 
is the set of points of $M$ whose path of 
steepest descent for the Morse function $\mu$ and
the K\"ahler metric have limit points in $F_\lambda$. One can also define 
the sets $D^M_\lambda$, the so-called {\em downward Morse flow of 
$F_\lambda$}, 
as the points of $M$ whose path of steepest descent
for the Morse function $-\mu$ and the K\"ahler metric have limit points in 
$F_\lambda$. $U^M_\lambda$ gives a stratification even in the 
non-compact case, however the set $\bigcup_\lambda
D^M_\lambda$ is not the whole space but a 
deformation retract of it. The set $\bigcup_\lambda D^M_\lambda$ is called the 
{\em downward Morse flow}. 

The other stratification is the {\em Bialynicki-Birula stratification}, where
the stratum $U^B_\lambda$ is the set of points $p\in M$ for which 
$\lim_{t\rightarrow 0}tp\in F_\lambda$. 
Similarly, as above, we can define
$D^B_\lambda$ as the points $p\in M$ for which 
$\lim_{t\rightarrow \infty}tp\in 
F_\lambda$. 

One of Kirwan's important results in \cite{kirwan1} Theorem $6.16$ 
asserts that the 
stratifications $U^M_\lambda$ and $U^B_\lambda$ coincide, and similarly 
$D^M_\lambda=D^B_\lambda=D_\lambda$. This result is important
because it shows that the strata $U_\lambda=U^M_\lambda=U^B_\lambda$ 
of the stratifications are total spaces
of affine bundles (so-called $\beta$-fibrations) on $F_\lambda$
(this follows 
from the Bialynicki-Birula picture) and moreover this stratification
is responsible for the topology of the space $M$ (this follows 
from the Morse 
picture). Thus we have the following theorem\footnote{Cf. 
Theorem $4.1$ of \cite{bialynicki} and also Theorem $1.12$ of 
\cite{thaddeus3}, for the statement about the downward Morse flow
cf. $\S$ 3 of \cite{simpson4}.}: 

\begin{theorem} $U_\lambda$ and $D_\lambda$ are complex submanifolds 
of $M$. They are isomorphic to total spaces of some $\beta$-fibrations
over $F_\lambda$, such that the normal bundles of $F_\lambda$ in these 
$\beta$-fibrations
are $E^+_\lambda$ and $E^-_\lambda$, respectively, 
where $E^+_\lambda$ is the positive and $E^-_\lambda$ is the negative 
subbundle of $D_M\mid_{F_\lambda}$ with respect to the $U(1)$-action. 
Moreover, the downward Morse flow $\bigcup_\lambda D_\lambda$ is a deformation
retract of $M$.
\label{stratification}
\end{theorem}

Recall that a $\beta$-fibration in our case 
is a fibration $E\rightarrow B^n$ with a $\C^*$-action on the 
total space which is locally like $\C^n \times V$, where $V$ is the 
$\C^*$-module $\beta:\C^*\rightarrow GL(V)$. Note that such a fibration is not
a vector bundle in general, but it is if $\beta$ is the sum of isomorphic,
one-dimensional non-trivial $\C^*$-modules.

\subsection{K\"ahler quotients} 
\label{rizsa2}

Whenever we are given a Hamiltonian 
$U(1)$-action on a K\"ahler manifold with a proper moment map, 
we can form the K\"ahler quotients
$Q_t=\mu^{-1}(t)/U(1)$, which are compact K\"ahler orbifolds at a regular value
$t$ of $\mu$. 

If this $U(1)$-action is induced from an action of $\C^*$ on $M$ as  
above, then we can relate the K\"ahler quotients to the 
quotients $M/\C^*$ as follows. First we define 
$M^{min}_t\subset M$ as the set of points in $M$ whose $\C^*$ orbit 
intersects $\mu^{-1}(t)$. Now Theorem 7.4 of \cite{kirwan1} states 
that it is possible to define a complex structure on the
orbit space $M^{min}_t/\C^*$, 
and she also proves that this space is homeomorphic to $Q_t$, 
defining the complex structure for the K\"ahler quotient $Q_t$. (Here again
we used the results of Kirwan for non-compact manifolds, but as above, these 
results can be easily modified for our situation.) It now
simply follows that $M^{min}_t$ only depends on that connected component
of the regular values of $\mu$ in which $t$ lies, and as a consequence of this 
we can see that the complex structure on $Q_t$ is the same as on $Q_{t'}$ if
the interval $[t,t']$ does not contain any critical value of $\mu$. 
We have as a conclusion the following theorem:

\begin{theorem} At a regular level $t \in \R$ of the moment map $\rm \mu$, 
we have the K\"ahler quotient $Q_t=\mu^{-1}(t)/U(1)$ which is a 
compact K\"ahler orbifold with $M^{min}_t$ as a holomorphic $\C^*$-principal 
 orbibundle on it. Moreover $M^{min}_t$ and the complex structure on $Q_t$
only depend on that connected component of the regular values of $\mu$ where
$t$ lies. 
\label{complex}   
\end{theorem}

It follows from the above theorem that there is a discrete family of complex
orbifolds which arise from the above construction. Moreover, at each level
we get a K\"ahler form on the corresponding complex orbifold. 
The evolution of the
different K\"ahler quotients has been well investigated\footnote{E.g. in the papers
\cite{duistermaat-heckman}, \cite{guillemin-sternberg}, cf. also
\cite{thaddeus3} and \cite{brion-procesi}.}. We can summarize these results in the following theorem:

\begin{theorem} The K\"ahler quotients $Q_t$ and $Q_{t'}$ are biholomorphic
if the interval $[t,t']$ does not contain a critical value of the moment
map. They are related by a blowup followed by a blow-down if the interval
$[t,t']$ contains exactly one critical point $c$ different 
from the endpoints. To be more precise, $Q_t$ blown up along the union of
submanifolds 
$\bigcup_{\mu(F_\lambda)=c}P_w(E^-_{\lambda})$ is isomorphic to
$Q_{t'}$ blown up along $\bigcup_{\mu(F_\lambda)=c} P_w(E^+_{\lambda})$ and
in both cases the exceptional divisor is 
$\bigcup_{\mu(F_\lambda)=c} P_w(E^+_{\lambda})\times_{F_\lambda}P_w(E^-_{\lambda})$
the fibre product of weighted projective bundles over $F_\lambda$.

Moreover, in a connected component of the regular values of $\mu$ 
the cohomology classes of the K\"ahler forms $\omega_t(Q_t)$ depend
linearly on $t$ according to the formula:
$$[\omega_t(Q_t)]-[\omega_{t'}(Q_{t'})]=(t-t')c_1(M_t^{min})=
(t-t')c_1(M_{t'}^{min}), $$
 where $c_1$ is the first Chern class of the $\C^*$-principal bundle.
\label{quotients}
\end{theorem}

\subsection{Symplectic cuts}
\label{rizsa3}

Now let us recall the construction of Lerman's symplectic
cut\footnote{Cf. \cite{larman} 
and also \cite{edidin-graham} for the algebraic case.}, 
first in the symplectic and then in the K\"ahler
category. 

If $(M,\omega)$
is a symplectic manifold with a Hamiltonian 
$U(1)$-action and proper
moment map $\mu$ with absolute minimum $0$,
then we can define the symplectic cut of $M$ at the regular level $t$ by a
symplectic quotient construction as follows. 

We let $U(1)$ act on the
symplectic manifold $M\times \C$ (where the symplectic structure
is the product of the symplectic structure on $M$ and the standard
symplectic structure on $\C$) by acting on the first factor according
to the above $U(1)$-action and on the second factor
by the standard multiplication.
This action is clearly Hamiltonian with proper moment map
$\mu+\mu_{\C}$, where
$\mu_\C$ is the standard moment map on $\C$: $\mu_C(z)={\left| z\right|}^2$.

Now if $t$ is a regular value of the moment map $\mu+\mu_\C$, such
that $U(1)$ acts with finite stabilizers on $M_t=\mu^{-1}(t)$ (i.e.
$M_{t}/U(1)$ gives a symplectic orbifold), then the
symplectic quotient $\overline{M}_{\mu<t}$ defined by
$$\overline{M}_{\mu<t}=\{(m,w)\in M\times \C: \mu(m)+\left| w\right|^2=t\}$$
will be a symplectic compactification of the symplectic manifold
$M_{\mu<t}$ in the sense that $$\overline{M}_{\mu<t}=M_{\mu<t}\cup
Q_{t},$$ and the inherited symplectic structure on $\overline{M}_{\mu<t}$
restricted to $M_{\mu<t}$ coincides with its original symplectic structure.
Moreover, if we restrict this structure onto $Q_t$, it coincides with
its quotient symplectic structure.

Now suppose that we are given a K\"ahler manifold $(M,I,\omega)$ and a
holomorphic $\C^*$-action on it, such that the induced $U(1)\subset \C^*$-action
preserves the K\"ahler form  and  is 
Hamiltonian with proper moment map. With 
these extra structures the symplectic cut construction will give us 
$\overline{M}_{\mu<t}$ a compact K\"ahler
orbifold with a $\C^*$-action, such that $\overline{M}_{\mu<t}\setminus Q_t$ is
symplectomorphic to $M_{\mu<t}$ as above and furthermore is biholomorphic 
to $\C^*(M_{\mu<t})$, the union of $\C^*$-orbits intersecting $M_{\mu<t}$.
This is actually an important point\footnote{Cf. \cite{larman}.}, 
as it shows
that $M_{\mu<t}$ is {\em not} K\"ahler isomorphic to
 $\overline{M}_{\mu<t}\setminus Q_t$. 
We can collect all these results into the next theorem:

\begin{theorem} The symplectic cut $\overline{M}_{\mu<t}=M_{\mu<t}\cup Q_t$ 
as a symplectic manifold compactifies the symplectic manifold $M_{\mu<t}$, 
such that the restricted symplectic structure on $Q_t$ coincides with the 
quotient symplectic structure.

Furthermore, if $M$ is a K\"ahler manifold with a $\C^*$-action as above, then
$\overline{M}_{\mu<t}$ will be a K\"ahler orbifold with a $\C^*$-action, such
that $Q_t$ with its quotient complex structure is a codimension $1$ 
complex suborbifold of
$\overline{M}_{\mu<t}$ whose complement is equivariantly biholomorphic to $\C^*(M_{\mu<t})$ with
its canonical $\C^*$-action.
\label{cut}
\end{theorem} 

\begin{remark} Note that if $t$ is higher than the highest critical value
(this assumes that we have finitely many of them),
then $\C^*(M_{\mu<t})=M$ is the whole space, 
therefore the symplectic
cutting in this case gives a holomorphic compactification of $M$
itself. The compactification is 
$\overline{M}_{\mu < t}$, which is equal to the quotient of 
$(M\times {\C} - N \times \{0 \})$ by the action of ${\C}^{\ast}$, where $N$ is the downward Morse flow.
This is the compactification we shall examine in Chapter~\ref{compact} for the case of $\M$, the
moduli space of stable Higgs bundles with fixed determinant of degree $1$.
\end{remark}

\newpage
\section{Cohomologies}

In this section we explain two --not so well known-- cohomology theories,
which will appear in the thesis: Equivariant cohomology of stratified
spaces and at the end hypercohomology of complexes. We have also
inserted some generalities about the Porteous' theorem, which will
enable us to geometrically present equivariant cohomology classes.

\subsection{Equivariant cohomology of stratified spaces}
\label{equivariantcoh}

Let $G$ be a topological group. Let its universal fibration be
$EG\rightarrow BG$. For any $G$-space $X$ we define $X_G:=(X\times
EG)/G$, where $G$ acts on $X\times EG$ with the diagonal action. The
$G$-equivariant cohomology of $X$ is then defined as
$H_G^*(X)=H^*(X_G)$. Since $G$ acts freely on $EG$ we have the
fibration \begin{eqnarray} X\to X_G\to BG. \label{eqfibration1}
\end{eqnarray}
An immediate consequence of this is that if $X$ is
contractible, then the fibration (\ref{eqfibration1}) has contractible
fibres, thus for a contractible space $H^*_G(X)=H^*(BG)$. In
particular $H^*_G(pt)=H^*(BG)$ is not trivial. Since $H^*_G$ is a
contravariant functor the map $X\to {pt}$ gives rise to a ring
homomorphism $H^*(BG)\rightarrow H^*_G(X)$. Thus we see that
$H^*_G(X)$ is a module over $H^*(BG)$. Because of this we see that for
a nontrivial $G$, 
equivariant cohomology is richer than ordinary cohomology.

In case $G$ acts freely on $X$ we have another fibration 
\begin{eqnarray} EG\to X_G\to X/G \label{eqfibration2}, \end{eqnarray}
which shows that for a free action: $$H_G^*(X)\cong H^*(X/G).$$ 

Equivariant cohomology can be calculated particularly well for $G$
stratified spaces. Following $\S 1$ of \cite{atiyah-bott}, we now
explain this.

Let $M$ be a manifold. A
disjoint set $\{M_\lambda\}_{\lambda \in I}$ of locally closed
$G$-invariant
submanifolds of $M$ indexed by a
partially ordered set $I$, with minimal element $0$, defines a {\em
stratification} of $M$   
if $$M=\bigcup_{\lambda\in I}M_\lambda$$ and 
\begin{eqnarray}\overline{M}_\lambda\subset \bigcup_{\mu\geq \lambda}
M_\mu,\label{strat}
\end{eqnarray} moreover we assume that $M_0\neq \emptyset$,
consequently it is the unique open stratum. Because we will encounter
stratifications of Banach manifolds with strata of 
{\em finite codimension}, such that $I$ is countable infinite, we
make two extra finiteness conditions on such stratifications:
\begin{condition}
For every finite subset $A\subset I$ there are a positive, finite number of
minimal elements in the complement $I\setminus A$.\label{condition1}\end{condition} This will ensure
that our inductive arguments still apply. Although the inductions
never terminate for our purposes the following condition will do:
\begin{condition} For each integer $q$ there are only
finitely many indices $\lambda\in I$ such that $\mbox{\rm
codim}(M_\lambda)< q$.  
\end{condition}

Given a $G$-equivariant stratification on $M$, as defined above, 
one can use Morse theory type arguments to get information about the $G$-equivariant
cohomology of $M$. 
To explain this define a subset $J\subset I$ of the indices to be {\em open}
if $\lambda\in J$ and $\mu\leq \lambda$ yields $\mu\in J$. It follows
from (\ref{strat}) that if $J$ is open then $$M_J=\bigcup_{\lambda\in J}M_\lambda$$
is open in $M$.
Now if $J$
is open and $\lambda \in I\setminus J$ is minimal then
$J^\prime=J\cup \lambda$ will be open.  From (\ref{strat}) it follows
that $M_\lambda=M_{J^\prime}\setminus M_J$ is a closed submanifold of 
$M_{J^\prime}$ of 
index, say, $k_\lambda$.
Using the Thom isomorphism we get the long exact sequence of the pair
$(M_{J^\prime},M_J)$ in the form: 
\begin{eqnarray}
\rightarrow H^{q-k_\lambda}_{G}(M_\lambda)\stackrel{(i_\lambda)_*}{\rightarrow} H_G^q(M_{J^\prime}) 
\stackrel{i_J^*} \rightarrow H_G^q(M_J)\rightarrow,\label{longexactofM}\end{eqnarray}
Now we say that the stratification is {\em $G$-perfect} if
\begin{eqnarray}i^*_J:H_G^*(M_{J^\prime})\to H_G^*(M_{J})\label{surjectiveM}\end{eqnarray} is a
surjection for all open $J$ and $\lambda$ minimal in $I\setminus J$.  
In this case it follows that the $G$-equivariant Poincar\'e polynomial of
$M$ \begin{eqnarray} GP_t(M)=\sum_{\lambda\in I} t^{k_\lambda}GP_t(M_\lambda)\label{poincare}\end{eqnarray} 
can be calculated in terms of the $G$-equivariant Poincar\'e
polynomials of the strata.

Thus a $G$-perfect stratification provides a good understanding of the
equivariant cohomology of the space. Consequently useful sufficient
conditions 
for
$G$-perfectness are important. One such criterion is due to
\cite{atiyah-bott}. To explain this consider the long exact sequence
(\ref{longexactofM}) and let $a\in H^*_G(M_\lambda)$. 
Observe that if $a\not= 0$ implies
$i^*_\lambda((i_\lambda)_*(a))\not=0$ then clearly $(i_\lambda)_*$ is
an injection. But  $i^*_\lambda((i_\lambda)_*(a))=ae_\lambda$, where
$e_\lambda$ denotes the equivariant Euler class of the normal bundle of
$M_\lambda$ in $M_{J\cup \lambda}$. Therefore if $e_\lambda$ is not a
zero divisor in $H^*_G(M_\lambda)$ then $(i_\lambda)_*$ is an
injection. Thus a sufficient condition for perfectness is to demand
$e_\lambda\in H^*_G(M_\lambda)$ to be a non zero divisor for each
$\lambda \in I$. We call such
a stratification a {\em strongly $G$-perfect} stratification. 

In the case of a strongly $G$-perfect stratification we have an even better
understanding of the cohomology of $M$. Following the paper \cite{kirwan2} of
Kirwan\footnote{Though \cite{kirwan2} only works with the Shatz
stratification, its methods are general enough to extend them to our
general setting.}, we now explain this. 

Let $J\subset
I$ be an open subset, 
$i_{J}:M_{J}\to M$ be the embedding, let $\lambda\in I\setminus J$
be minimal and $J^\prime=J\cup \lambda$. 
In the case of a strongly $G$-perfect stratification
(\ref{longexactofM}) gives that the restriction map
$$H^*_G(M_{J^\prime})\to H_G^*(M_\lambda)\oplus H_G^*(M_{J})$$ is
injective. An inductive argument then shows that the restriction map
\begin{eqnarray}
H^*_G(M_{J^\prime})\to \bigoplus_{\mu \in
 J^\prime}H^*_G(M_\mu)\label{injective}
\end{eqnarray} is injective.

Now define for any subset $A\subset I$:
$$\calK_A=\ker\left(i^*_{A}:H^*_G(M)\to
H^*_G(M_{A})\right),$$ where $$M_A=\bigcup_{\lambda \in A} M_\lambda$$
and we set $\calK_{\emptyset}=H^*_G(M)$.
Then $\calK_A$ is an ideal of $H^*_G(M)$. Moreover for an open
$J\subset I$ we have $$\calK_J=\bigcap_{\lambda \in J} \calK_\lambda$$
from (\ref{injective}).
If now $J$ and $\lambda$ are as
above then the long exact sequence (\ref{longexactofM}) gives that
$$i_\lambda(\calK_J)\subset \langle e_\lambda \rangle\subset
H^*_G(M_\lambda),$$ where $\langle e_\lambda \rangle$ is the ideal of $H^*_G(M_\lambda)$
generated by $e_\lambda$. Since the 
stratification is $G$-perfect the map $i^*_{J^\prime}:H^*_G(M)\to
H^*_G(M_{J^\prime})$ is surjective consequently \begin{eqnarray}
i^*_{\lambda}(\calK_J)=\langle e_\lambda \rangle\subset
H^*_G(M_\lambda). 
\label{condition} \end{eqnarray}
The following proposition shows that, in the case of a strongly 
$G$-perfect stratification, in some sense $\calK_J$ is
unique with respect to this property.

\begin{proposition}[Kirwan]
Let $M=\bigcup_{\lambda \in I} M_\lambda$ be a strongly $G$-perfect
stratification of $M$. 
Suppose we are given a subset $\calR_{\lambda} \subset H^*_G(M)$ for each
$\lambda \in I$, with the following property: \begin{eqnarray}
i_\mu^*(\calR_\lambda)=0 \mbox{ if } \mu\not\ge \lambda, \mbox{ and } i^*_\lambda(\calR_\lambda)=
\langle e_\lambda\rangle \subset
H^*_G(M_\lambda),\label{kirwancondition} \end{eqnarray}
where $\langle e_\lambda \rangle$ is the ideal of $H^*_G(M_\lambda)$
generated by the equivariant Euler class of the normal
bundle\footnote{We set $e_0=1$.} 
of $M_\lambda$ in $M$.
Then for any open subset $J\subset I$ we have 
$$\langle \calR_{\bar{J}}\rangle_\Qu=\calK_J=\ker\left( i^*_{J}:H^*_G(M)\rightarrow
H^*_G(M_{J})\right),$$ where $\langle \calR_{\bar{J}}
\rangle_\Qu$ denotes the $\Qu$ vector subspace of
$H^*_G(M)$, generated additively by $\calR_{\bar{J}}=
\bigcup_{\mu\not\in J} \calR_{\mu}$.    
\label{propositionkirwan}
\end{proposition}
\begin{proof} Though the proof is essentially the same as the
proof\footnote{Cf. also Proposition 11 of \cite{earl}.} 
of Proposition 1 of \cite{kirwan1}, we give it here for
the sake of completeness. 

Let us 
suppose that $a\in H^*_G(M)$ such that $i^*_{J}(a)=0$, i.e. $a\in \calK_J$. 
Let $\lambda$ be a
smallest element of $I\setminus J$ and set $J^\prime=J\cup \lambda$. 
By the assumption (\ref{kirwancondition}) we have $a_\lambda \in \calR_{\lambda}$ such
that $i^*_{\lambda}(a_\lambda)=i^*_{\lambda}(a)$.   
It follows that $i^*_\mu(a-a_\lambda)=0$ for every
$\mu\in {J^\prime}$, hence $i_{J^\prime}^*(a-a_\lambda)=0$ from (\ref{injective}). 
An inductive argument\footnote{That this works even in the infinite
dimensional case is ensured by Condition 1 and 2 above.} now gives $a_\mu\in R_\mu$ for each $\mu\not\in
J$ such that $a=\sum_{\mu\not\in J}a_\mu$.

The result follows. 
\end{proof}

\paragraph{Remarks.} 1. If $J=\{0\}$ then the above result gives that
$$\left\langle \bigcup_{\mu\not= 0}  \calR_\mu\right\rangle_\Qu
=\calK_{\{0\}}=
\ker\left(i^*_0:H^*_G(M)\to
H^*_G(M_0)\right),$$ thus if one can choose $\calR$ in a simple
form, satisfying the condition (\ref{kirwancondition}), 
then one has information about the
relations 
in $H^*_G(M_0)$.  This form was actually stated and used by
Kirwan in \cite{kirwan2} in the proof of the 
Mumford conjecture\footnote{Cf. Section~\ref{cohomologyofN}.}. In
Section~\ref{mumfordproof} we provide a purely geometric proof of
the Mumford conjecture using this special case.

2. We will also use another special case of the above proposition in
Chapter~\ref{cohomology}. Namely if $J=\emptyset$, then the
proposition says that $$\left\langle\bigcup_{\mu\in I}
\calR_\mu\right\rangle_\Qu=\calK_\emptyset=H^*_G(M).$$ Thus if one can find the
sets $\calR_\mu$ being generated by a subset of $H^*_G(M)$, then the
whole $H^*_G(M)$ is generated by this subset. We now give another,
closely related, application of this special case of
Proposition~\ref{propositionkirwan}:

\begin{corollary}  Let $M=\bigcup_{\lambda\in I} M_\lambda$ and 
$M^\prime =\bigcup_{\lambda\in I} M^\prime_\lambda$ be two
strongly $G$-perfect stratified spaces. 
Suppose further that a map $f:M^\prime\to M$ is such that
$f^{-1}(M_\lambda)=M^\prime_\lambda$ and 
$$f^*_\lambda: H_G^*(M_\lambda)\to H_G^*(M^\prime_\lambda)$$ is surjective
for each $\lambda \in I$. Then $$f^*:H_G^*(M)\to
H_G^*(M^\prime)$$ is surjective.
\label{surjection}
\end{corollary}

\begin{proof}  For each $\lambda\in I$ set $\calR^\prime_\lambda=f^*(\calK_\lambda)\subset
H^*_G(M^\prime)$, which is an additively closed subspace of
$H^*(M^\prime)$. From $f_\lambda^*(e_\lambda)=e^\prime_\lambda$ and
(\ref{condition}) it follows that the sets
$\left\{\calR^\prime_\lambda\right\}_{\lambda\in I}$ satisfy the
conditions of Proposition~\ref{propositionkirwan}. In particular
$$f^*(H^*_G(M))=\langle \calR^\prime_I\rangle_\Qu=\calK^\prime_\emptyset=H^*_G(M^\prime).$$
The result follows. 
\end{proof}

Many examples of $G$-perfect stratifications arise in symplectic
geometry. We explain one special case in detail:

\subsubsection{The special case $G=U(1)$.}  
First note that $BU(1)\sim \C P^\infty$ thus $H^*(BU(1))\cong \Qu [u]$
is a
free polynomial algebra on a degree $2$ generator $u$. For convenience
we will write $H^*_\circ$ instead of $H^*_{U(1)}$. 

Assume that $M$ is a symplectic manifold and $U(1)$ acts on $M$
in a Hamiltonian way with proper moment map as in
Subsection~\ref{rizsa1}. Recall from Subsection~\ref{rizsa1} that in
this case we get a stratification $M=\bigcup_{\lambda\in I} U_\lambda$,
which is a stratification in the sense
we defined above, as Kirwan proves in 
\cite{kirwan1}. Moreover an important result of Kirwan in
\cite{kirwan1} shows that this stratification is always strongly 
$U(1)$-perfect. 
Thus we can calculate the $U(1)$-equivariant cohomology of $M$ from
(\ref{poincare}), and even say something about the ring $H^*_\circ(M)$ from 
Proposition~\ref{propositionkirwan}. However in the present case we
have an even 
better understanding of the ring $H^*_{\circ}(M)$ through the so called {\em
Localization Theorem}: 

\begin{theorem}[Localization Theorem] For any $\psi\in H^*_\circ(M)$ one
gets the following equation\footnote{Recall from
Subsection~\ref{rizsa1} that $\{F_\lambda\}_{\lambda \in I}$ denotes the
fixed point set of the circle action.} in the localized ring 
$H^*_\circ(M)\otimes_{\Qu [u]} \Qu (u)$:
$$\psi=\sum_{\lambda\in I} \frac{(i_{F_\lambda})_* i^*_{F_\lambda}\psi}{E_\circ(\nu_{F_\lambda})},$$
where $E_\circ(\nu_{F_\lambda})\in H^*_\circ(F_\lambda)\cong
H^*(F_\lambda)\otimes \Qu [u]$ 
is the equivariant Euler class of the normal bundle of $F_\lambda$ in $M$.
\label{localization}
\end{theorem}

An easy consequence of this theorem is the following\footnote{Note
that it
follows also from (\ref{injective}).}  
\begin{corollary}[Kirwan] The restriction map
\begin{eqnarray} H^*_\circ(M)\to\bigoplus_{\lambda \in I}H^*_\circ(F_\lambda)\label{injection}\end{eqnarray}
is an injection. 
\end{corollary}

Thus in order to understand $H^*_\circ(M)$ it is enough to
understand the restriction maps $H^*_\circ(M)\rightarrow
H^*_\circ(F_\lambda)$ and the ring structures\footnote{Since $U(1)$
acts trivially on $F_\lambda$, we have $H^*_\circ(F_\lambda)\cong
H^*(F_\lambda)\otimes H^*(BU(1))$.} of $H^*(F_\lambda)$. This will be
our method of describing the cohomology ring of the moduli space of
Higgs bundles in Chapter~\ref{cohomology}.

\subsection{Porteous' theorem}
\label{porteoussection}

First we introduce some notation concerning cohomology classes of
subvarieties, which will be used throughout this thesis:

\begin{notation} If $X$ is an irreducible locally closed
subvariety of a smooth algebraic variety $Y$ of
codimension d, 
then $\eta^Y_X\in H^{2d}(Y)$ denotes the cohomology class of $\overline{X}$ in
$Y$. Moreover if $G$ acts on $Y$ and $X$ is $G$-invariant, then
$\eta^{G,Y}_X\in H^{2d}_G(Y)$ denotes the $G$-equivariant cohomology
class of $\overline{X}$. 

If $X$ is an irreducible locally closed and relatively 
complete\footnote{It means that $\overline{X}$ is complete i.e. compact.} subvariety of $Y$ then
$\overline{\eta}^Y_X\in H_{cpt}^{2d}(Y)$ denotes the compactly supported
cohomology class of $\overline{X}$ in $Y$. 
\label{cohclass}
\end{notation}

A frequently used technique of the present thesis will be to find some
geometric way to calculate the cohomology class $\eta_X^Y$. The
prototype of such a presentation is when a closed subvariety $X$ is 
the zero locus of a section
$s$ of a vector bundle $W$ of rank $d={\rm codim} X$. 
In differential topology if the section $s$ is transversal, then we have the
geometric formula
\begin{eqnarray}c_d(W)=\eta^X_Y\label{geometric}.\end{eqnarray} 
In algebraic geometry the transversality
assumption is replaced by demanding $X$ to have codimension
$d=\rank(W)$ in order the geometric formula (\ref{geometric}) to be
true. In this thesis we will need some ramification of this idea. 

To explain it let us reformulate the above construction, in a manner more suitable for
this thesis, as follows: Let $V$ be the trivial line bundle on
$Y$. Then $s$ is equivalent to a map $s:V\to W$. Let us denote by
$\mathbf F$ the virtual bundle $W-V\in K(Y)$. Then we get that
$\mathbf F$ outside
$X$ is an honest vector bundle.  

More generally suppose that we are given a virtual bundle $\mathbf F$ in $K(Y)$ of rank\footnote{I.e. $0$th
Chern class.} $d-1>0$. Suppose also that ${\mathbf F}\mid_{Y\setminus X}$ is an honest vector
bundle outside a closed smooth subvariety $X\subset Y$ of complex codimension
$d$. Then consider the following bit of the cohomology long exact sequence of the pair
$(Y,Y\setminus X)$: \begin{eqnarray}
H^0(X)\stackrel{\tau}{\to} H^{2d}(Y)\to H^{2d}(Y\setminus
X),\label{porteousexact}\end{eqnarray} 
where $\tau$ is the Thom  
map. Since ${\mathbf F}\mid_{Y\setminus
X}$ is an honest vector bundle of rank $d-1$, we have that
$c_d({\mathbf F})\mid_{Y\setminus X}=0$ vanishes. From the exactness of
(\ref{porteousexact}), we have that $c_d({\mathbf F})=\tau(q)=q\cdot \eta_X^Y$
for some $q\in H^0(X)\cong \Qu$. If $q$ was non-zero, or especially
$1$, we would find a nice geometric way expressing $\eta_X^Y$ as
$c_d({\mathbf F})$. The following special case of 
Porteous' theorem\footnote{Cf. (4.2) of \cite{arbarello-et-al}.} states
that if ${\mathbf F}=W-V$ is a difference of two vector bundles 
then $q=1$, 
choosing $X$ to be the degeneration locus of a homomorphism $f:V\to W$.

\begin{theorem}[Porteous] Let $f:V\to W$ be a homomorphism of vector
bundles over a smooth algebraic variety $Y$. Let $X$ be the
degeneration locus of $f$, i.e. the subvariety of $Y$ consisting 
of points, where $f$ fails
to be an injection\footnote{For a rigorous construction of degeneracy
loci cf. \cite{arbarello-et-al} p.83. Our degeneracy locus is the
$k$-th degeneracy locus of \cite{arbarello-et-al}, where
$k=\rank(V)-1$.} Then we have $$\eta_X^Y=c_d(W-V),$$
if the codimension of $X$ coincides with
$d=\rank(W)-\rank(V)+1$.
\label{porteous}
\end{theorem}

Since (\ref{porteousexact}) exists in the equivariant cohomology too,
the next corollary is an easy consequence of this Theorem~\ref{porteous}.

\begin{corollary} Let  $G$ act on $Y$, and $V$ and $W$ be
$G$-equivariant vector bundles, with  $f$
a $G$-equivariant homomorphism. We have  
$$\eta_X^{G,Y}=c_d^G(W-V),$$ if $X$, as in the above
Theorem~\ref{porteous}, has the expected codimension
$d=\rank(W)-\rank(V)+1$.
\label{eqporteous} 
\end{corollary}

\newpage
 \subsection{Hypercohomology of complexes}
\label{hypercohomology}

In this subsection we recall the notion of hypercohomology of a complex from 
\cite{griffiths-harris}, and list some properties of it, which we will
use in Chapter~\ref{intersection} to describe rigorously the virtual Dirac
bundle, already mentioned in Subsection~\ref{diracequ}. 

\begin{definition} Let $${\cal A}=(A_0\stackrel{d}{\longrightarrow}A_1
\stackrel{d}{\longrightarrow}A_2\longrightarrow...)$$
be a complex of coherent sheaves $A_i$ over an algebraic variety $X$. For a
covering $\underline{U}=\{U_\lambda\}$ of $X$ and each
$A_i$ we get the \v Cech cochain complex with boundary operator $\delta$:
$$(C^0(\underline{U},A_i)\stackrel{\delta}{\longrightarrow} 
C^1(\underline{U},A_i)\stackrel{\delta}{\longrightarrow}...).$$
Clearly $d$ induces operators 
$$(C^j(\underline{U},A_i)\stackrel{d}{\longrightarrow} 
C^j(\underline{U},A_j)),$$
satisfying $\delta^2=d^2=d\delta+\delta d=0$: and hence gives rise
to a double complex $$\{C^{p,q}=C^p(\underline{U},A_q);\delta,d\}.$$
The hypercohomology of the complex $\cal A$ is given by the cohomology
of the total complex of the double complex $C^{p,q}$:
$$\Hy^*(X,{\cal A})=\lim_{\underline{U}}H^*(C^*(\underline{U}),D).$$

Moreover if $\cal A$ is a complex over $X$ and $f:X\rightarrow Y$ is 
a projective morphism then for every non-negative integer $i$ 
define the sheaf $\R^i f_*({\cal A})$ over $Y$
by $$\R^i f_*({\cal A})(U)=\Hy^i(f^{-1}(U),{\cal A}).$$ 

Finally, define the pushforward of a complex to be:

$$f_!({\cal A})=
\R^0 f_*({\cal A})-\R^1 f_*({\cal A})+\R^2 f_*({\cal A})-\dots\in K(Y).$$
\end{definition}

\begin{remark} 
In the present thesis we will work only with two-term complexes. 

There is one important property of hypercohomology which we will
make constant use of. If $$0\rightarrow {\cal A}\rightarrow
{\cal B} \rightarrow {\cal C}\rightarrow 0$$ 
is a short exact sequence of complexes
then there is a long exact sequence of hypercohomology vector spaces:
\begin{eqnarray} 0\rightarrow \Hy^0(X,{\cal A})
\rightarrow \Hy^0(X,{\cal B})\rightarrow \Hy^0(X,{\cal C})
\rightarrow \Hy^1(X,{\cal A}) \rightarrow 
\dots 
\end{eqnarray}
As an example consider the short exact sequence of two-term
complexes: 
$$
\begin{array}{c}0\longrightarrow 0\\
		\uparrow \hskip1cm \uparrow \\
		 0\stackrel{0}{\longrightarrow} A_2 \\
		 \hskip.36cm \uparrow  \hskip1cm \uparrow 
			{\mbox{\scriptsize{$\cong$}}}\\
		A_1\stackrel{d}{\longrightarrow} A_2\\
		\mbox{\scriptsize{$\cong$}}\uparrow\hskip1cm\uparrow 
					\hskip.36cm \\
		A_1\stackrel{0}{\longrightarrow} 0\\
		\uparrow\hskip1cm\uparrow\\
		0\longrightarrow 0
\end{array}		
$$
The long exact sequence in this case is:
\begin{eqnarray}0\rightarrow \Hy^0(X,{\cal A})
\rightarrow H^0(X,A_1)\rightarrow H^0(X,A_2)\rightarrow 
\Hy^1(X,{\cal A}) \rightarrow \dots \label{longexact}
\end{eqnarray}
which we will call the hypercohomology long exact sequence of the
two-term 
complex ${\cal A}=A_1\stackrel{d}\rightarrow A_2$.

Consequently if ${\cal A}=A_1\stackrel{d}{\rightarrow} A_2$
is a two-term complex over $X$ and $f:X\rightarrow Y$ is a projective 
morphism
then we have:

$$0\rightarrow  \R^0f_*(X,{\cal A})
\rightarrow R^0f_*(X,A_1)\rightarrow R^0f_*(X,A_2)\rightarrow 
\R^1f_*(X,{\cal A}) \rightarrow
\dots, $$
a long exact sequence of sheaves over $Y$.
\end{remark}

\newpage
\thispagestyle{empty}

\chapter{Cohomology of moduli spaces}
\label{cohomologyofmoduli}

In general the
cohomology of moduli spaces obviously plays an  important role in the
algebraic geometric understanding of the moduli problem. As was
recently discovered by
physicists\footnote{Cf. Section~\ref{motivation}}
the cohomology
of the moduli spaces of some physical theories contain physically
relevant information about the global aspects of the theory. For
an example intersection numbers on these moduli spaces sometimes
can be identified with correlation functions, which are the experimentally
measurable quantities of the theory.   

The main subject of this thesis is the determination of the rational
cohomology of $\M$, the moduli space of Higgs bundles. Since we will
use the cohomological properties of other moduli spaces, and we would like to
convey a general view of the cohomology of moduli spaces of
objects on an algebraic curve $\Sigma$, 
we describe in detail the cohomology of some such moduli spaces here.

First we summarize the general picture which is emerging from the rest
of the section. Let $M$ be a smooth manifold (not necessarily compact),
the moduli space of a moduli problem. In each of the following cases
we will find an infinite tower of spaces $M_k$, (which usually will be
the moduli spaces of some deformed moduli problem):     
 $$M=M_0\subset M_1 \subset
\dots \subset M_k \subset \dots,$$ with inclusions $i_k:M_k\rightarrow
M_{k+1}$, such that $i_k^*:H^*(M_{k+1})\rightarrow H^*(M_k)$ is surjective. From this data
we have the direct limit: $$M_\infty=\lim_{\longrightarrow} M_k,$$ 
and the inverse limit: $$H^*(M_\infty)=\lim_{\longleftarrow} H^*(M_k),$$
since $H^*$ is a contravariant functor. It will turn out that
$H^*(M_\infty)$ will be a free graded commutative\footnote{It means that
even degree classes commute with any other class, while odd degree
classes anticommute among each other. Note that the cohomology ring of
any space is automatically graded commutative.} algebra on a finite set of  universal generators. 
Furthermore, it follows that the map $$i^*_M:H^*(M_\infty)\rightarrow H^*(M)$$ is
surjective. Putting everything together, the cohomology ring $H^*(M)$ is
generated 
by the images of the 
universal generators and the kernel of $i^*_M$ provides the
relations, in other words $i^*_M$ is a resolution of the cohomology
ring of $M$. Because of this we call a tower having the above properties a
{\em resolution tower} of $M$.  

In Chapter~\ref{minfty} we will show that the cohomology ring of $\M$,
the moduli space of Higgs bundles, can be understood in this general
framework too.

\newpage
\section{The curve $\Sigma$}
\label{sigma}

The basic object of this thesis is a  
fixed, smooth and complex projective curve $\Sigma$ of genus $g\geq 2$. We
also fix a point $p\in \Sigma$.

An additive basis of $H^*(\Sigma)$: $1\in H^0(\Sigma)$, 
$\xisi_i\in H^1(\Sigma),\  i=1,..,2g$ and the fundamental cohomology class 
$\sisi\in H^2(\Sigma)$ with the properties that 
$\xisi_i\cdot \xisi_{i+g} = -\xisi_{i+g}\cdot \xisi_i =\sisi$ for $i=1,..,g$ and otherwise
$\xisi_i\cdot \xisi_j=0$, will be fixed throughout the thesis.  

As a matter of fact the curve $\Sigma$ is itself a moduli space, the
moduli space of its points, or in a more sophisticated way: the moduli space
of degree one divisors on $\Sigma$, or in other words $\Sigma_1$
the first symmetric product.

\newpage
\section{The Jacobian $\J$}
\label{cohomologyofJ}

The moduli space of 
line bundles of degree $k$ over $\Sigma$ is the Jacobian $\J_k$. This is an 
Abelian variety of dimension $g$. Tensoring by a fixed line bundle of degree
$k-l$ gives an isomorphism between $\J_l$ and $\J_k$.
We will write $\J$ for $\J_1$. 

Being a torus 
$H^*(\J_k)$ is a free exterior algebra --in particular a free graded
commutative algebra-- on $2g$ classes $\tau_i\in H^1(\J_k)$
defined by the formula $$c_1(\bP_k)=k\otimes\sisi + 
\sum_{i=1}^{g}(\tau_i
\otimes \xisi_{i+g}-\tau_{i+g}\otimes \xisi_i)\in H^2(\J_k\times \Sigma)\cong
\sum_{r=0}^{2} H^r(\J_k)\otimes H^{2-r}(\Sigma).$$ Here $\bP_k$ is
the normalized Poincar\'e bundle, or universal line bundle over 
$\J_k\times \Sigma$. Universal means that for any $L\in \J_k$: 
$$\bP_k\mid_{\{L\}\times \Sigma}\cong L$$
and normalized means that $\bP_k\mid_{\J_k\times \{p\}}$ is
trivial\footnote{Cf. p.166 \cite{arbarello-et-al}.}.

Being a free exterior algebra, the Poincar\'e polynomial of $H^*(\J)$
is given by $$P_t(\J_k)=(1+t)^{2g}.$$

\newpage
\section{The symmetric product $\Sigma_n$}
\label{symmetric}

The $n$-th symmetric product $\Sigma_n$ is the moduli space of degree
$n$ effective divisors. It is a smooth projective variety of dimension $n$. 
Clearly $\Sigma_1=\Sigma$.

The cohomology ring $H^*(\Sigma_n)$ 
is multiplicatively generated by $\xi_i\in H^1(\Sigma_n),$ for
$i=1\dots 2g$ and
$\eta\in H^2(\Sigma_n)$ defined by the formula:      
\begin{eqnarray}c_1(\Delta_k)=k\otimes \sisi + \sum_{i=1}^{g}(\xi_i
\otimes \xisi_{i+g}-\xi_{i+g}\otimes \xisi_i)  + \eta\otimes 1
\label{cherndelta}\end{eqnarray}
in the decomposition $$ H^2(\Sigma_n\times \Sigma)\cong
\sum_{r=0}^{2} H^r(\Sigma_n)\otimes H^{2-r}(\Sigma).$$

Here $\Delta_n\in Div(\Sigma_n\times \Sigma)$ is the universal
divisor\footnote{Cf. \cite{arbarello-et-al}.},
i.e. $\Delta_n\mid_{\{D\}\times \Sigma}=D$ for every divisor 
$D\in \Sigma_n$.
The relation set of the ring $H^*(\Sigma_n)$ and thus
a complete description of it is given in \cite{macdonald}. We let
$\sigma_i=\xi_i\xi_{i+g}$ and $\sigma=\sum_{i=1}^g\sigma_i \in H^2(\Sigma_n)$.
  
The Poincar\'e polynomial of $H^*(\Sigma_n)$ was calculated in
\cite{macdonald} in the form: \begin{eqnarray}P_t(\Sigma_n)=
\Coeff_{x^n}\left(\frac{(1+xt)^{2g}}{(1-x)(1-xt^2)}\right).\label{poincaresigman}\end{eqnarray}
For the case $n>2g-2$ we have the Abel-Jacobi map $\Sigma_n\to \J_n$
being a locally trivial fibration with fibre $\Proj^{n-g}$, which for
the Poincar\'e polynomial gives $$P_t(\Sigma_n)=\frac{(1+t)^{2g}(1-t^{2(n-g+1)})}{(1-t^2)}.$$

\subsection{The resolution tower of $\Sigma$}

There are embeddings $i_n:\Sigma_n\rightarrow \Sigma_{n+1}$
given by $i_n(D)=D+p$, yielding the tower \begin{eqnarray}
\Sigma_1\subset \Sigma_2 \subset \dots \subset \Sigma_n \subset
\dots. \label{resolutionsymmetric}
\end{eqnarray}  We
consider the direct limit of them: $$\Sigma_\infty=\lim_{n\to
\infty}\Sigma_n.$$ It is a $\Proj^\infty$ bundle over $\J$,
thus its Poincar\'e polynomial is 
\begin{eqnarray}P_t(\Sigma^\infty)=\frac
{(1+t)^{2g}}{(1-t^2)}.\label{poincaresigmainfty}
\end{eqnarray}
The pullback map \begin{eqnarray} i_n^*:H^*(\Sigma_{n+1})\rightarrow
H^*(\Sigma_n)\label{surjectsymmetric} \end{eqnarray}
is surjective because these rings are generated by universal classes.
Since the cohomology ring of $\Sigma_\infty$ is the inverse limit
of the cohomology rings of the $\Sigma_n$'s, i.e. 
$$H^*(\Sigma_\infty)=\lim_{\infty \leftarrow n} H^*(\Sigma_n),$$
it is also generated by the tautological classes $\xi_i$ and
$\eta$. In fact it is a free graded commutative
algebra on these
generators, as (\ref{poincaresigmainfty}) shows.

\subsection{Some results about $H^*(\Sigma_n)$}
It is convenient to insert the following lemmata here, which will be
needed in Section~\ref{relations}. They are taken from \cite{hausel-thaddeus}.

\begin{lemma} Let us denote by $H_I^*(\Sigma_n)$ the subring of
$H^*(\Sigma)$ generated by $\eta$ and $\sigma$. For $n\leq 2g-2$  
its Poincar\'e polynomial $P_T^I(\Sigma_n)$ can be written in the form:
\begin{eqnarray}
P_T^I(\Sigma_n)=\sum_\stack{q,s\geq 0}{q+2s\leq n}
T^{q+s}\label{invpoinsigman} \end{eqnarray}
\end{lemma}
\begin{proof} It follows from B-1 and B-2 on p. 328-329 of
\cite{arbarello-et-al} that the classes $\eta^q\cdot \sigma^{s}$ for $q,s\geq 0$
and $q+2s\leq n$ form an additive basis for $H^*_I(\Sigma_n)$. The
result follows.
\end{proof}

\begin{lemma} Let $n,k,m$ and $l$ be non-negative integers. If
$n-g+m\leq l$ and $g+k-m<l$ then $$\left( \frac{\exp(\si) \eta^{k}}{(1+\eta)^{m}} 
\right)_{l} = 0$$ over $\Sigma_n$. Here $(\cdots)_l$ denotes the part of total
degree $2l$.  
\label{symmetricvanish}
\end{lemma}

\begin{proof} First we show that Poincar\'e duality still holds in the
subring generated by $\eta$ and $\si$.  
By Poincar\'e
duality in $H^*(\Sigma_n)$, the vanishing of a polynomial in $\eta$ and
$\si$ of total degree $2d$ is equivalent
to the vanishing of its product with every monomial in $\eta \in H^2$ and
$\xi_i \in H^1(\Sigma_n)$ of total degree $2(n-d)$.  If the power of $\xi_i$
is greater than $1$, this certainly vanishes; if the power of $\xi_i$ is
$1$ but the power of $\xi_{i \pm g}$ is $0$, this again vanishes since 
$H^*(\Sigma_n)\cong H^*(\Sigma^n)^{S_n}\subset H^*(\Sigma^n)$ and the corresponding
expression 
in $H^*(\Sigma^n)$ is certainly zero. So it suffices to
consider monomials in $\eta$ and $\si_i$ with degree $1$ in the latter;
then by symmetry it suffices to consider monomials in $\eta$ and $\si$
only.

As pointed out by Zagier in (7.2) of \cite{thaddeus2}, for any power series
$A(x)$ and $B(x)$, 
$$A(\eta) \exp(B(\eta)\si)[\Sigma_n] 
= \Res_{\eta=0} \left( \frac{A(\eta)(1+\eta B(\eta))^g}{\eta^{n+1}}
\right).$$
We multiply our expression by the generating function 
$\exp(s\si)/(1+t\eta)$ for the monomials in $\eta$ and $\si$, and ask
the
coefficient of $s^i t^j$ to be 0 whenever $i + j =n-l$:
\begin{eqnarray*}
\Coeff_{s^i t^j} \left(\frac{\exp((s+1)\si)
\eta^{k}}{(1+\eta)^{m}(1+t\eta)}[\Sigma_n]\right) 
& = &  
\Coeff_{s^i t^j} \Res_{\eta=0} 
\left(
\frac{\eta^{k}((\eta+1)+s\eta)^g}{(1+\eta)^{m}(1+t\eta)\eta^{n+1}}
\right) \\
& = &  \const \Res_{\eta=0} 
\left( \frac{\eta^{k}(1+\eta)^{g-i}\eta^i\eta^j}{(1+\eta)^{m}\eta^{n+1}}
\right) \\
& = &  \const \Res_{\eta=0} 
\left( \eta^{k+i+j-n-1}(1+\eta)^{g-i-m}
\right).
\end{eqnarray*}
But this is $0$ because $$g-i-m\geq g-m-(n-l)\geq 0,$$ from the first condition, 
thus $(1+\eta)^{g-i-m}$ is a polynomial of degree $g-i-m$, therefore 
the highest 
degree term of $\eta^{k+i+j-n-1}(1+\eta)^{g-i-m}$ has degree
$$(k+i+j-n-1)+(g-i-m)\leq k+(n-l)-n-1+g-m=k-l-1+g-m<-1$$ from the
second condition.  

The result follows. 
\end{proof}

\newpage
\section{The moduli space of rank $2$ stable bundles $\N$}
\label{cohomologyofN}

Collecting results from the literature, 
in this section we describe $H^*(\N)$.
First we explain how Atiyah and Bott calculated the Poincar\'e
polynomial of $\N$ and how they found generators for the ring. Then we
describe the method of Kirwan for proving Mumford's conjecture and
thus providing in principle a complete description of the cohomology ring
$H^*(\N)$. Then we state a calculationally useful 
description of the ring structure of $H^*(\N)$ and cite some formulae of Zagier.

\subsection{The Poincar\'e polynomial and generators}

Here we explain how the idea of Subsection~\ref{yangmills2} was made
rigorous in \cite{atiyah-bott} in the rank $2$ case. 
\subsubsection{The Shatz stratification}
Let $\calV$ be a fixed rank $2$ smooth complex vector bundle of
degree $1$. Let $\calC$ denote the infinite dimensional complex affine
space of holomorphic structures on $\calV$.  Fixing a hermitian structure on $\calV$, we have the
gauge group $\G$ of smooth unitary automorphisms of $\calV$ acting
naturally on $\calC$.
Its complexification  $\G^c=\Gamma(\Sigma, {\rm Aut}(\calV))$, the
group of 
complex automorphisms of $\calV$, also acts on
$\calC$. Moreover $\G^c/\G$ is the contractible space of Hermitian
structures on $\calV$, thus $B\G\sim B\G^c$ and 
so for the purpose of equivariant cohomology
they are equivalent: $H^*_{\G^c}=H^*_{\G}$. For convenience 
we will always use $\G$-equivariant cohomology, even where 
$\G^c$-equivariant cohomology is understood.  

We also let 
$\calC_s\subset\calC$ (and
frequently $\calC_0$) denote the open subspace of stable bundles. By
Theorem~\ref{harder-narasimhan} of Harder and Narasimhan, 
if $E\in \calC$ is not stable it has a unique destabilizing subbundle
 of degree $d>0$. For $d>0$ we let
$\calC_d\subset\calC$ denote the subspace of unstable 
vector bundles with destabilizing bundle of degree $d$. Since the
Harder-Narashiman filtration is canonical, $\calC_d$ is invariant under
$\G^c$, thus each $\calC_d$ is a union of orbits. Atiyah and Bott
prove that for $d>0$, the space $\calC_d$ is locally a
Banach submanifold of $\calC$ of
finite codimension $2g+4d-4$ and that \begin{eqnarray}
\calC=\bigcup_{d=0}^\infty \calC_d \label{atiyahbottstrat}\end{eqnarray}
 is a $\G^c$-equivariant 
stratification in the sense of Subsection~\ref{equivariantcoh}.
It is called the {\em Shatz stratification}.   

The paper \cite{atiyah-bott} uses Morse theory, in a manner as we
explained in Subsection~\ref{equivariantcoh}, for the Shatz
stratification in order to calculate the Poincar\'e polynomial of 
$\tN=\calC_0/\G^c$. They prove that the stratification is strongly
$\G^c$-perfect\footnote{Cf. Subsection~\ref{equivariantcoh}.}. It
follows from (\ref{poincare}) that for 
the $\G$-equivariant Poincar\'e polynomials we have: 
$$\G P_t(\calC)=\G P_t(\calC_0)+\sum_{i=1}^\infty t^{2g+4i-4}\G P_t(\calC_i),$$
or in a more suitable form
\begin{eqnarray}
\G P_t(\calC_0)=\G P_t(\calC)-\sum_{i=1}^\infty t^{2g+4i-4}\G P_t(\calC_i).\label{poincarec}
\end{eqnarray}

Now $\calC$ is contractible and so $H^*_{\G}(\calC)\cong
H^*(B\G)$. The cohomology ring of $B\G$
is described in $\S 2$ of \cite{atiyah-bott}:

\subsubsection{The cohomology of $B\G$}

The ring $H^*(B\G)$ is freely
generated as a graded commutative algebra 
by classes \begin{eqnarray*} a_r\in H^{2r}(B\G), \hspace{1cm} 
b_r^j\in H^{2r-1}(B\G),\hspace{1cm} f_2\in H^2(B\G), 
\end{eqnarray*} for $1\leq r \leq 2$ and
$1\leq j\leq 2g$. 
These classes appear as the K\"unneth components of a certain
(universal)
rank $2$
vector bundle $\U$ on $B\G\times \Sigma$. Namely $$c_1(\U)=a_1\otimes
1+ \sum_{j=1}^{2g} b_1^j\otimes \xisi_j,$$
and $$c_2(\U)=a_2\otimes
1+ \sum_{j=1}^{2g} b_2^j\otimes \xisi_j+f_2\otimes \sigma.$$

Since $H^*(B\G)$ is a freely generated graded commutative algebra on
the above mentioned classes,
its Poincar\'e polynomial
is \begin{eqnarray}
P_t(B\G)=\frac{\left\{(1+t)(1+t^3)\right\}^{2g}}{(1-t^2)^2(1-t^4)}.\label{poincarebg}
\end{eqnarray}

Consider the constant central $U(1)\subset \G$ subgroup of $\G$. Let
$\bG=\G/U(1)$ denote the quotient group. Then we have the fibration
\begin{eqnarray} BU(1)\rightarrow B\G\rightarrow
B\bG.\label{groupfibration} \end{eqnarray} It is shown in $\S 9$
of \cite{atiyah-bott} that this fibration is actually a product: 
\begin{eqnarray}B\G\sim BU(1)\times B\bG,\label{product}\end{eqnarray}
and so the generators of $H^*(B\G)$ give generators for $H^*(B\bG)$, with
$a_1$ being redundant, they are:  degree $1$ generators $b_1^j$, a degree
$2$ generator $f_2$, degree $3$ generators $b_3^j$ and a degree $4$
generator $a_2$. Now $B\bG$ is a free graded commutative algebra on
these generators and consequently (or equivalently from
(\ref{product})) we have:
\begin{eqnarray}
P_t(B\bG)=(1-t^2)P_t(B\G)=\frac{\left\{(1+t)(1+t^3)\right\}^{2g}}{(1-t^2)(1-t^4)}.\label{poincarebbg}
\end{eqnarray}

\subsubsection{The $\G$-equivariant cohomology of $\calC_d$}
For $d>0$ the ring $H_\G^*(\calC_d)$ 
can be described explicitly\footnote{For more details see
\cite{atiyah-bott}. For the homotopy type of $(\calC_d)_\G$ see
Corollary~\ref{homotopycdg}.}  
as a freely generated
graded commutative algebra by
degree $1$ elements $b_1^j$ and $b^j_2$ for $1\leq j\leq 2g$ and
degree $2$ elements $a_1^1$ and $a^2_1$. 
Thus the $\G$-equivariant Poincar\'e polynomial of $\calC_d$ is 
\begin{eqnarray}\G P_t(\calC_d)=\left(\frac{(1+t)^{2g}}{(1-t^2)}\right)^2.\label{poincarecd}\end{eqnarray}

\subsubsection{The $\G$-equivariant cohomology of $\calC_s$}

Since the constant central gauge transformations act trivially on
$\calC$ the factor group
$\bG^c=\G^c/\C^*$ acts on $\calC$ and moreover it acts freely on $\calC_s$. 
Now since (\ref{product}) is a product we have that as rings:
\begin{eqnarray*}
H^*_{\G}(\calC_s) &\cong &H^*(BU(1))\otimes
H^*_{\bG}(\calC_s)\\ &\cong& H^*(BU(1))\otimes
H^*(\calC_s/\bG)
\\ &\cong& H^*(BU(1))\otimes H^*(\tN),
\end{eqnarray*}
which for the Poincar\'e polynomials gives:
\begin{eqnarray} P_t(\tN)=(1-t^2)\G P_t(\calC_s).\label{poincaretn1}\end{eqnarray} 

Now putting (\ref{poincarec}), (\ref{poincarebg}), (\ref{poincarecd})
and (\ref{poincaretn1}) together yields the desired formula for the Poincar\'e
polynomial of $\tN$: 
\begin{eqnarray}
P_t(\tN) & = & P_t(B\bG)-\sum_{d=1}^{\infty}
t^{2(g+2d-2)}\frac{(1+t)^{4g}}{(1-t^2)} \nonumber \\ & = & 
\frac{\left\{(1+t)(1+t^3)\right\}^{2g}}{(1-t^2)(1-t^4)}-\frac{t^{2g}(1+t)^{4g}}{(1-t^2)(1-t^4)} \nonumber \\ & = &
(1+t)^{2g}\left(\frac{(1+t^3)^{2g}-t^{2g}(1+t)^{2g}}{(1-t^2)(1-t^4)}\right)
\label{poincaretn}
\end{eqnarray}
 
Finally, since the Shatz stratification is $\G$-perfect, the map $$H^*(B\G)\cong
H_\G^*(\calC)\to H_\G^*(\calC)$$ is surjective, thus the images of the
generators of $H^*(B\G)$  
give generators in $H^*_\G(\calC_s)$. As we saw above
$$H^*_\G(\calC_s)\cong H^*(BU(1))\otimes H^*(\tN),$$
thus the images of the 
generators of $H^*(B\G)$ give generators in
$H^*(\tN)$. Since one of them $a_1$ is redundant, we have the
following list of generators for the cohomology of $\tN$: degree $1$ generators $b_1^j$, a degree
$2$ generator $f_2$, degree $3$ generators $b_3^j$ and a degree $4$
generator $a_2$.   

In the next paragraph we explain the consequences of these result for
the cohomology of $\N$.   

\subsubsection{The cohomology of $\N$}

Let $\Gamma:=H^1(\Sigma,\Z_2)\cong \Z_2^{2g}\cong \ker(\sigma_2)$, 
where $\sigma_2:\J_0\rightarrow \J_0$ is given by $\sigma_2(L)=L^2$. Now
$\Gamma$ acts on $\N$ and $\J$ by tensoring with the corresponding line bundle in 
$\ker(\sigma_2)$ and also on $\N\times \J$ by the diagonal action. 
Then we have 
\begin{eqnarray}
\tN=(\N\times \J)/\Gamma.\label{split}
\end{eqnarray} 
Because 
$\Gamma$ acts trivially on $H^*(\J)$ and on $H^*(\N)$ (the latter was
first proved in \cite{harder-narasimhan}) we see that as rings 
\begin{eqnarray}
H^*(\tN)\cong (H^*(\N)\otimes H^*(\J))^\Gamma\cong H^*(\N)\otimes
H^*(\J).
\label{splitN}\end{eqnarray}
Thus for understanding the cohomology ring $H^*(\tN)$ it is enough
to know the cohomology ring $H^*(\N)$. Since we know a generator set
for the ring $H^*(\tN)$, it gives one for $H^*(\N)$. However the
classes $b_1^j$ will go to $0$. Thus a generator set is provided by: a
degree $2$ generator $f_2$, degree $3$ generators $b_2^s$ and a degree
$4$ generator $a_2$.

The generators we have just described differ from the original
generator set defined by Newstead in \cite{newstead2}.
He defined cohomology classes:
$\alpha\in H^2(\N)$, $\psi_i\in H^3(\N)$ and $\beta\in
H^4(\N)$, which appear in the K\"unneth decomposition of $c_2(\End(\E_\N))$:
\begin{eqnarray} 
c_2(\End(\E_\N))=2\alpha\otimes\sisi +
\sum_{i=1}^{2g}4\psi_i\otimes \xisi_i -
\beta\otimes 1 \label{kunneth}
\end{eqnarray}
in $H^4(\N\times \Sigma)\cong
\sum_{r=0}^{4} H^r(\N)\otimes H^{4-r}(\Sigma)$.
Here $\E_\N$ is the normalized 
rank $2$ universal bundle over $\N\times \Sigma$, i.e.
$c_1(\E_\N)=\alpha$ and $\E_\N\mid_{\{E\}\times \Sigma}\cong E$ for
every $E \in \N$. We will later use the following notations: 
$\gamma=-2\sum_{i=1}^g\psi_i\psi_{i+g}\in
H^6(\N)$,
$\gamma^*=2\gamma+\alpha\beta\in H^6(\N)$.

The relation between the Newstead and Atiyah-Bott generators can be
traced back by using the fact\footnote{Cf. \cite{atiyah-bott} p. 579-580} 
that the universal bundle $\U$ over $B\G$ 
restricts to $\E_\N$. The correspondence  is
given\footnote{Cf. \cite{earl}.}
by the formulae:
$$\al=2f_2-a_1,\hspace{1cm} \be=(a_1)^2-4a_2, \hspace{1cm}
\psi_i=b_2^i.$$

Finally, another consequence of (\ref{splitN}) and (\ref{poincaretn})
is the so-called Harder-Narasimhan formula for the
Poincar\'e polynomial of $\N$:
\begin{eqnarray}
P_t(\N)=\frac{(1+t^3)^{2g}-t^{2g}(1+t)^{2g}}{(1-t^2)(1-t^4)}
\label{poincaren}
\end{eqnarray}

\subsection{The resolution tower of $\tN$} 
We promised the existence
of a resolution tower for $\tN$ at the beginning of the present
section. It exists up to homotopy equivalence. Namely $\tN\sim
(\calC_0)_\bG$ and let us choose $\tN_k\sim (\calC_{\leq k})_\bG$ for
$k>0$. Thus homotopically we have the tower:
$$\tN=\tN_0\subset \tN_1\subset \dots \subset \tN_k \subset \dots,$$
which has the property that $$i_k^*:H^*(\tN_{k+1})\to H^*(\tN_k)$$ is
surjective\footnote{It is a consequence of the $\G$-perfectness of the Shatz
stratification.}. 
Thus if we consider $$\tN_\infty=\lim_{\longrightarrow}
\tN_k\sim \lim_{\longrightarrow}(\calC_{\leq k})_\bG \sim(\calC)_\bG\sim B\bG,$$
then we have that $$H^*(\tN_\infty)\to H^*(\tN)$$ is surjective, and
moreover $H^*(\tN_\infty)\cong H^*(B\bG)$ is a free commutative graded
algebra on universal classes, providing the picture we described at
the beginning of the present chapter. 

In the next section we explain how Kirwan used this resolution tower, or more
generally Proposition~\ref{propositionkirwan}, in order
to settle the Mumford conjecture, providing a complete set of
relations for $\tN$.

\subsection{Complete set of relations}

The ring structure of $H^*(\N)$ is described in terms of the so-called 
{\em Mumford relations} of the generators $\alpha,\beta$ and the
$\psi_i$'s. 
To explain this
consider the {\em virtual Mumford bundle} over $\tN$:
\begin{eqnarray*} \bM & = &
-{\pr_\tN}_!(\E_\tN\otimes\pr_\Sigma^*(L_p^{-1}))\\ & = &
-R^0{\pr_\tN}_*(\E_\tN\otimes \pr_\Sigma^*(L_p^{-1}))+R^1{\pr_\tN}_*
(\E_\tN\otimes 
\pr_\Sigma^*(L_p^{-1}))\in K(\tN).\end{eqnarray*} 
Using standard properties of stable bundles it can be shown that $R^0$ vanishes.
Thus $\bM$ is a vector bundle of rank $2g-1$. Its total Chern
class is a complicated\footnote{It was calculated by Zagier in
\cite{zagier}.} 
polynomial of the universal classes.
Since $\rank(\bM)=2g-1$, the Chern class
$c_{2g+r}(\bM)\in H^{4g+2r}(\tN)$ vanishes for $r\geq 0$. According to
(\ref{splitN}), the cohomology of $\tN$ is the tensor product of
$H^*(\J)$ and $H^*(\N)$. Thus if we write $\tau_S=\prod_{i\in S}\tau_i$
for $S\subset \{1\dots 2g\}$ and 
$$c_{2g+r}(\bM)=\sum_{S\subset \{1\dots 2g\}} \ze^r_S\otimes \tau_S$$ in the
K\"unneth decomposition of (\ref{splitN}) then we get the vanishing of
each $\ze^r_S$. Thus for every $r\geq 0$ and $S\subset \{1\dots 2g\}$ 
we get a relation
\begin{eqnarray}\ze^r_S\in \Qu[\alpha,\beta,\psi_i]\label{mumfordrelation}\end{eqnarray} 
of degree $d=4g+2r-\deg(\tau_S)$. The classes $\ze^r_S$ are called the
{\em Mumford relations}.

Mumford conjectured\footnote{Cf. \cite{atiyah-bott} p. 582.}  that the Mumford relations constitute
a complete set of relations of the cohomology ring of $\N$. Mumford's
conjecture was first settled by Kirwan in \cite{kirwan2}, by  using 
the method of Remark 1 after Proposition~\ref{propositionkirwan} for 
the Shatz stratification. The proof of \cite{kirwan2} goes by
building $\calR_d$ from the classes $\ze^r_S$ for $d>0$, and for $d=0$
setting $\calR_0=H^*_\G(\calC).$ The heart of the proof is to
show that $\calR_d$ satisfies the conditions of
Proposition~\ref{propositionkirwan}, which takes some pages of
calculation. Now Proposition~\ref{propositionkirwan} proves the
Mumford conjecture as explained in Remark 1 after it. 
We will give a purely geometric proof   
of the Mumford conjecture in Section~\ref{mumfordproof}. 

Now we explain an explicit description of the cohomology ring.

\subsection{Explicit description}

Recently, the following very explicit characterization of the
ring
$H^*(\N)$ has
been obtained by several authors 
\cite{baranovsky},\cite{king-newstead},\cite{siebert-tian} and
\cite{zagier}. To explain it first note that there is a natural action
of $Sp(2g,\Z)$ on $H^*(\N)$ induced by the obvious action on
$H^3(\N)$. The above mentioned sources prove that as an
$Sp(2g,\Z)$-algebra
$$H^*(\N)\cong \bigoplus_{k=0}^g\Lambda^k_0 H^3(\N)\otimes \Qu[\al,\be,\ga]/I_{g-k},$$
where $$\Lambda_0^k=\ker\left(\gamma^{g-k+1}:\Lambda^kH^3(\N)\to
\Lambda^{2g-k+2}H^3(\N)\right)$$
and $I_{g}$ is the relation ideal of the $Sp(2g,\Z)$-invariant
part 
$H_I^*(\N)$, i.e. $$H^*_I(\N)\cong
\Qu[\al,\be,\ga]/I_{g}.$$ The ideal $I_{g-k}$ is moreover described
as being generated by polynomials
$$\ze_{g-k},\ze_{g-k+1},\ze_{g-k+2}\in \Qu[\al,\be,\ga],$$
which are given  recursively 
by the following rule: \begin{eqnarray}(r+1)\ze_{r+1}=\alpha \ze_r + r \beta
\ze_{r-1} + 2\gamma \ze_{r-2},\label{recursive}\end{eqnarray} 
with initial conditions
$\ze_0=1,\ze_r=0$ for $r<0$.

Moreover an additive basis for $I_g$ is given by Zagier in \cite{zagier} of the form
\begin{eqnarray}\ze_{r,s,t} \mbox{\rm \ for all \ } r,s,t\geq 0 \mbox{\rm \ and \ }
r+s+t\leq g-1,\label{additivebasis}\end{eqnarray} 
where $\ze_{r,s,t}=\ze_{r,s}(2\ga)^t/t!$, and the polynomials
$\ze_{r,s}$ are given by a generating function:
\begin{eqnarray}\sum_{r,s\geq 0} \ze_{r,s}x^r y^s = 
\frac{e^{-2\ga x}/\be}{\sqrt{(1-\be y)^2-\be
x^2}}\left(\frac{1+x\sqrt{\be}-\be y}{1-x\sqrt{\be}-\be
y}\right)^{\ga^*/2\be \sqrt{\be}}. \label{generatingxirs} 
\end{eqnarray} 
It follows from (\ref{additivebasis}) that the Poincar\'e polynomial
of $H^*_I(\N)$ equals:
\begin{eqnarray} P_t^I(\N)=\sum_\stack{r,s,t\geq 0}{r+s+t\leq g-1} T^{r+2s+3t}.\label{invpoincareofN}\end{eqnarray}

\newpage
\section{The moduli space of Abelian Higgs bundles $T^*_\J$}
\label{toy}

As an instructive  example for the discussions in Section~\ref{statement}, we
consider here the moduli space of Abelian Higgs bundles. 

The tangent bundle of $\J$ is canonically isomorphic to $\J\times 
H^1(\Sigma,\calO_\Sigma)$. Thus by Serre duality 
$T^*_{\J}\cong \J\times
H^0(\Sigma,K)$ canonically. An element $\Phi\in (T^*_\J)_L\cong
H^0(\Sigma,K)$, can be thought of as a rank $1$ Higgs bundle\footnote{Cf. Definition~\ref{higgs}.}: 
$\cL=L\stackrel{\Phi}{\rightarrow} LK$.
Thus
we can think of $T^*_\J$ as the moduli space of rank $1$ Higgs
bundles. 

The cohomology of $T^*_\J$ is isomorphic to that of $\J$. However there
is an extra piece of cohomological information, namely the intersection
numbers in the compactly supported cohomology or in other words the
map:
$$j_\J:H^*_{cpt}(T^*_\J)\rightarrow H^*(T^*_\J).$$ 
Clearly this map is interesting only in the middle dimension, where 
both $H^{2g}_{cpt}(T^*_\J)$ and $H^{2g}(T^*_\J)$ are
one-dimensional. However the Euler characteristic of $\J$ is clearly
$0$, thus the self-intersection number of the zero section of $T^*_\J$ is
$0$, which shows that $j_\J$ vanishes. 

We can also determine the $L^2$-cohomology of $T^*_\J$ for the Riemann 
metric on $T^*_\J\cong  \J\times H^0(\Sigma,K)$ which is the product 
of the flat metrics on the two terms (this is the metric which we get
if we perform Hitchin's work in \cite{hitchin1} for the Abelian case).
From the Weitzenb\"ock decomposition of the Hodge Laplacian 
and the $L^2$-vanishing theorem of
Dodziuk \cite{dodziuk}, 
since the
metric is flat 
there are no non-trivial $L^2$ harmonic forms on $T^*_\J$. Thus
in the Abelian Higgs case the topology gives the harmonic space, as 
conjectured for the rank $2$ Higgs moduli space in 
Conjecture~\ref{conjecture}
and the moduli space of magnetic monopoles in \cite{sen}.

\newpage
\section{Moduli space of rank $2$ stable Higgs bundles $\M$}
\label{cohomologyofM}

Recall from Section~\ref{nagysrac} that $\tM$ denotes the moduli
space of rank $2$ stable Higgs bundles of degree $1$. 
The determinant gives a map 
$det_\M:\tM\rightarrow T^*_\J$, defined by 
$det_\M(E,\Phi)=(\Lambda^2E,\trace(\Phi))$. For any 
$\cL\in T^*_\J$ the fibre $det_\M^{-1}(\cL)$ will be
denoted by $\M_{\cL}$. Just as in the case of $\tN$ any two fibres of
$det_\M$
are 
isomorphic. Usually we will write $\M$ for
$\M_\cL$, when the Abelian Higgs bundle $\cL$ has zero Higgs field. 

Our main concern in the present thesis
is $\M$. Recall from Subsection~\ref{nagysrac} that it 
is a non-projective, smooth quasi-projective variety of 
dimension $6g-6$.  

Similarly to (\ref{split}) we have a $\Gamma$-action on $\tM$ and on
$T^*_\J$ such that:
$$\tM=(\M\times T^*_\J)/\Gamma.$$ This on the level of cohomology gives 
\begin{eqnarray}
H^*(\tM)\cong (H^*(\M)\otimes H^*(T^*_\J))^\Gamma\cong (H^*(\M))^\Gamma
\otimes H^*(\J) .\label{splitM}
\end{eqnarray}
In the case of $\M$ however we do not have the triviality of the action 
of $\Gamma$ on $H^*(\M)$, but nevertheless 
the cohomology ring of $\tM$ is determined
by the ring $(H^*(\M))^\Gamma$. 

There is quite little known about the
ring $H^*(\M)$. 
The Poincar\'e polynomial of it is calculated in \cite{hitchin1} using
Morse theory. We
now explain this:

\subsection{The Poincar\'e polynomial of $\M$}
We can outline Hitchin's Morse theory calculation in the language of
stratifications of Subsection~\ref{equivariantcoh} as follows. As we explained
in Subsection~\ref{rizsa1} the Morse function $\mu$ 
defines an upward stratification of $\M=\bigcup_{d=0}^{g-1}U_d$,
where $U_d=\left\{x\in \M:\lim_{z\to 0}z\cdot x\in F_d\right\}$, with $F_d$
denoting a component of the fixed point set of the $U(1)$-action. 
We call this the {\em Hitchin stratification}, 
because Hitchin used the
perfectness of this stratification to calculate the Poincar\'e
polynomial of $\M$ in Theorem 7.6 of \cite{hitchin1}. In particular
he proved in Proposition 7.1 of \cite{hitchin1}
that for $d>0$ the index of $F_d$ in $\M$ is $2(g+2d-2)$. On
the other hand the index of $F_d$ in $\M$ is the same as the real
codimension of $U_d$ in $\M$. The perfectness of the stratification
follows from \cite{kirwan1}. Now \cite{hitchin1} calculates the
Poincar\'e polynomial of $\M$ using the short exact sequences
(\ref{longexactofM}):
 \begin{eqnarray} 0\to H^*(U_d)\stackrel{(i_{d})_*}{\to}
H^*(\M_{\leq d})\stackrel{i^*_{<d}}{\to} H^*(\M_{<d})\to
0\label{short}, \end{eqnarray} 
and in turn the formula (\ref{poincare})
giving: 
\begin{eqnarray}
P_t(\M)=P_t(\N)+
\sum_{d=1}^{g-1} t^{2(g+2d-2)}P_t(U_d)=P_t(\N)+\sum_{d=1}^{g-1}
t^{2(g+2d-2)}P_t(F_d),\label{poincarem}
\end{eqnarray}
as $F_d$ is a deformation retract of $U_d$. For the $\Gamma$-invariant
part we have: \begin{eqnarray}
(P_t(\M))^\Gamma&=& (P_t(\N))^\Gamma+
\sum_{d=1}^{g-1} t^{2(g+2d-2)}(P_t(U_d))^\Gamma\nonumber
\\ &=&(P_t(\N))^\Gamma+\sum_{d=1}^{g-1}
t^{2(g+2d-2)}(P_t(F_d))^\Gamma\label{invpoinvarem}\end{eqnarray} from
the short exact sequences \begin{eqnarray} 0\to H^*(U_d)^\Gamma\stackrel{(i_{d})_*}{\to}
H^*(\M_{\leq d})^\Gamma\stackrel{i^*_{<d}}{\to} H^*(\M_{<d})^\Gamma\to
0\label{gammashortexact}, \end{eqnarray}
which can be obtained by noticing that the Hitchin
stratification is $\Gamma$-invariant and in turn that the short exact sequence
(\ref{short}) is in fact a sequence of $\Gamma$-modules.    
The next step is thus to
understand $H^*(F_d)$ and the action of $\Gamma$ on it:

\subsection{The cohomology of $F_d$}

We can construct $F_d$ as the moduli space of rank 2 degree 1
Higgs bundles of the form : $\compE$, where $$E=L\oplus L^{-1}\Lambda,$$ and 
\begin{eqnarray}\Phi=\left(\begin{array}{cc} 0&0\\ \phi &
0\end{array}\right),
\label{descriptionnd}\end{eqnarray}
with $\deg(L)=d$ and 
$0 \not= \phi\in H^0(\Sigma;L^{-2} \Lambda K).$ Note that
$\deg(L^{-2}\Lambda K)=2g-2d-1$, thus modulo non-zero
scalars $\phi$ is a point in $\Sigma_{\bd}$, where we used the
notation $\bd=2g-2d-1$.  

Now we determine $H^*(F_d)$. From the description (\ref{descriptionnd})
it follows that $F_d$ is isomorphic to the moduli
space of complexes: $L\stackrel{\phi}{\to}L^{-1}\Lambda K$, with
$\deg(L)=d$ and $\phi\in
H^0(\Sigma;L^{-2} \Lambda K).$ Thus if we
define
maps
$\J_d\to \J_{\bd}$ by sending 
$L\mapsto L^{-2}\Lambda K$ and the Abel-Jacobi map 
$\Sigma_{\bd}\to \J_{\bd}$ by sending $D\mapsto L(D)$, then the fibred product of the these
maps is:
\begin{eqnarray} F_d=\Sigma_{\bd}\times_{\J_{\bd}}
\J_d,\label{fibredproduct} \end{eqnarray}
We have the two projections $\pr_{\J_d}:F_d\to \J_d$ and 
$\pr_{\Sigma_{\bd}}:F_d\to \Sigma_{\bd}$. This last one is a
$2^{2g}$-fold cover, and it is induced by the action of $\Gamma\cong \Z_2^{2g}$ on
$F_d\subset \M$. Now Hitchin calculates the Poincar\'e polynomial of
$F_d$ by understanding the action $\Gamma$ on $H^*(F_d)$. He finds in
(7.13) of \cite{hitchin1}
that this action is not trivial and finds the cohomology of $F_d$ in
the form: $$H^*(F_d)=(H^*(F_d))^\Gamma\oplus 
V^{\bd}=H^*(\Sigma_{\bd})\oplus V^{\bd},$$ 
where $V^{\bd}$ is a faithful representation of $\Gamma$ 
in the middle $\bd$ degree. Moreover its dimension is
$(2^{2g}-1)\binomial{2g-2}{\bd}$. Thus for the Poincar\'e polynomial
he finds \begin{eqnarray}P_t(F_d) &=&(P_t(F_d))^\Gamma+(2^{2g}-1)\binomial{2g-2}{\bd}\nonumber\\ &=&
P_t(\Sigma_{\bd})+(2^{2g}-1)\binomial{2g-2}{\bd}.\label{poincarend}\end{eqnarray}

Now Hitchin calculates\footnote{There is, however, a small calculational mistake in 
(7.16) and (7.17) of \cite{hitchin1}, namely they should be multiplied
by $t^2$ in order to get the correct residue.\label{mistake}} 
the Poincar\'e polynomial of $\M$ from
(\ref{poincarem}), (\ref{poincarend}) and
(\ref{poincaresigman}). By subtracting the contributions of the non
$\Gamma$-invariant parts it follows\footnote{Keeping in mind the above
Footnote~\ref{mistake}.} from the formula of Hitchin that: 
$$(P_t(\M))^\Gamma=\frac{(1+t^3)^{2g}-t^{4g-2}(P(t))}{(1-t^2)(1-t^4)},$$
where $P(t)$ is some complicated polynomial of $t$. Consequently from
(\ref{splitM}) 
we have \begin{eqnarray}
P_t(\tM)=\frac{(1+t)^{2g}(1+t^3)^{2g}-t^{4g-2}\widetilde{P}(t)}{(1-t^2)(1-t^4)}=P_t(B\bG)-
t^{4g-2}\frac{\widetilde{P}(t)}{(1-t^2)(1-t^4)},
\label{poincaretm} \end{eqnarray} where $\widetilde{P}(t)=(1+t)^{2g}P(t)$
is some polynomial.
Comparing this result to (\ref{poincaretn}) we see that the Poincar\'e
polynomial of $\tM$ approximates the Poincar\'e polynomial of the
classifying space of $\bG$ roughly twice better than the Poincar\'e polynomial
of $\tN$! An explanation of this phenomenon will be provided in
Chapter~\ref{minfty}.

\subsection{Contribution of the present thesis}

In Chapter~\ref{compact} we make a detailed study of the
$\C^*$-action, 
compactify $\M$, and calculate the Poincar\'e polynomial of
the compactification. By doing so we establish the geometric
background for the cohomological calculations of the following
chapters, where we deal with two aspects of the cohomology of $\M$: 
the intersection numbers and the cohomology ring structure of $\M$. 

In the non-compact case the intersection numbers are in the compactly
supported cohomology. We will calculate all intersection numbers of
$H^*_{cpt}(\M)$ in Chapter~\ref{intersection}, by proving Theorem~\ref{main}.  

As we already mentioned in Section~\ref{statement}, for the calculation of the
sigma model of \cite{bershadsky-et-al} one needs to have a good
understanding of the ring structure of $(H^*(\M))^\Gamma$.
There is, however, no result about the cohomology ring $H^*(\tM)$ in
the literature. We will attempt to fill this gap in
Chapters~\ref{cohomology} and \ref{minfty}: We  find generators
for $H^*(\tM)$, and conjecture a complete set of relations for
$H^*_I(\M)$ in Chapter~\ref{cohomology}. In
Chapter~\ref{minfty}, we approach $H^*(\tM)$ in the general framework,
described at the beginning of the present chapter. Namely we show that
the tower of $U(1)$-manifolds $$\tM=\tM_0\subset \tM_1\subset \dots
\subset \tM_k \subset\dots$$ of moduli
spaces\footnote{Cf. Definition~\ref{poles}} 
of Higgs $k$-bundles with poles
is the right candidate for resolving the cohomology ring
$H^*(\tM)$. 

As byproducts of our considerations in Chapter~\ref{minfty}, we give a simple 
geometric proof of Mumford's conjecture, following ideas of
Kirwan, and show that the homotopy type of the above tower of Higgs
$k$-bundle moduli spaces shares the
Atiyah-Jones
property appearing in the theory of instantons.

\part{New results}

\chapter{Compactification}

\label{compact}

The main aim of this chapter is to investigate a canonical compactification of
$\M$: among other things we show that the compactification is projective,
calculate its Picard group, and calculate the Poincar\'e polynomial for the
cohomology.

We use a simple method to compactify non-compact
K\"ahler manifolds with a nice proper Hamiltonian $U(1)$-action via
Lerman's construction of symplectic cutting \cite{larman}. We use
this method to compactify $\M$.
Our approach is symplectic in nature and eventually
produces some fundamental results about the spaces occurring, using existing
techniques from the theory of symplectic quotients. 

We show that the compactification described here 
is a good example of  Yau's problem of
finding a complete Ricci flat metric on the complement of 
a nef anticanonical divisor
in a projective variety. 

Many of the results of this chapter can be easily generalized
to  other Higgs bundle moduli spaces, which have been extensively investigated
(see e.g. \cite{nitsure} and \cite{simpson}). As a matter of fact 
Simpson gave a 
definition of a similar compactificitation for these more general Higgs bundle
moduli spaces in Theorems 11.2 and 11.1 of \cite{simpson2} and 
in Proposition 17 of \cite{simpson3}, without investigating it in detail. 
For example, the projectiveness of the compactification 
is not clear from these definitions. One novelty of this chapter is the
proof of the projectiveness of the compactification in our case.

Since the compactification method used here is fairly 
general it is possible
to apply it to other K\"ahler manifolds with the above properties. It could
be interesting for instance to see how this method works for
the toric hyperk\"ahler manifolds of Goto \cite{goto} and Bielawski and Dancer
\cite{bielawski-dancer}. 

Finally we note that the compactification described
in the subsequent sections 
solves one half of the problem
of compactifying the moduli space $\M$, namely the `outer' half,
i.e.  shows
what the resulting spaces look like. The other half of the problem
the `inner' part, i.e. how this fits into the moduli space description
of $\M$, is treated in the recent paper of Schmitt \cite{schmitt}. Schmitt's
approach is algebro-geometric in nature, and concerns mainly the construction
of the right notion of moduli to produce $\cM$, thus complements 
our results.

\newpage
\section{Statement of results}

In this section we describe the structure of the chapter and list 
the results.

In Section~\ref{toymodel} we describe $\M_{toy}$ the moduli space of
parabolic Higgs bundle on $\Proj^1_4$, which will serve as a toy
example throughout this chapter.

In Section~\ref{stratificationM} we start to apply the ideas of
Section~\ref{rizsa}. Here we show that the Hitchin stratification,
constructed from the upward flows of Subsection~\ref{rizsa1} coincides
with the Shatz stratification coming from the stratification
(\ref{atiyahbottstrat}) on $\calC$.

In Section~\ref{nilpotent} (following ideas of Subsection~\ref{rizsa1}) 
we describe the nilpotent cone after Thaddeus
\cite{thaddeus1} 
and show that it coincides with the downward Morse
flow (Theorem~\ref{morse}). We reprove 
Laumon's theorem in our case, that the nilpotent cone is Lagrangian 
(Corollary~\ref{lagrange}). 

In Section~\ref{kaehler} we describe $Z$, the highest level K\"ahler quotient
of $\M$, while in \ref{kompakt} we analyse $\cM=\M\cup Z$. 
Here we follow the approaches of Subsection~\ref{rizsa2} and
Subsection~\ref{rizsa3}, respectively. 
Among others, we prove the following
statements:

\begin{itemize}

\item

$\cM$ is a compactification of $\M$, the moduli space of stable Higgs bundles
with fixed determinant and degree $1$ (Theorem~\ref{compactification}).

\item

$Z$ is a symplectic quotient of $\cal M$ by the circle action
$(E,\Phi)\mapsto (E,e^{i\theta}\cdot \Phi)$. $\overline{\cal M}$ is a symplectic
quotient  of ${\cal M}\times{\mathbf C}$ with respect to the circle action,
which is the usual one on $\cal M$ and multiplication on $\mathbf C$.

\item
While $\M$ is a smooth manifold, $Z$ is an orbifold, with only
${\mathbf Z}_2$
singularities corresponding to the fixed point set of the map
$(E,\Phi)\mapsto (E,-\Phi)$ on $\M$ (Theorem~\ref{Zorbi}), 
while similarly $\cM$ is
an orbifold with only ${\mathbf Z}_2$ singularities, and the singular locus of
$\cM$ coincides with that of $Z$ (Theorem~\ref{Morbi}).

\item

The Hitchin map $$\chi : \M \rightarrow {\mathbf C}^{3g-3}$$ extends to
a map $$\overline{\chi}:\cM\rightarrow {\mathbf P}^{3g-3}$$ which when
restricted
to $Z$ gives a map $$\overline{\chi}: Z \rightarrow {\mathbf P}^{3g-4}$$ whose generic
fibre is a Kummer variety corresponding to the Prym variety of the generic
fibre of the Hitchin map (Theorem~\ref{Zchi}, Theorem~\ref{Mchi}).
\item

$\cM$ is a projective variety (Theorem~\ref{Mproj}), with divisor $Z$ such that
$$(3g-2)Z=-K_{\cM},$$
the anticanonical divisor of $\cM$ (Corollary~\ref{vonat}).

\item

Moreover, $Z$ itself is a projective variety (Theorem~\ref{proj})
with an inherited holomorphic contact structure with
contact line bundle $L_Z$ (Theorem~\ref{contact})
and a one-parameter family
of K\"ahler forms $\omega_t(Z)$ (Theorem~\ref{csalad}). The Picard group
of $Z$ is described in Corollary~\ref{picard}. Moreover, the normal bundle
of $Z$ in $\cM$ is $L_Z$ which is nef by Corollary~\ref{triv}.

\item

Furthermore, $\cM$ has a one-parameter family of
K\"ahler forms $\omega_t(\cM)$,
which when restricted $Z$ gives the above $\omega_t(Z)$. 

\item

$Z$ is birationally equivalent to $P(T^*_\N)$ the projectivized
cotangent bundle of the moduli space of rank 2 stable bundles with fixed
determinant and odd degree (Corollary~\ref{Zjeno}). 
$\overline{\M}$ is birationally equivalent to
$P(T^*_\N\oplus {\cal O}_\N)$, the canonical compactification of
$T^*_{\N}$ (Corollary~\ref{Mjeno}).

\item

We calculate certain sheaf cohomology groups in Corollaries~\ref{H0Z} and
\ref{H1Z} and interpret some of these results as the equality
of certain infinitesimal deformation spaces. 

\item

The Poincar\'e polynomial of $Z$ is described in Corollary~\ref{PZ}, the
Poincar\'e polynomial of $\cM$ is described in Theorem~\ref{PM}. 

\item

We end Section~\ref{kompakt} by showing an interesting isomorphism between
two vector spaces: one contains information about the 
intersection of the components
of the nilpotent cone, the other says something 
about the contact line bundle $L_Z$ on $Z$.  
\end{itemize}

\newpage
\section{A toy model $\M_{toy}$}
\label{toymodel}

Unfortunately, even when $g=2$ the moduli space $\M$ is 
already $6$ dimensional, too big to serve as an instructive example. 
We rather choose $\M_{toy}$, the moduli space\footnote{These moduli 
spaces were considered by Yokogawa  \cite{yokogawa}.} of stable 
parabolic Higgs bundles on 
$\Proj^1$, with four marked points,  
 in order to show how
our later constructions work. 
We choose this example because
it is a complex surface, and can be constructed explicitly. 

We fix four distinct points on $\Proj^1$ and denote by $\Proj^1_4$ the 
corresponding complex orbifold. 
Let $P$ be the elliptic curve corresponding to $\Proj^1_4$.
Let $\sigma_P$ be the involution $\sigma_P(x)=-x$ on $P$. 
Thus, $P/\sigma_P$ is just 
the complex orbifold $\Proj^1_4$. The four fixed
points of the involution $x_1,x_2,x_3,x_4\in P$ correspond to the four marked
points on $\Proj^1_4$. Furthermore, let $\tau$ be the involution 
$\tau(z)=-z$ on $\C$. 

Consider now the quotient space 
$(P\times \C)/(\sigma_P\times \tau)$. 
This is a complex orbifold of dimension $2$
with four isolated $\Z_2$ quotient singularities at the points 
$x_i\times 0$. Blowing up these singularities we get a smooth complex 
surface $\M_{toy}$ with four exceptional divisors $D_1,D_2,D_3$ and $D_4$.  
Moreover the map $\chi:(P\times \C)\rightarrow \C$ sending 
$(x,z)\mapsto z^2$, 
descends to the quotient $(P\times \C)/(\sigma_P\times \tau)$ and sending 
the exceptional divisors to zero one obtains a map 
$\chi_{toy}:\M_{toy}\rightarrow\C$, 
with generic fibre $P$. 
The map $\chi_{toy}$ will serve as our toy Hitchin map.

There is a $\C^*$-action on $\M_{toy}$, coming from the standard
action on $\C$. 
The fixed point set of $U(1)\subset\C^*$ has five
 components: one is 
$\N_{toy}\subset \M_{toy}$ (the moduli space of stable parabolic
bundles on $\Proj_4^1$) which is the proper transform of 
$(P\times 0)/(\tau\times \sigma_P)=
P^1_4\subset (P\times \C)/(\sigma_P\times \tau)$ in $\M_{toy}$. 
The other four components consist of single points 
$\tilde{x_i}\in D_i,\ i=1,2,3,4$.

The fixed point set of the involution $\sigma:\M_{toy}\rightarrow \M_{toy}$
has five components, one of which is $\N_{toy}$, the others $E^2_i$ 
are the proper
transforms of the sets 
$(x_i\times \C)/(\sigma_P\times \tau)
\subset (P\times \C)/(\sigma_P\times \tau).$

\newpage
\section{The Shatz stratification on $\M$}
\label{stratificationM}

The results in Section~\ref{nagysrac} show that the K\"ahler manifold 
 $(\M,I,\omega)$ is equipped with a $\C^*$-action which restricts
to an $U(1)\subset \C^*$-action which is Hamiltonian with proper
moment map $\mu$. Moreover, $0$ is an absolute minimum for $\mu$. Therefore
we are in the situation described in 
Section~\ref{rizsa}. In the following sections
we will apply the ideas developed there to our situation and deduce
important properties of the spaces $\M$, $Z$ and $\cM$.

In the present and the following section we apply the general results of
Section~\ref{rizsa1} to $\M$. Moreover we identify the downward and
upward flows with important objects in the algebraic geometric
understanding of the Higgs moduli problem. 
We show that the Hitchin stratification coincides with the
Shatz stratification and that the downward Morse flow coincides with the
nilpotent cone. First we deal with the Hitchin stratification given
by the upward flows:

Recall from Subsection~\ref{rizsa1} that the upward flows give a
stratification on $\M$: $$\M=\bigcup_{d=0}^{g-1} U_d,$$
we call this the {\em Hitchin stratification} because
Hitchin calculated\footnote{Cf.  Subsection~\ref{cohomologyofM}.} the Poincar\'e polynomial of $\M$
using the perfectness of this stratification. 

We show in this section
that the Hitchin stratification coincides with the Shatz
stratification on $\M$. First we define the latter:

\begin{definition} Let $U^\prime_0\subset \M$ be the locus of points
$(E,\Phi)\in \M$ such that $E$ is stable, and moreover for $d>0$ let
$U^\prime_d\subset \M$ be the locus of points $(E,\Phi)\in \M$ such that
the  destabilizing line bundle of $E$ is of degree $d$. 
\end{definition}
  
\paragraph{{\it Remark.}} 1. This stratification can be easily
constructed in the gauge theory setting of
Subsection~\ref{gaugemk}. Namely one can pullback the Shatz
stratification $$\calC=\bigcup_{d=0}^\infty \calC_d,$$  by
$\pr_0:\calB_0\to \calC$ from $\calC$ to
$\calB_0$, restrict it to the open subset $(\calB_0)^s$, which, being
a $\bG^c$-invariant stratification 
induces a stratification on the quotient $\tM$ and in turn on $\M$. This is exactly
the same what we defined above.  

2. From (2) of Proposition 3.2 of \cite{nitsure} it
follows that if $d>g-1$ the locus $U^\prime_d$ is empty. 

3. The gauge theoretic construction in Remark 1 above implies 
that $\M=\bigcup_{d=0}^{g-1}U^\prime_d$ is a
(perfect) stratification in the sense of Subsection~\ref{equivariantcoh} and the
stratum $U^\prime_d$ has real codimension $2g+4d-4$ for
$d>0$. Moreover from the description (\ref{descriptionnd})
it easily follows that $F_d\subset U^\prime_d$.  

It follows that $U_d$ has the same
codimension as $U^\prime_d$. Moreover 
if $(E,\Phi)\in U^\prime_d$ then clearly its entire $\C^*$-orbit is contained in
$U^\prime_d$ thus it follows from (\ref{strat}) that $$\lim_{z\to 0} (
E,z\cdot \Phi)\in F\cap \overline{U^\prime_d}\subset F\cap \bigcup_{i\geq d} U^\prime_i=\bigcup_{i\geq d}F_d,$$ 
and in turn that $$\overline{U^\prime_d}\subset \bigcup_{i\geq d} U_d.$$ 
It follows that $U_{g-1}^\prime\subset U_{g-1}$ and hence
$U^\prime_{g-1}=U_{g-1}$ because they are closed submanifolds of
$\M$ with the same codimension and $U_{g-1}$ is connected. 
An inductive argument proves the following

\begin{proposition} The Shatz stratification coincides with the
Hitchin stratification i.e. $U^\prime_d=U_d$.
\end{proposition} 

\begin{remark} Thus the perfectness of the Shatz
stratification implies the perfectness of the Hitchin
stratification. Putting it in other words: Hitchin's
calculation\footnote{Cf. Section~\ref{cohomologyofM}.} of
the Poincar\'e polynomial of $\M$ is analogous to the Atiyah-Bott's
calculation\footnote{Cf. Section~\ref{cohomologyofN}.} of the
Poincar\'e polynomial of $\N$! 
In Section~\ref{homotopy} we will show why the two
calculations are profoundly related. 
\end{remark}

\newpage
\section{The nilpotent cone $N$}
\label{nilpotent}

We saw in Theorem~\ref{stratification} that the 
downward Morse flow is a deformation retract of $\M$, so it is
responsible for the topology, and as such it is an important object. 
On the other hand we will prove that
the downward Morse flow coincides with the nilpotent cone.

\begin{definition}
The nilpotent cone is the preimage of zero of the Hitchin map
$N=\chi^{-1}(0)$.
\end{definition}

The name `nilpotent cone' was given by Laumon, to emphasize the analogy with 
the nilpotent cone in a Lie algebra. 
In our context this is the most important fibre of the Hitchin map,
and the most singular one at the same time. We will show below that the
nilpotent cone is a central notion in our considerations. 

Laumon in \cite{laumon} investigated the nilpotent cone in a much more
general context and showed its importance in the Geometric Langlands 
Correspondence. Thaddeus in \cite{thaddeus1} concentrated on our case, 
and gave the exact description of the nilpotent cone. In what follows
we will reprove some of their results.  

The following assertion was already stated in \cite{thaddeus1} 
which will turn out to be crucial in some of our considerations.

\begin{theorem} The downward Morse flow coincides with the nilpotent
cone.
\label{morse}
\end{theorem}

\begin{proof} As we saw in Theorem~\ref{stratification} the downward Morse flow can be identified with 
the set of points in $\M$ whose $\C^*$-orbit
is relatively compact in $\M$.

Since the nilpotent cone is invariant under the $\C^*$-action
and compact ($\chi$ is proper) we immediately get that the nilpotent
cone is a subset of the downward Morse flow.

On the other hand if a point in $\M$ is not in the nilpotent cone then the
image of its $\C^*$-orbit by the Hitchin map is a line in $\C^{3g-3}$, 
therefore cannot be relatively compact.
\end{proof}

Laumon's  main
result is the following assertion\footnote{Cf. Theorem $3.1$ in \cite{laumon}.}, 
which we prove in our case:

\begin{corollary}[Laumon] The nilpotent cone is a Lagrangian subvariety 
of $\M$
with respect to the holomorphic symplectic form $\omega_h$.
\label{lagrange}
\end{corollary}

\begin{proof} The Hitchin map is a completely integrable Hamiltonian
system, and the nilpotent cone is a fibre of this map, so it is
coisotropic. Therefore it is Lagrangian if and only if its dimension is
$3g-3$.

On the other hand the nilpotent cone is exactly the downward Morse flow and
we can use Hitchin's description of the critical submanifolds in 
\cite{hitchin1}, giving that
the sum of the index and the real dimension of
any critical submanifold is $6g-6$.
We therefore conclude that the complex dimension of the downward Morse flow 
(i.e. the nilpotent cone) is $3g-3$.
\end{proof}

\begin{remark} Nakajima's Proposition $7.1$ in \cite{nakajima} states that
if $X$ is a K\"ahler manifold with a $\C^*$-action and a 
holomorphic symplectic form $\omega_h$ of homogeneity $1$ then the downward
Morse flow of $X$ is Lagrangian with respect to $\omega_h$. Thus Nakajima's 
result and Theorem~\ref{morse} together give an alternative proof of the 
theorem. We preferred the one above for it concentrates on the specific 
properties of $\M$. 
\end{remark}

From the above proof we can see that for higher rank Higgs bundles
Laumon's theorem is equivalent to the assertion that every critical
submanifold contributes to the middle dimensional cohomology, i.e the
sum of the index and the real dimension of
any critical submanifold should always
be half of the real dimension of the corresponding moduli space.

Using the results of \cite{gothen1} one easily shows 
that the above statement also holds 
for the rank $3$ case. 
Gothen could show directly the above statement
for any rank and therefore gave an alternative proof of Laumon's theorem
in these cases \cite{gothen2}. 

\begin{corollary} The middle dimensional homology $H_{6g-6}(\M)$ of $\M$ is
freely generated by the homology classes of irreducible components
of the nilpotent cone and therefore has dimension $g$.
\label{middle}
\end{corollary}

\begin{proof} We know that each component of $N$ is a projective variety of
dimension $3g-3$. $N$ is a deformation retract of $\M$, therefore
the middle dimensional homology of $\M$ is generated by the homology classes
of the components of $N$. Furthermore, from the Morse picture,
components of $N$ are in a one to
one correspondence with the critical manifolds of $\M$, so 
there are $g$ of them.
The result follows. 
\end{proof}

We finish this section with Thaddeus's description\footnote{See
\cite{thaddeus1} and Proposition~\ref{cone}, cf. also \cite{laumon}.} 
of the nilpotent cone.

\begin{theorem} The nilpotent cone is the union of $D_0=\N$ and the
downward flows $D_d$, which are total
spaces  of vector bundles $E^-_d$ over $F_d$, where $E^-_d$ is the negative
subbundle of $T_\M\mid_{F_d}$.

Moreover, the restricted action of $\C^*$ on $N$ is just the inverse
multiplication 
on the fibres. 
\label{th-nil}
\end{theorem}

\begin{proof} This follows directly from Theorem~\ref{stratification} and
Theorem~\ref{morse}, noting that by Hitchin's description of the weights
of the circle action on $T_\M\mid F_d$ in the proof of 
Proposition $7.1$ of \cite{hitchin1},
we have that there is only one negative weight. Therefore the $\beta$-fibration
of Theorem~\ref{stratification} is a vector bundle in this case. 
The result follows. 
\end{proof}

\begin{remark} From the description of $E^-_d$ in \cite{thaddeus1} and that of
$E^2_d$, a component of the fixed point set
of the involution $\sigma(E,\Phi)=(E,-\Phi)$,
in \cite{hitchin1}, one obtains the remarkable fact 
that the vector bundle $E^-_d$ is actually dual to $E^2_d$. This is
not so surprising if we observe that $E^-_d$ is the weight $-1$ and
$E^2_d$ is the weight $2$ component of the $\C^*$-equivariant bundle
$T_\M\mid_{F_d}$ and these are naturally dual, because of the
homogeneity $1$ holomorphic symplectic structure on $\M$! 
\end{remark}

\begin{example} In our toy example we have the elliptic fibration 
$\chi_{toy}:\M_{toy}\rightarrow \C$, with the only singular fibre 
$N_{toy}=\chi_{toy}^{-1}(0)$, the toy nilpotent cone. We have now
the decomposition $$N_{toy}=\N_{toy}\cup \bigcup^4_{i=1}D_i,$$ where
we think of $D_i$ as the closure of $E_i$, the total space of the trivial
line bundle on $\tilde{x_i}$. 
 
The possible singular fibres of elliptic fibrations
have been 
classified by Kodaira\footnote{Cf. \cite{barth-peters-van}, p. 150.}. According
to this classification $N_{toy}$ is of type $I_0^*(\tilde{D}_4)$.     
\end{example}

\newpage
\section{The highest level K\"ahler quotient $Z$}
\label{kaehler}

In this section we apply the ideas of Subsection~\ref{rizsa2} to our 
situation. 

\begin{definition} Define for every non negative $t$ the K\"ahler
quotient $$Q_t=\mu^{-1}(t)/U(1).$$ 
As the complex structure of the K\"ahler quotient
depends only on the connected component of the regular values of $\mu$, we
can define $Z_d=Q_t $ for $c_{d}<t<c_{d+1}$ as a complex orbifold\footnote{We set
$c_g=\infty$.}. Similarly, we define $X_{Z_d}$ to be $\M^{min}_t$ for 
$c_{d}<t<c_{d+1}$.

For simplicity let the highest level quotient $Z_{g-1}$
be denoted by $Z$ and the corresponding principal $\C^*$-bundle $X_{Z_{g-1}}$ 
by $X_{Z}$. 
\end{definition}

In the spirit of Theorem~\ref{quotients} we have the following

\begin{theorem}$Z_d$ is a complex orbifold with only
$\Z_2$-singularities, the singular
locus is diffeomorphic to some union of projectivized vector bundles
$P(E^2_i)$:
$$\sing(Z_d)=\bigcup_{0<i\leq d} P(E^2_i), $$
where $E^2_i\subset\M$ is the total space of a vector bundle over $F_i$ and 
is a component of the fixed point set of the involution 
$\sigma(E,\Phi)=(E,-\Phi)$.  

\label{Zorbi}
\end{theorem}

\begin{proof} 
The induced action of $U(1)$ on $\C^{3g-3}$ by the Hitchin map
is multiplication by $e^{2i\theta}$ so an orbit of $U(1)$ on
$\M\setminus N$ is a non trivial double cover 
of the image orbit on $\C^{3g-3}$.
On the other hand by Thaddeus' description of $N$ (Theorem~\ref{th-nil}) 
it is clear 
that if a point of $N$
is not a fixed point of the circle action, then the stabilizer is trivial at
that point.

Summarizing these two observations we obtain that if a point of $\M$ is
not fixed by $U(1)$, then its stabilizer is either trivial or $\Z_2$ . The
latter case occurs exactly at the fixed point set of the involution $\sigma$.
The statement now follows from Theorem~\ref{quotients}.
\end{proof}

\begin{proposition}
$Z_d$ and $Z_{d+1}$ are related by a blowup followed by a blowdown.
Namely, $Z_d$ blown up along $P(E^-_d)$ is the
same as the singular quotient $Q_{c_d}$ blown up along $F_d$ (its singular
locus),
which in turn gives $Z_{d+1}$ blown up at $P(E^+_d)$.
Moreover, this birational equivalence is an isomorphism outside an analytic set
of codimension at least $3$.
\label{birational}
\end{proposition}

\begin{proof} The first bit is just the 
restatement of Theorem~\ref{quotients} in
our setting.

The second part follows because $$\dim(P(E_d^-))=3g-3-1<6g-6-2$$ and
$$\dim(P(E^+_d))=3g-3+2g-2d-1-1<6g-6-2$$ for $g>1$.
\end{proof}

\begin{corollary} $Z=Z_{g-1}$  is birationally equivalent to $P(T^*_\N)=Z_0$.
Moreover this gives an isomorphism in codimension $>2$.
\label{Zjeno}
\end{corollary}

\begin{proof} Obviously $X_{Z_0}$ is $T^*_\N$, and therefore by 
Theorem~\ref{complex} $Z_0$ is isomorphic 
to the projectivized cotangent bundle $P(T^*_\N)$. The
statement follows from the previous theorem. 
\end{proof}

\begin{corollary} $Z$ has  Poincar\'e polynomial
$$ P_t(Z)=\frac{t^{6g-6}-1}{t^2-1}P_t(\N) +
\sum _{d=1}^{g-1}\frac{t^{6g-6}-t^{2g-4+4d}}{t^2-1}P_t(F_d),$$
where $F_d$ is a $2^{2g}$-fold cover of $\Sigma_{\bd}$.
\label{PZ}
\end{corollary}

\begin{proof} One way to derive this formula is through Kirwan's formula
in \cite{kirwan1}. We use the above blowup, blowdown picture instead. 
This approach\footnote{Cf. \cite{thaddeus2}.} is due to Thaddeus.

Applying the formula in \cite{griffiths-harris},p.605 twice we get that
$$P_t(Z_{d+1})-P_t(Z_d)=P_t(PE^+_d)-P_t(PE^-_d).$$

On the other hand for a projective bundle on a manifold
$P\rightarrow M$ with fibre $\Proj^n$ one has\footnote{Cf. \cite{griffiths-harris} p.606.}
$$P_t(P)=\frac{t^{2n+2}-1}{t^2-1}P_t(M).$$

Hence the formula follows.
\end{proof}

\begin{remark} All the Poincar\'e polynomials on the right hand side
of the above formula were calculated\footnote{For $P_t(\N)$ see
(\ref{poincaren}), 
for $P_t({F_d})$ see (\ref{poincarend}).} in Chapter~\ref{cohomologyofmoduli}.
\end{remark}

We will determine the Picard group of $Z$ exactly. First we
define some line bundles on several spaces.

\begin{notation} Let 

\begin{itemize} 
\item $\calL_\N$ denote the ample generator\footnote{Cf. \cite{drezet-narasimhan}.}
of the Picard group of $\N$,

\item
$\calL_{PT^*_\N}$ be its pullback 
to $PT^*_\N$, 
\item
$\calL_Z$ 
denote the corresponding line
bundle\footnote{Cf. Corollary~\ref{Zjeno}.}
on $Z$,

\item 
$L_{PT^*_\N}$ be the dual of the tautological line bundle on $PT^*_\N$,

\item 
$L_Z=X^*_Z\times_{\C^*}\C$ denote the corresponding line orbibundle
on $Z$. 
\end{itemize}
\end{notation}

\begin{corollary} $\Pic(Z)$, the Picard group of $Z$, is of rank $2$ over
$\Z$ and is freely generated by $\calL_Z$ and $L_Z$.
\label{picard}
\end{corollary}

\begin{remark} The Picard group of $Z$ is the group of invertible sheaves on
$Z$. As the singular locus of $Z$ has codimension $\geq 2$, this group can
be thought of as the group of holomorphic line orbibundles on $Z$. Namely, in
this case the restriction of a holomorphic line orbibundle to 
$Z\setminus \sing(Z)$ gives a one-to-one correspondence between holomorphic
line orbibundles on $Z$ and holomorphic line bundles on $Z\setminus \sing(Z)$,
by the appropriate version of Hartog's theorem.
\end{remark}

\begin{proof} It is well known\footnote{Cf. \cite{drezet-narasimhan}.} that $\Pic(\N)$ is freely generated by one ample
line bundle $\calL_\N$ therefore is of rank $1$.
Thus $\Pic(P(T^*_\N))$
is  of rank $2$ and freely generated by $\calL_{PT^*_\N}$
the pullback of $\calL_\N$ and the dual of the 
tautological line bundle $L_{PT^*_\N}$.
From Corollary~\ref{Zjeno} $\Pic(Z)$ is
isomorphic with $\Pic(P(T^*_\N))$ therefore is of rank $2$, and freely
generated by $\calL_Z$ and $L_Z$, where $\calL_Z$ is isomorphic
to $\calL_{PT^*_\N}$ and $L_Z$ is isomorphic to $L_{PT^*_\N}$ outside
the codimension $2$ subset of Corollary~\ref{Zjeno}.
\end{proof}

\begin{definition} A contact structure on a compact complex orbifold
$Z$ of complex dimension $2n-1$ is given by the following data: 

\begin{enumerate}
\item
a contact line orbibundle $L_Z$ such that $L_Z^n=K_Z^{-1}$, where  $K_Z$ is
the line orbibundle of the canonical divisor of $Z$,

\item

a complex contact form $\theta\in H^0(Z,\Omega^1(Z)\otimes L_Z)$ a 
holomorphic $L_Z$ valued $1$-form, such that 
\begin{eqnarray}
0\neq \theta\wedge (d\theta)^{n-1}\in H^0(Z,\Omega^{2n-1}(Z)\otimes K_Z^{-1})=
H^0(Z,\calO_Z)=\C 
\label{degeneralt}
\end{eqnarray}
is a nonzero constant.
\end{enumerate} 
\end{definition}
 
\begin{theorem} There is a canonical holomorphic contact structure
on $Z$ with contact line orbibundle $L_Z$.
\label{contact}
\end{theorem}

\begin{proof} This contact structure can be defined by the construction of 
Lebrun as in \cite{lebrun} Remark $2.2$. 
We only have to note that the holomorphic 
symplectic form $\omega_h$ on $\M$ is of homogeneity $1$. 

The construction goes as follows. If $\pr_Z:X^*_Z\rightarrow Z$ denotes the
canonical projection of the principal $\C^*$-orbibundle $X^*_Z$ the
dual of $X_Z$, then $\pr_Z^*(L_Z)$ is canonically trivial with the 
canonical section having homogeneity $1$. Thus in order to give a complex
contact form $\theta\in H^0(Z,\Omega^1(Z)\otimes L_Z)$ it is sufficient
to give a $1$-form $\pr_Z^*\theta$ on $X^*$ of homogeneity $1$. This can
be defined by $\pr_Z^*\theta=i(\xi)\omega_h$, where $\xi\in H^0(\M,T_\M)$
is the holomorphic vector field generated by the $\C^*$-action. The 
non-degeneracy condition (\ref{degeneralt}) is exactly equivalent to requiring
that the closed holomorphic $2$ form $\omega_h$ satisfy $\omega_h^n\neq 0$. 
This is the case as $\omega_h$ is a holomorphic symplectic form. 

The result follows.
\end{proof}

We will be able to determine the line orbibundle $L_Z$ explicitly. 
For this, consider the Hitchin map $\chi:\M\rightarrow \C^{3g-3}$. 
As it is equivariant
with respect to the $\C^*$-action, $\chi$ induces a map 
$$\overline{\chi}:Z\rightarrow \Proj^{3g-4}$$
on Z. The generic fibre of this map is easily seen to be the Kummer
variety corresponding to the Prym variety (the Kummer variety of an
Abelian variety is the quotient of the Abelian variety by the involution
$x\rightarrow -x$), the generic fibre of the Hitchin map. 
Thus we have proved 

\begin{lemma} There exists a map $\chi:Z\rightarrow \Proj^{3g-4}$ the reduction
of the Hitchin map onto $Z$, for which the generic fibre is a Kummer variety.
\label{Zchi}
\end{lemma}

\begin{remark} This observation was already implicit in Oxbury's
thesis\footnote{Cf. $2.17$a of \cite{oxbury}.}. 
\end{remark}

The following theorem determines the line bundle $L_Z$ in terms of the
Hitchin map.

\begin{theorem} $L^2_Z=\overline{\chi}^*{\cal H}_{3g-4}$ where 
${\cal H}_{3g-4}$ is the hyperplane bundle on $\Proj^{3g-4}$.
\label{feco}
\end{theorem}

\begin{proof} We know from Corollary~\ref{picard} that
$\overline{\chi}^* {\cal H}_{3g-4}={\calL}_Z^k\otimes L_Z^l$ for some 
integers $k$ and $l$.

We show that $k=0$. For this consider the pullback of $\calL_Z$ onto
$\M\setminus N$ the total space of the principal $\C^*$-orbibundle $X^*_Z$. 
This line orbibundle extends to $\M$ as $\calL_\M$
and restricts to $T^*_\N$ as the pullback
of $\calL_{PT^*_\N}$ by construction. $c_1(\calL_\M)$ is not trivial when
restricted to $\N$ 
(namely it is $c_1(\calL_\N)$, since this bundle is ample) 
therefore is not trivial when
restricted to a generic fibre of the Hitchin map. We can deduce
that $c_1(\calL_Z)$ is not trivial on the generic fibre of $\overline{\chi}$.

On the other hand $L_Z$ restricted to a generic fibre of
$\overline{\chi}$ can be described as follows. Let this Kummer variety
be denoted by $K$, the corresponding Prym variety by $P$. Form
the space $P\times \C^*$, the trivial principal $\C^*$-bundle
on $P$ and quotient it out by the involution
$\tau(p,z)=(-p,-z)$. The resulting space is easily seen to be
the $\C^*$-orbit of the Prym $P$ in $\M$, therefore the total space
of the principal $\C^*$-orbibundle $L^*_Z\setminus (L^*_Z)_0$ on $K$. 
Hence $L_Z^2$ is the trivial line orbibundle on $K$. Thus $c_1(L_Z\mid_K)=0$.

Now
$\overline{\chi}^*{\cal H}_{3g-4}$ is trivial on the Kummer variety.
Hence the assertion
$k=0$.

The rest of the proof will follow the lines of Hitchin's proof of Theorem
$6.2$ in \cite{hitchin2}. We show that $l=2$.

The sections of $L_Z$ can be 
identified with holomorphic functions
homogeneous of degree $2$ on the principal
$\C^*$-orbibundle $X_Z=L_Z^*\setminus (L_Z^*)_0$. As $N$ is of codimension
$\geq 2$ such functions extend to $\M$. Since the Hitchin map is proper,
these functions are constant on the fibres of the Hitchin map, therefore
are the pullbacks of holomorphic functions on $\C^{3g-3}$ of
homogeneity $1$ which can be identified with the holomorphic sections of
the hyperplane bundle ${\cal H}_{3g-4}$ on $P(\C^{3g-3})=\Proj^{3g-4}$.
\end{proof}

\begin{corollary} If $n$ is odd, there are natural isomorphisms

$$H^0(Z,L_Z^n)\cong H^0(\N,S^n{T_\N})\cong0,$$ whereas if $n$ is even, then
$$ H^0(Z,L_Z^n)\cong H^0(\N,S^n{T_\N})
\cong H^0(\Proj^{3g-4},{\cal H}^{\frac{n}{2}}_{3g-4}).$$
\label{H0Z}
\end{corollary}

\begin{proof} We show that $H^0(Z,L^n_Z)\cong H^0(\N,S^n(T_\N))$ for every
$n$, the rest of the theorem will follow from Theorem $6.2$ of 
\cite{hitchin2}. 

By Proposition~\ref{birational} we get that 
$H^0(Z,L^n_Z)\cong H^0(PT^*_\N,L^n_{PT^*_\N})$. Let $\pr_\N:PT^*_\N\rightarrow\N$
denote the projection. It is well known that the Leray spectral sequence
for $\pr_\N$ degenerates at the $E^2$ term. Moreover, we have\footnote{Cf. \cite{hartshorne} Theorem $5.1$b.} that 
$R^i(\pr_\N)_*(L^n_{PT^*_\N})=0$ if $0<i<{3g-4}$. Therefore 
$H^0(PT^*_\N,L^n_{PT^*_\N})\cong H^0(\N,(\pr_\N)_*(L^n_{PT^*_\N}))$. Finally 
the sheaf $(\pr_\N)_*(L^n_{PT^*_\N})$ is $S^n(T_\N)$, 
which proves the statement.
\end{proof}   

We can moreover determine the first cohomology group corresponding to
the infinitesimal deformations of the holomorphic contact structure on $Z$
and can interpret it in a nice way.

\begin{corollary} There are canonical isomorphisms 
$$H^1(Z,L_Z) \cong (H^1(\M,{\calO}_{\M}))_1\cong H^1(\N,T_\N)
\cong H^1(\Sigma, K_\Sigma^{-1}),$$
where $(H^1(\M,{\calO}_{\M}))_1\subset H^1(\M,{\calO}_{\M})$ is 
the vector space of elements of $H^1(\M,{\calO}_{\M})$ homogeneous of
degree $1$.
\label{H1Z}
\end{corollary}

\begin{proof} We may use the cohomological version\footnote{Cf. \cite{scheja}.} of Hartog's theorem 
to show that $$H^1(Z,L_Z)\cong H^1(PT^*_\N,L_{PT^*_\N}),$$
as $Z$ and  $PT^*_\N$ are isomorphic on an analytic set of codimension $\geq 3$
from Proposition~\ref{birational}.

The proof of the other isomorphisms can be found in \cite{hitchin3}.
\end{proof}

\begin{remark} We can interpret this result as saying that the deformation of
the complex structure on $\Sigma$ corresponds to the deformation of complex
structure on $\N$, to the deformation\footnote{Cf. \cite{lebrun}.} of holomorphic contact structure on $Z$
and to the deformation of the holomorphic symplectic
structure of homogeneity $1$ on $\M$.
\end{remark}

As an easy corollary of the above we note the following

\begin{corollary} The line orbibundle $L_Z$ is nef but neither trivial
nor ample.
\label{triv}
\end{corollary}

\begin{proof} The line bundle $L_Z$ is certainly not ample since $c_1(L_Z)$ is trivial on the
Kummer variety.
On the other hand $L^2_Z$ being the pullback of an ample bundle is not
trivial and is nef itself,
hence the result. 
\end{proof}

The next theorem will describe the inherited K\"ahler structures of $Z$.
Considering the one-parameter family of
K\"ahler quotients $Q_t$, for $t>c_{g-1}$ we
get a one-parameter family of K\"ahler forms $\omega_t$ on $Z$.
Theorem $1.1$ from \cite{duistermaat-heckman} gives the
following result for our case\footnote{Cf. Theorem~\ref{quotients}.}. 

\begin{theorem}[Duistermaat, Heckman]
The complex orbifold $Z$ 
has a one-para\-meter family of K\"ah\-ler forms $\omega_t$,
$t>c_{g-1}$ such that 
$$[\omega_{t_1}(Z)]-[\omega_{t_2}(Z)]=(t_1-t_2)c_1(L_Z)$$
where $t_1,t_2>c_{g-1}$ and $[\omega_t]\in H^2(Z,\R)$ is the cohomology
class of $\omega_t$.
\label{csalad}
\end{theorem}

Many of the above results will help us to prove the following theorem.

\begin{theorem} $Z$ is a projective algebraic variety.
\label{proj}
\end{theorem}

\begin{proof} 
By the Kodaira embedding theorem for orbifolds\footnote{Cf. \cite{bailey}.} we have only
to show that $Z$ with a suitable K\"ahler form is a Hodge orbifold, i.e.
the K\"ahler form is integer. For this to see we show that the
K\"ahler cone of $Z$ contains a subcone, which is open in $H^2(Z,\R)$.
This is sufficient since such an open subcone must contain an
integer K\"ahler
form i.e. a Hodge form.

Since Corollary~\ref{picard} shows that $\Pic_0(Z)$ is trivial, 
by Corollary~\ref{triv} we see that $c_1(L_Z)\neq 0$. Therefore the 
 previous theorem exhibited a half line in the K\"ahler cone of $Z$. Thus
to find an open subcone in the $2$ dimensional vector space $H^2(Z,\R)$
(Corollary~\ref{picard}) it is sufficient to show that this line 
does not go through
the origin or in other words $c_1(L)$ is not on the line.
But this follows from Corollary~\ref{triv}, because $L$ being not ample
$c_1(L)$ cannot contain a K\"ahler form. Hence the result.
\end{proof}

\begin{remark} We see from this proof that $c_1(L_Z)$ lies on the closure
of the K\"ahler cone, thus $L_Z$ is nef. This reproves a statement of
Corollary~\ref{triv}.
\end{remark}

\begin{example}  In the case of the toy example the lowest level
K\"ahler quotient $Z_0$ is    
the projectivized cotangent bundle 
$PT^*_{\N_{toy}}$ of $\N_{toy}$, which is isomorphic to 
$\N_{toy}=\Proj^1$, and
the blowups and blowdowns 
add the four marked points to $\Proj^1$.  
Therefore $Z_{toy}$ is isomorphic to the orbifold $\Proj^1_4$, where
the marked points correspond to the fixed point set of the involution 
$\sigma$, namely these are the projectivized bundles $PE^2_i$, i.e. points.

Moreover the principal $\C^*$-orbibundle $X_{Z_{toy}}$ on $\Proj^1_4$ has the
form 
$$X_{Z_{toy}}=(P\times \C^*)/(\sigma_P\times \tau).$$ 

Thus in the toy example, 
unlike in the ordinary Higgs case,  
we have $c_1(L_{Z_{toy}})=0$. This latter assertion can be seen using 
\ref{feco} and noting that the target of the reduced toy Hitchin map
$\overline{\chi}_{toy}:Z_{toy}\rightarrow \Proj^0$ is a point. 

There is another difference, namely the Picard group of $Z_{toy}$ is of 
rank $1$, because $L_{Z_{toy}}^2$ is the trivial bundle on $Z_{toy}$. 
\end{example}

In the next section we show how to compactify $\M$ by sewing in $Z$ at
infinity.

\newpage
\section{The compactification $\cM$}
\label{kompakt}

In this section we compactify $\M$ by adding to each non-relatively 
compact
$\C^*$-orbit an extra point i.e. sewing in $Z$ at infinity. 
Another way of saying
the same is to glue together $\M$ and $E$ the total space of $L_Z$ along the
principal $\C^*$-orbibundle $X^*_Z=E\setminus E_0=\M\setminus N$. 
To be more
precise we use the construction of Lerman, called the symplectic 
cut\footnote{Cf. Subsection~\ref{rizsa3} and \cite{larman}.}.

Since the complex structure on the K\"ahler quotients depends only on the
connected component of the level, we can make the following definition.

\begin{definition} Let $\overline{\M}_d$ denote
the underlying compact complex orbifold of
the K\"ahler quotient of $\M\times \C$ by  the product $U(1)$-action
$$\overline{\M}_{\mu<t}=(\mu+\mu_\C)^{-1}(t)/{U(1)},$$ with $c_d<t<c_{d+1}$.

Let $X_{\cM_d}$ denote the corresponding principal $\C^*$-bundle on $\cM_d$.
For simplicity we let $\overline{\M}$ denote $\overline{\M}_{g-1}$ and $X_\cM$ denote
$X_{\cM_{g-1}}$.
\end{definition}

As a consequence of the construction of symplectic cutting we have the
following theorem\footnote{Cf. Theorem~\ref{cut}.}:

\begin{theorem} The compact orbifold
$\cM=\M\cup Z$ is a compactification of $\M$ such
that $\M$ is an open complex submanifold and $Z$ is a codimension one 
suborbifold, i.e. a divisor.

Moreover, $\C^*$ acts on $\cM$ extending the action on $\M$ with the points
of $Z$ being fixed.
\label{compactification}
\end{theorem}

In addition to the above we see that we have another
decomposition $\cM=N\cup E$ of $\cM$ into the
nilpotent cone and the total space $E$
of the contact line bundle $L_Z$ on $Z$. Thus the
compactification by symplectic cutting produced the same orbifold
as the two constructions we started this section with.

We start to list the properties of $\cM$. We will mention properties
analogous to properties of $Z$ (these correspond to the fact that both
spaces were constructed by a K\"ahler quotient procedure)
and we will clarify the
relation between $Z$ and $\cM$.

Theorem~\ref{cut} and Theorem~\ref{quotients} give the following result
in our case. 

\begin{theorem} $\cM_d$ is  a compact orbifold. It has a decomposition
$\cM_d=\M_d \cup Z_d$ into an open complex suborbifold $\M_d$ (which
is actually a complex manifold) and a codimension one suborbifold $Z_d$,
i.e. a divisor.
The singular locus of $\cM_d$ coincides with that of $Z_d$:

$$\sing(\cM_d)=\sing(Z_d)=\bigcup_{0<i\leq d} P(E^2_i)$$
where $E^2_i$ is a component of the fixed points set of the involution
$\sigma(E,\Phi)=(E,-\Phi)$.

Furthermore, the $\C^*$-action on $\M_d$ extends to $\cM_d$ with an extra
component $Z_d$ of the fixed point set.
\label{Morbi}
\end{theorem}

We have the corresponding statement of Theorem~\ref{Zjeno}.

\begin{theorem} $\cM=\cM_{g-1}$ is birationally isomorphic to
$\cM_0=P(T^*_\N\oplus{\calO}_\N)$. Moreover, they are isomorphic outside
an analytic subset of codimension at least $3$.
\label{Mjeno}
\end{theorem}

\begin{proof} In a similar manner to the proof 
of Corollary~\ref{Zjeno} we can argue
by noting that $X_{\cM_0}$ is obviously isomorphic to $T^*_\N\oplus 
\calO_\N$ with
the standard action of $\C^*$. Hence indeed 
$\M_0=P(T^*_\N\oplus{\calO}_\N)$. 

By Theorem~\ref{quotients} it is clear that $\cM$ and $\cM_0$ are related
by a sequence of blowups and blowdowns. 
The codimensions of the submanifolds
we apply the blowups are at least $3$ by a calculation analogous to the one
in the proof of Proposition~\ref{birational}.  
\end{proof}

\begin{notation} Let 

\begin{itemize} \item
$\calL_{P(T^*_\N \oplus\calO_\N)}$ denote
the pullback of 
$\calL_\N$ to $P(T^*_\N\oplus \calO_\N)$, 
\item 
$\calL_\cM$ be the corresponding 
line bundle on $\cM$,
\item
$L_{P(T^*_\N\oplus\calO_\N)}$ be the dual of the tautological line bundle
on the projective bundle $P(T^*_\N\oplus \calO_\N)$, 
\item  
$L_\cM=X_\cM\times_{\C^*}\C$ be the corresponding line orbibundle on $\cM$.  
\end{itemize}

\end{notation}

\begin{corollary} $\Pic{\cM}$ is isomorphic to
$\Pic(P(T^*_\N\oplus{\calO}_\N))$ and therefore is of rank $2$ and freely
generated by $L_\cM$ and
$\calL_\cM$.
\label{picM}
\end{corollary}

\begin{proof} 
The previous theorem shows that $\cM$ and $P(T^*_\N\oplus \calO_\N)$ are 
isomorphic outside an analytic subset of codimension at least $2$, thus their
Picard groups are naturally isomorphic. 

However, $\Pic(P(T^*_\N\oplus \calO_\N))$ is freely generated by 
$L_{P(T^*_\N\oplus\calO_\N)}$ and 
$\calL_{P(T^*_\N\oplus\calO_\N)}$. The result follows. 
\end{proof}

\begin{corollary}
The canonical line orbibundle $K_\cM$ of $\cM$ coincides with 
$L_\cM^{-(3g-2)}$. More\-over, $L_\cM$ is the line bundle of the divisor $Z$,
therefore (3g-2)Z is the anticanonical divisor of $\cM$. Finally, $L_\cM$
restricts to $L_Z$ to $Z$. 
\label{vonat}
\end{corollary}

\begin{proof} 
$L_\cM$ by its construction clearly restricts to $L_Z$ on $Z$ and
it is the line bundle of $Z$, as the corresponding statement
 is obviously true for
$P(T^*_\N\oplus \calO_\N).$

The restriction of $K_\cM$ to $\M$ has a non-zero section, namely the 
holomorphic Liouville form $\omega_h^{3g-3}$, thus trivial. Hence
$K_\cM=L_\cM^k$ for some $k\in \Z$.

By the adjunction formula $K_Z=(K_\cM\otimes [Z])\mid_Z$. The right
hand side equals $L_Z^{-(3g-3)}$ as $L_Z$ is a contact line
bundle\footnote{Cf. Theorem~\ref{contact}.}. 
The left hand side can be written as 
$(L_\cM^k\otimes L_\cM)\mid_Z=L_Z^{k+1}$, therefore $k=-(3g-2)$.
\end{proof} 

\begin{lemma} $\chi$ has an extension to $\cM$,  
$$\overline{\chi}:\cM\rightarrow \Proj^{3g-3}$$ such that $\overline{\chi}$ restricted
to $Z$ gives the map of Lemma~\ref{Zchi}.
\label{Mchi}
\end{lemma}

\begin{proof}
We let $\C^*$ act on $\C^{3g-3}\times \C$ by $\lambda(x,z)=
(\lambda^2 x,\lambda z)$. With respect to this action the map 
$(\chi,id_\C):\M\times \C\rightarrow \C^{3g-3}\times \C$ is equivariant. 
Therefore making the symplectic cut 
it reduces to a map $\overline{\chi}:\cM\rightarrow \Proj^{3g-3}$ since the
quotient space 
$(\C^{3g-3}\setminus 0)\times\C /\C^*$ is isomorphic to $\Proj^{3g-3}$. 

The result follows. 
\end{proof}

\begin{remark} In the higher rank case, where $\C^*$ acts on the
target space of the Hitchin map with different weights, the target space
of the compactified Hitchin map is a weighted projective space.
\end{remark}

\begin{corollary} $L^2_\cM=\overline{\chi}^*{\cal H}_{3g-3}$.
\label{hurka}
\end{corollary}

\begin{proof} Obviously, $\overline{\chi}^*{\cal H}_{3g-3}\mid_\M$ is trivial,
therefore $\overline{\chi}^*{\cal H}_{3g-3}$ is some power of $L_\cM$. 
By \ref{feco} this power is $2$. 
\end{proof}

\begin{theorem}[Duistermaat, Heckman] $\cM$ has a one-para\-meter family
of K\"ahler forms $\omega_t(\cM)$, $t>c_{g-1}$ such that

$$[\omega_{t_1}(\cM)]-[\omega_{t_2}(\cM)]=(t_1-t_2)c_1(L_{\cM}).$$

Furthermore this one-parameter family of K\"ahler forms restricts to Z as
the one-parameter family of K\"ahler forms of Theorem~\ref{csalad}.
\label{DuHeM}
\end{theorem}

\begin{proof} This is just the application of Theorem~\ref{quotients} and
Theorem~\ref{cut} to our situation.
\end{proof}

\begin{corollary} $\cM$ is a projective algebraic variety.
\label{Mproj}
\end{corollary}

\begin{proof} The argument is the same as for Theorem~\ref{proj}, noting that
by Corollary~\ref{picM} $H^2(\cM,\R)$ is two dimensional and $L_\cM$ is neither
trivial nor ample since $L_\cM\mid Z=L_Z$ (by Corollary~\ref{vonat}) is neither
trivial nor ample (by Corollary~\ref{triv}).
\end{proof}

\begin{remark} 1. The above proof yields that the cohomology class 
$c_1(L_\cM)$ sits
in the closure of the K\"ahler cone of $\cM$, hence $L_\cM$ is nef.

2.  From the previous remark and Corollary~\ref{hurka} 
we can deduce that there is a complete hyperk\"ahler (hence Ricci flat) 
metric on $\M=\cM\setminus Z$, the complement of a nef anticanonical divisor 
of a compact orbifold. 

Therefore our compactification of 
$\M$ is compatible with Yau's problem, which addresses the question: which
non-compact complex manifolds possess a complete Ricci flat metric? 
Tian and Yau in \cite{yau-tian} 
could show that this is the case for the complement of
an ample  anticanonical divisor in a compact complex manifold. (Such manifolds
are called Fano manifolds.) 

The similar statement with ample replaced by nef is an 
unsolved problem.
\end{remark}

\begin{theorem} $\cM$ has Poincar\'e polynomial
$$P_t(\cM)=P_t(\M)+t^2P_t(Z).$$
\label{PM}
\end{theorem}

\begin{proof} We have three different ways of calculating the Poincar\'e
polynomial of $\cM$. The first is through Kirwan's formula in \cite{kirwan1},
the second is due to Thaddeus in \cite{thaddeus3}, 
which we used to calculate the Poincar\'e
polynomial of $Z$.

For $\cM$ there is a third method, namely direct Morse theory. All we have to
note is that the $U(1)$-action $\cM$ is Hamiltonian with respect to
any K\"ahler form of Theorem~\ref{DuHeM}, and the critical 
submanifolds and corresponding
indices are the same as for $\M$ with one extra critical submanifold $Z$
of index $2$. Hence the result.
\end{proof}

\begin{example} We can describe $\cM_{toy}=\M_{toy}\cup Z_{toy}$ as follows.
As we saw above $\M_{toy}\setminus N_{toy}=X_{Z_{toy}}$. Thus gluing together
$\M_{toy}$ and $E_{toy}$, the total space of the line orbibundle 
$L_{Z_{toy}}$, along $X_{Z_{toy}}$ yields 
$$\cM_{toy}=\M_{toy}\cup_{X_{Z_{toy}}} E_{toy}.$$

One can construct $\cM_{toy}$ directly, as follows. Take 
$\Proj^1=\C\cup \infty$ extending  the involution $\tau$ from $\C$ to 
$\Proj^1$. 
Consider
the quotient $(P\times \Proj^1)/(\sigma_P\times \tau)$. This is a compact
orbifold with eight $\Z_2$-quotient singularities. Blow up four of them
corresponding to $0\in \C$. The resulting space will be isomorphic
to $\cM_{toy}$. The remaining four isolated $\Z_2$ quotient singularities will
just be the four marked points of $Z_{toy}\subset \cM_{toy}$, the singular
locus of $\cM_{toy}$. 
\end{example}

We finish this section with a result which gives an interesting relation
between the intersections of the component of the nilpotent cone $N$
in $\M$ 
(equivalently the intersection form on the middle compact cohomology 
$H_{cpt}^{6g-6}(\M)$ from  Corollary~\ref{middle})  
and the contact structure of $Z$.

\begin{theorem}
There is a canonical isomorphism between the cokernel of $j_\M$ and the
cokernel of $L$, where
$$j_\M:H_{cpt}^{6g-6}(\M)\rightarrow H^{6g-6}(\M)$$ is the
canonical map and
$$L:H^{6g-8}(Z)\rightarrow H^{6g-6}(Z)$$
is multiplication with $c_1(L_Z)$.
\label{fura}
\end{theorem}

\begin{proof} We will read off the statement from the following diagram. 

$$
\begin{array}{ccccccccc}
 &  &     &       &       0           &        &         &        & \\
 &  &     &       &      \downarrow   &        &         &        & \\
   &    &     &   & H^{6g-8}(Z) &        &         &        & \\ 
   &    & & &\downarrow     &{\searrow}^L &         &  & \\

0 & \rightarrow & H^{6g-6}_{cpt}(\M)&\rightarrow & H^{6g-6}(\cM)
&\rightarrow &H^{6g-6}(Z)&\rightarrow & 0 \\
 & &  &{\searrow}^j&\downarrow&   &  &  &\\
& & & & H^{6g-6}(\M)&&&&\\  
&&&&\downarrow&&&& \\
&&&&0&&&&\\
\end{array}
$$

We show that both the vertical and horizontal sequences are exact and the
two triangles commute. 

From  the Bialynicki-Birula decomposition of $\cM$ 
we get the short exact sequence of
middle dimensional cohomology groups\footnote{Recall that $E\subset\cM$ denotes
the total space of the contact line bundle $L_Z$ on $Z$.}:
$$0\rightarrow H^{6g-6}_{cpt}(E)\rightarrow H^{6g-6}(\cM)\rightarrow
H^{6g-6}(\M)\rightarrow 0.$$

Applying the Thom isomorphism\footnote{This also holds in the orbifold
category with coefficients from $\Qu$.} 
we can identify $H^{6g-6}_{cpt}(E)$ with $H^{6g-8}(Z)$, this gives
the vertical short exact sequence of the diagram. The horizontal one
is just its dual short exact sequence. 

Finally, the left triangle clearly commutes as all the maps are natural, 
while the right triangle commutes because the original triangle commuted as
above and the canonical map $j_E:H^{6g-6}_{cpt}(E)\rightarrow H^{6g-6}(E)$ 
transforms to $L:H^{6g-8}(Z)\rightarrow H^{6g-6}(Z)$ by
the Thom isomorphism. 

Now the theorem is the consequence of the Butterfly
lemma\footnote{Cf. \cite{lang} IV.$4$ p.102.}, 
or can be proved
by an easy diagram chasing.

Hence the result follows. 
\end{proof}

\begin{remark} 

1. In the next Chapter we shall prove Theorem~\ref{main}, that $j_\M$
is $0$. Combined with the above theorem we have that
the cokernel of $L$ is $g$-dimensional! We will say more about this in
Section~\ref{compactthesis}.

2. If the line bundle $L_Z$ was ample, the map $L$ would just be the
Lefschetz isomorphism, and therefore the cokernel would be trivial. In
our case we have $L_Z$ being only nef and the map is certainly not an isomorphism,
the cohomology class of the Kummer variety lies in the kernel. Therefore
the cokernel measures how far $L_Z$ is from being ample. 

3. Notice that the proof did not use any particular property of 
$\M$, therefore the statement is true in the general setting of 
Section~\ref{rizsa}.
\end{remark}

\begin{example} 
We can simply calculate the dimension of the cokernels in our toy example.
Namely, the dimension of $\coker(L_{toy})$ is clearly $1$, as the map 
$L_{toy}:H^0(Z_{toy})\rightarrow H^2(Z_{toy})$ is the 
multiplication\footnote{Cf. the example at the end of Section~\ref{kaehler}.} with
$c_1(L_{Z_{toy}})=0$. 

Thus, by the above theorem, we have that $\coker(j_{\M_{toy}})$ is 
$1$-dimensional. It can be seen directly, using Zariski's
lemma\footnote{Lemma 8.2 in \cite{barth-peters-van} 
p. 90.}, that the kernel of the
map $j_{\M_{toy}}$ is generated by the cohomology class of the 
elliptic curve $P$, the generic fibre of the toy Hitchin map, hence it is
$1$-dimensional, indeed.         

Thus $j_{\M_{toy}}$ is not $0$ unlike $j_\M$. This indicates a profound
difference between $\M$ and $\M_{toy}$.  

\end{example}

\newpage
\thispagestyle{empty}

\chapter{Intersection numbers}
\label{intersection}

The aim of the present chapter is to prove Theorem~\ref{main}. We
discussed in the Introduction the physical motivation
for Theorem~\ref{main}. 

From an algebraic geometric point of view 
Theorem~\ref{main} can be interpreted as follows. First
of all it is really about middle dimensional cohomology,
because we know that $\M$ does not have cohomology beyond the
middle dimension, and equivalently by Poincar\'e duality 
$\M$ does not have compactly supported cohomology below the middle dimension. 
Thus the main content of Theorem~\ref{main} is the vanishing of
the canonical map $j_\M:H_{cpt}^{6g-6}(\M)\rightarrow H^{6g-6}(\M)$
between $g$-dimensional spaces\footnote{Cf. Corollary~\ref{middle}.}. 
This in turn is
equivalent to the vanishing of the intersection form on
$H_{cpt}^{6g-6}(\M)$. 

There are $g+1$ intersection numbers whose 
vanishing follows easily.
 One vanishing is obtained by recalling that the moduli space $\N$ of stable
bundles of real dimension $6g-6$ 
sits inside $\M$ with normal bundle $T^*_\N$, thus its
self-intersection number is its Euler characteristic up to sign, which is
known\footnote{Substitute $t=-1$ into (\ref{poincaren})!} to vanish. 
 
The other $g$ vanishings follow from the fact that the ordinary
cohomology class of the Prym variety, the generic fibre of the Hitchin
map, is $0$ i.e. $j_\M(\overline{\eta}_P^\M)=0$. This can be seen
by thinking of the Hitchin map as a section of the trivial
rank $3g-3$ vector bundle on $\M$ and considering the ordinary cohomology class
of the Prym variety as the Euler class of this trivial vector bundle, and as
such, the ordinary cohomology class of the Prym variety is trivial
indeed.   
Note that for the case $g=2$, the above vanishings are already enough
to have $j_\M=0$. 
 
The vanishing of the
intersection form on $\M$ for any genus, 
proved in this chapter, can be considered as a generalization of these facts.  
We should also mention that as we explained at the end of the previous
chapter, the intersection form is not trivial in the case of the toy
example $\M_{toy}$. 

The structure of this section is as follows: In the next section we
develop the theory of stable Higgs bundles analogously to the stable
case, and prove an important vanishing theorem. In
Section~\ref{universalbundles}
we prove that $\tM$ is a fine moduli space, and define certain
universal bundles. In Section~\ref{virtualdirac} we construct the
virtual Dirac bundle\footnote{For its gauge theoretic construction
see Subsection~\ref{diracequ}.}, as the analogue of the virtual
Mumford bundle,  and show that it can be considered as the degeneracy
sheaf of a homomorphism of vector bundles. In
Section~\ref{downdegeneracylocus} we determine the degeneracy locus of the
above homomorphism in terms of the components of the nilpotent
cone. Finally in Section~\ref{proof} we prove our main
Theorem~\ref{main} using Porteous' formula for the degeneracy locus of
the virtual Dirac bundle.

\newpage
\section{A vanishing theorem}

\label{vanthe}

\begin{definition} 
The complex $\compE$ with $E$ a
vector bundle on $\Sigma$, $K$ the canonical bundle of $\Sigma$, 
and $\Phi\in H^0(\Sigma,\Hom(E,E\otimes K))$, 
is called a {\em Higgs bundle}, while $\Phi$ is called the 
{\em Higgs field}. 

A morphism $\Psi:{\cal E}_1 \rightarrow {\cal E}_2$ 
between two Higgs bundles 
\mbox{${\cal E}_1=E_1\stackrel{\Phi_1}{\rightarrow} E_1\otimes K$} and  
\mbox{${\cal E}_2=E_2\stackrel{\Phi_2}{\rightarrow} E_2\otimes K$} 
is defined to be a homomorphism of vector bundles 
$\Psi\in \Hom(E_1,E_2)$ such that the following diagram commutes:
$$
\begin{array}{c}
 \  E_1\stackrel{\Phi_1}{\longrightarrow} E_1\otimes K \ \\
\hskip.35cm\mbox{\scriptsize{$\Psi$}}\downarrow\hskip1.8cm \downarrow 
\mbox{\scriptsize{$\Psi\otimes id_K$}} \\
 \  E_2\stackrel{\Phi_2}{\longrightarrow} E_2\otimes K\ 
\end{array}
$$ 

Moreover we say that 
${\cal E}_1$ is a {\em Higgs
subbundle} of ${\cal E}_2$ if 
$\Psi\in \Hom(E_1,E_2)$ is injective and a morphism of Higgs
bundles. We denote this by ${\cal E}_1\subset {\cal E}_2$. 
In this case we can easily construct  
the quotient Higgs bundle ${\cal E}_2/{\cal E}_1$ together with a
surjective morphism of Higgs bundles $\pi:{\cal E}_2\rightarrow {\cal
E}_2/{\cal E}_1$ whose kernel is exactly ${\cal E}_1$.
\label{higgs}
\end{definition}

\begin{remark} It is a tautology that morphisms of Higgs bundles
form the hypercohomology\footnote{In connection
with Higgs bundles the language of 
hypercohomology was first used in \cite{simpson}.
In \cite{biswas-ramanan} 
it was used to describe the tangent space to $\M$.} 
vector space  
$$\Hy^0(\Sigma, E_1^*\otimes E_2\stackrel{\left[\Phi_1,\Phi_2\right]}
{\longrightarrow}
E_1^*\otimes E_2\otimes K),$$ where the homomorphism 
$\left[\Phi_1,\Phi_2\right]$ is given
by: $${\left[\Phi_1,\Phi_2\right](\Psi):=(\Psi\otimes
id_K)\Phi_1-\Phi_2\Psi}$$ for $\Psi\in \Hom(E_1,E_2).$ 
\end{remark}

Now recall Definition~\ref{stablehiggs} where we defined the notion  
of stability of Higgs bundles. 
The main result\footnote{The second part of which is Proposition
(3.15) in \cite{hitchin1}.} of this section is the following theorem about
morphisms between stable Higgs bundles.

\begin{theorem} Let ${\cal E}=E\stackrel{\Phi}{\rightarrow}E\otimes K$
and ${\cal F}=F\stackrel{\Psi}{\rightarrow}{F\otimes K}$ be stable Higgs
bundles with $\mu({\cal F})<\mu({\cal E})$. Then the only morphism
from ${\cal E}$ to ${\cal F}$ is the trivial one. In other words 
$$\Hy^0(\Sigma, E^*\otimes F\stackrel{[\Phi,\Psi]}{\longrightarrow}
E^*\otimes F\otimes K)=0.$$

Moreover if $\mu({\cal F})=\mu({\cal E})$, then there is a non-trivial
morphism $f:{\cal E} \rightarrow {\cal F}$ if and only if 
${\cal E}\cong {\cal F}$ in which case every non-trivial morphism 
$ f:{\cal E} \rightarrow {\cal F}$ is an isomorphism and 
\begin{eqnarray}
\dim(\Hy^0(\Sigma, E^*\otimes F\stackrel{[\Phi,\Psi]}{\longrightarrow}
E^*\otimes F\otimes K))=1.
\label{iso}
\end{eqnarray}
\label{van}
\end{theorem}

\begin{proof}

A characteristic property of the ring $\C[x_1]$ is that it is the only
principal ideal domain among the rings $\C[x_1,\dots,x_n]$. 
It follows that $\calO_\Sigma$ is a sheaf of principal ideal domains,
and that every torsion free sheaf\footnote{Every sheaf over $\Sigma$ is assumed to
be a sheaf of $\calO_\Sigma$-modules.} over $\Sigma$, 
such as a subsheaf of a locally free sheaf,
will be locally free. An easy consequence of this is  
the lemma of Narasimhan and Seshadri\footnote{Cf. section 4 in \cite{narasimhan-seshadri}.}:

\begin{lemma}
Let $E$ and $F$ be two vector bundles over the Riemann surface
$\Sigma$ with a non-zero homomorphism $f:E\rightarrow F$, then $f$ has
the following canonical factorisation: 
$$
\begin{array}{c}
 0\longrightarrow E_1 \longrightarrow E
 \stackrel{\eta}{\longrightarrow} E_2 \longrightarrow 0 \\
\hskip1.65cm \downarrow {\mbox{\scriptsize{f}}} \hskip.75cm 
\downarrow {\mbox{\scriptsize{g}}} \\
0 \longleftarrow F_2 \longleftarrow F \stackrel{i}{\longleftarrow} F_1
 \longleftarrow 0
\end{array}
$$
where $E_1,E_2,F_1$ and $F_2$ are vector bundles, each row is exact, 
$f=ig\eta$ and $g$ is of maximal rank, i.e. $\rank(E_2)=\rank(F_1)=n$ and 
$\Lambda^n(g):\Lambda^n(E_2)\rightarrow \Lambda^n(F_1)$ is a non-zero 
homomorphism. In other words $g$ is an isomorphism on a Zariski open
subset
$U$ of $\Sigma$. $F_1$ is called the subbundle of $F$ generated by the
image of $f$. $\square$
\label{narasimhan-seshadri}
\end{lemma}

Let $f:{\cal E}\rightarrow {\cal F}$ be a non-zero morphism of Higgs
bundles. In particular $f:E\rightarrow F$ is a homomorphism of
vector bundles. 

Construct the canonical factorisation of $f$ of
Lemma~\ref{narasimhan-seshadri}.  Consider the Zariski open subset $U$ of $\Sigma$ where
$g$ is an isomorphism. Here clearly $\ker(f\mid_U)=\ker(\eta\mid_U)
= E_1\mid_U$. Now $\ker(f\mid_U)$ being the kernel of a morphism of
Higgs bundles, is $\Phi$-invariant, i.e. a Higgs subbundle of ${\cal E}\mid_U$. 
Thus $E_1\mid_U$ is a Higgs subbundle of ${\cal E}\mid_U$. This means
that $\Phi(E_1)$ is contained in $E_1\otimes K\subset E\otimes K$ on
$U$. Because $U$ is Zariski open in $\Sigma$, it follows that 
${\cal E}_1=E_1\stackrel{\Phi}{\rightarrow} E_1\otimes K$ is
a Higgs subbundle of ${\cal E}$. Let ${\cal E}_2=E_2
\stackrel{\Phi}{\rightarrow} E_2\otimes K$ denote the quotient
Higgs bundle.  

Similarly $\im(\alpha)\mid_U= F_1\mid_U$ is $\Psi$-invariant, thus
${\cal F}_1=F_1\stackrel{\Psi}{\rightarrow} F_1\otimes K$ is a
Higgs subbundle of $F$. 

By assumption $\mu({\cal F)}<\mu({\cal E})$, 
stability of ${\cal E}$ gives $\mu({\cal E})\leq \mu ({\cal E}_2)$ (it
may happen that $E=E_2$) and
because $g$ is of maximal rank we get 
$\mu({\cal E}_2)=\mu(E_2)\leq \mu(F_1)=\mu({\cal F}_1)$.
Thus $\mu({\cal F})< \mu({\cal F}_1)$ contradicting the stability of
$\cal F$. 

If $\mu(\cE)=\mu({\cal F})$ then the above argument leaves the only
possibility that $\eta$, $g$ and $i$ are isomorphisms, showing that
$f$ must be an isomorphism. Suppose that we have such an isomorphism
$f$ of Higgs bundles. Then consider $h:\cE \rightarrow {\cal F}$ 
another non-zero morphism of Higgs bundles. In particular $h\in 
\Hom(E,F)$. Let $\lambda$ be an eigenvalue of the homomorphism 
$f_p^{-1}h_p\in\Hom(E_p,E_p)$. Then the homomorphism $h-\lambda f$ is
not an isomorphism, though clearly a morphism of Higgs bundles. From
the above argument  $h-\lambda f=0$.     

The result follows. 
\end{proof}

\begin{corollary} 
For any stable Higgs bundle ${\cal E}$ with $\mu({\cal E})<0$:
 \begin{eqnarray}\Hy^0(\Sigma,{\cal E})=0,\label{h0}\end{eqnarray}
for any stable Higgs bundle ${\cal E}$ with $\mu({\cal E})>0$:
\begin{eqnarray} \Hy^2(\Sigma,{\cal E})=0.\label{h2} \end{eqnarray}
If $\cE$ is a stable Higgs bundle with $\mu(\cE)=0$ and $\cE \ncong
\cE_0=\calO_{\Sigma}\stackrel{0}{\rightarrow}\calO_\Sigma \otimes K$
then both (\ref{h0}) and (\ref{h2}) hold. 
\label{vanish}
\end{corollary}
\begin{proof} 
For the first part consider the Higgs bundle 
${\cal E}_0=\calO_{\Sigma}\stackrel{0}{\rightarrow}\calO_\Sigma \otimes
K$. 
Being of rank $1$ it is obviously stable, with $\mu({\cal
E}_0)=0$. Now the previous theorem yields that there are no nontrivial
morphisms
from ${\cal E}_0$ to ${\cal E}$, which in the language of
hypercohomology is exactly $\Hy^0(\Sigma,{\cal E})=0$, which we had to
prove. 

For the second part Serre duality gives that 
$\Hy^2(\Sigma,\cE)\cong (\Hy^0(\Sigma,\cE^*\otimes K))^*$. Now 
 clearly $\cE^*\otimes K$ is stable and 
$\mu(\cE^*\otimes K)=-\mu(\cE)<0$. Thus the first part gives the
second. 

Likewise, the third statement follows by referring to the last part of 
Theorem~\ref{van}.
\end{proof}

It is convenient to include here the analogue of the 
Harder-Narasimhan filtration for Higgs bundles, which will be used in Chapter~\ref{cohomology}:

\begin{corollary}
 Every Higgs bundle $\cE$ has
a {\em canonical filtration}: \begin{eqnarray}
0=\cE_0\subset \cE_1 \subset \cE_2 \subset \dots \subset
\cE_r=\cE,\label{higgsfiltration} \end{eqnarray}
with $\cD_i=\cE_i/\cE_{i-1}$ semi-stable and $$\mu(\cD_1)>\mu(\cD_2)>\dots >
\mu(\cD_r).$$
\label{higgsfiltrationtheorem}
\end{corollary}

\begin{proof} It follows from Theorem~\ref{van} just as in the case of
vector bundles. \end{proof}

\newpage
\section{Universal bundles}
\label{universalbundles}

Nitsure showed that $\tM$ is a coarse moduli space. 
Here we show that 
$\tM$ is in fact a {\em fine} moduli space. We closely follow
the proof of Theorem 5.12 in \cite{newstead} and (1.19) of
\cite{thaddeus3}. All the ingredients have already appeared in the 
unpublished \cite{thaddeus1}. 

\begin{definition} Two families $\cE_T$ and $\cE_T^\prime$ 
of stable Higgs bundles over $T\times \Sigma$ are said to be 
{\em equivalent}, (in symbols $\cE_T\sim \cE_T^\prime$) 
if there exists a line
bundle $L$ on $T$ such that 
 $\cE^\prime_T\cong\cE_T\otimes\pr_T^*(L)$. 
\end{definition}

The next lemma, which is taken from \cite{thaddeus1}, 
shows that two families are equivalent iff they give rise to the same
map to the coarse moduli space $\tM$.

\begin{lemma} If 
$\cE_T=\E_T\stackrel{\uPhi}{\rightarrow}\E_T\otimes
\K$ and $\cE_T^\prime=
\E^\prime_{T}\stackrel{\uPhi^\prime}{\rightarrow}
\E^\prime_T\otimes
\K$ are families of stable
Higgs bundles
over $T\times \Sigma$ such 
that \begin{eqnarray}\cE_T\mid_{\{t\}\times \Sigma}\cong
\cE^\prime_T\mid_{\{t\}\times \Sigma}\label{con}
\end{eqnarray} for each $t\in T$, then $\cE_T\sim \cE_T^\prime$.
\label{ramanan}
\end{lemma}
  
\begin{proof} Let $${\cal F}:= \E^*_T\otimes\E_T^\prime
\stackrel{[\uPhi_T,\uPhi_T^\prime]}{\longrightarrow} 
\E^*_T\otimes\E_T^\prime\otimes\K.$$      
We define $L=\R^0{\pr_T}_*({\cal F})$. By (\ref{con}) and 
(\ref{iso}) this is a line
bundle over $T$. It follows from the projection formula\footnote{Cf.  
p.124 of \cite{hartshorne}.} that 
the sheaf $\R^0{\pr_T}_*({\cal F}\otimes 
\pr_T^*(L^*))$ is just ${\calO}_T$, the structure sheaf. A
non-zero section $${\mathbf \Psi}\in H^0(T,\R^0{\pr_T}_*({\cal F}\otimes 
\pr_T^*(L^*)))$$ for every $t\in T$ gives $${\mathbf
\Psi}\mid_{\{t\}\times \Sigma}:
(\cE_T\otimes \pr_T^*(L))\mid_{\{t\}\times \Sigma} \rightarrow 
\cE^\prime_T\mid_{\{t\}\times \Sigma}$$ a non-zero morphism of  
Higgs bundles, which is by Theorem~\ref{van} an isomorphism. 

The result follows.
\end{proof}

Now we prove the existence of universal Higgs bundles\footnote{Cf. \cite{thaddeus1}.}:

\begin{proposition} Universal Higgs bundles $\cE_\tM=\ucompE$ over
$\tM\times \Sigma$ do exist. 
\label{universal}
\end{proposition}

\begin{proof} We construct a {\em holomorphic} universal Higgs bundle over
$\tM\times \Sigma$ by using the gauge theoretic construction of $\tM$ from
Subsection~\ref{gaugemk}. An analogous construction in the
GIT framework of \cite{nitsure} however gives an {\em algebraic} universal
Higgs bundle. We preferred here the gauge theoretic
proof, because it fits better into this thesis. 

To construct $\E_\tM$ we proceed similarly to p.579-580 of
\cite{atiyah-bott}. First we note that there is an obvious tautological rank
$2$ holomorphic
bundle $\E_\calC$ over $\calC\times \Sigma$ with the constant scalars $\C^*\subset
\G^c$ acting trivially on the base and as scalars in the fibre of
$\E_\calC$. 

We also need a $\G^c$-invariant holomorphic line bundle $\bL$ over
$\calC_{< g}$ such that $\C^*\subset \G^c$ acts via scalar
multiplication. To construct such a line bundle  we choose $k$ large enough such that
$H^1(M,E\otimes L_p^k)=0$ for each $E\in \calC_{<g}$. Then
$(\pr_{\calC})_*(\E_\calC\otimes L_p^k)$ is a $\G^c$-equivariant 
holomorphic vector bundle over
$\calC_{<g}$ of degree $2k+1-2(g-1)$. Taking determinants gives a
$\G^c$-equivariant line bundle $A_k$ on $\calC_{<g}$, with the group $\C^*$ of scalar
automorphisms of $E$ acting on this with weight $m=2k+1-2(g-1)$. Now
taking the determinant of $ (\pr_{\calC})_*(\E_\calC\otimes L_p^{k+1})$ gives a
line bundle $A_{k+1}$ over $\calC_{<g}$ such that  the weight of the 
$\C^*$-action is $2k+3-2(g-1)$. It follows that we can find $a$ and $b$ such
that $\bL=A^a_k\otimes A^b_{k+1}$ is a holomorphic $\G^c$-equivariant 
line bundle over $\calC_{<g}$, with $\C^*\subset \G^c$ acting on it
with weight $1$. 

From Subsection~\ref{gaugemk}  we have $(\calB)^s\subset \calB$
the subspace of stable Higgs bundles and a map $\pr:\calB\to
\calC$. The restricted map $\pr^s:(\calB)^s\to \calC$ has image in
$\calC_{<g}$. Thus we can pullback the bundle $\E_\calC\otimes \bL^{-1}$ from
$\calC_{<g}\times \Sigma$ to $(\calB)^s\times \Sigma$ gaining
$\E_\calB$. This bundle is a priori a 
$\G^c$-equivariant bundle. However the subgroup of constant scalar
automorphisms $\C^*\subset \G^c$ acts trivially on $\E_\calB$, therefore it
reduces to a $\bG^c=\G^c/\C^*$-equivariant bundle. Since $\bG^c$
acts freely on $(\calB)^s$ it follows that $\E_\calC\otimes \bL^{-1}$ reduces to
a rank $2$ holomorphic vector bundle $\E_{\tM}$ over $\tM\otimes \Sigma$, with the
property that $\E_{\tM}\mid_{(E,\Phi)}\cong E$. 

Now the fibre of $\pr:\calB\to \calC$ over a point $E$, can be
identified canonically with $H^0(\Sigma,\End(E)\otimes K)$, thus there
is an obvious tautological section $$\uPhi_\calB\in
H^0(\calB;\End(\pr^*(\E_\calC))\otimes K).$$ However
$\End(\pr^*(\E_\calC)\cong \End(\E_\calB)$ canonically. It follows
that we have a tautological Higgs bundle
$\E_\calB\stackrel{\uPhi_\calB}{\longrightarrow}\E_\calB\otimes K$ over $(\calB)^s\times
\Sigma$, which is a priori a $\G^c$-equivariant complex but reduces to
a $\bG^c$-equivariant complex as proved above. It follows that it reduces to $\tM\times
\Sigma$ to give a universal Higgs bundle
$\E_{\tM}\stackrel{\uPhi_{\tM}}{\longrightarrow}\E_{\tM}\otimes K$.
\end{proof}

As in Theorem 5.12 of \cite{newstead} and (1.19) of \cite{thaddeus3} 
our Lemma~\ref{ramanan} and Proposition~\ref{universal} gives:

\begin{corollary} The space $\tM$ is a fine moduli space for 
rank $2$ stable Higgs bundles of degree $1$ with respect to the
equivalence $\sim$ of families of stable Higgs bundles.   
\end{corollary}

As another  consequence of Proposition~\ref{universal} and 
Lemma~\ref{ramanan}, we see that although $\E_\tM$ is not unique
$\End(\E_\tM)$ is.  Moreover it is clear that 
by setting $\E_\M=\E_\tM\mid_{\M\times \Sigma}$ we have
\begin{eqnarray}
c(\End(\E_\tM))=c(\End(\E_\M))\otimes 1\label{c}
\end{eqnarray}
in the decomposition (\ref{splitM}).
 
Thus from the K\"unneth decomposition of $\End(\E_\M)$ we get 
universal classes
\begin{eqnarray} c_2(\End(\E_\M))=2\alpha_\M\otimes \sisi+\sum_{i=1}^{2g}4\psi^i_\M 
\otimes \xisi_i -
\beta_\M\otimes 1\label{universalclassesM} \end{eqnarray}
in $H^4(\M\times \Sigma)\cong
\sum_{r=0}^{4} H^r(\M)\otimes H^{4-r}(\Sigma)$
for some $\alpha_\M\in H^2(\M)$, $\psi^i_\M\in H^3(\M)$ and $\beta_\M\in
H^4(\M)$.

Though $\E_\M$ is not unique we can still write its Chern classes 
in the K\"unneth decomposition\footnote{Cf. proof of Newstead's theorem in \cite{thaddeus2}).}, 
getting\footnote{Note that $\M$ being simply connected by
\cite{hitchin1} $H^1(\M)=0$.} 
$$c_1(\E_\M)=1\otimes \sisi +\beta_1\otimes 1,$$
where
$\beta_1\in H^2(\M)$
and $$c_2(\E_\M)=\alpha_2 \otimes \sisi + \sum_{i=1}^{2g}a_i \otimes
\xisi_i + 
\beta_2\otimes 1,$$
where $\alpha_2\in H^2(\M)$, $a_i\in H^3(\M)$ and $\beta_2\in
H^4(\M)$. 
Since
$$4c_2(\E_\M)-c_1^2(\E_\M)=c_2(\End(\E_\M)),$$ we get
$\alpha_\M=2\alpha_2-\beta_1$ and $\beta=\beta_1^2-4\beta_2$.
Because $\Pic(\M)\cong H^2(\M,\Z)$ we can
normalize $\E_\M$ uniquely such that $c_1((\E_\M)_p)=\alpha_\M$, where
$$(\E_\M)_p:=\E_\M\mid_{\M\times \{p\}}.$$  

\begin{definition} The universal Higgs bundle $\cE_\M$ is
normalized if $c_1((\E_\M)_p)=\alpha_\M$.
\end{definition}

We also need to work out the Chern classes of $\E_\tM$.
It is easy to
see that $c(\E_\tM)$ in the decomposition (\ref{splitM}) is the product
of $c(\E_\tM)\mid_{\M\times \Sigma}$ and $c({\bL}_1)$, where 
${\bL}_1$ is some universal line bundle over $\J\times \Sigma$. 

\begin{definition} We call the universal
Higgs bundle $\cE_\tM$ {\em normalized} if in the decomposition 
(\ref{splitM}) 
\begin{eqnarray} 
c_1((\E_\tM)_p)=\alpha_\M \label{c1},
\end{eqnarray}
where $(\E_\tM)_p=\E_\tM\mid_{\tM\times \{p\}}.$
\end{definition}

\begin{remark} Since
$$4c_2((\E_\tM)_p)-c_1((\E_\tM)_p)^2=c_2\left(\End((\E_\tM)_p)\right),$$
for a normalized universal Higgs bundle over $\tM\times \Sigma$
(\ref{c}) and (\ref{c1}) yield:
\begin{eqnarray}c_2((\E_\tM)_p)=\frac{(\alpha_\M^2-\beta_\M)}{4}
  \label{c2}
\end{eqnarray} 
\end{remark}

Finally, given a universal Higgs bundle $\cE_\tM$ over $\tM\times
\Sigma$, we introduce
a {\em universal Higgs bundle of degree $2k-1$} by setting 
$\cE_{\tM}^k:=\cE_{\tM}\otimes \pr_\Sigma^*(L_p^{k-1})$, where $L_p$ is
the line bundle of the divisor of the point 
$p\in \Sigma$. It is called normalized
if $\cE_\tM$ is normalized.

\newpage
\section{The virtual Dirac bundle $\D_k$}
\label{virtualdirac}

The strategy of the proof of Theorem~\ref{main} will be to examine the
virtual Dirac bundle  $\D_k$ which is defined in the following:

\begin{definition} The {\em virtual Dirac bundle} is\footnote{Recall
the definition of the pushforward of a complex from Subsection~\ref{hypercohomology}.}  
$$\D_k:=-{\pr_\tM}_!(\cE_{\tM}^k)\in K(\tM),$$ 
where $\cE^k_\tM$ is a normalized
universal Higgs bundle of degree $2k-1$ 
and $\pr_\tM:\tM\times \Sigma\rightarrow \tM$ is
the projection to $\tM$.
\end{definition}

The name is justified by Hitchin's construction\footnote{Cf. Subsection~\ref{gaugedirac}.}
 of
$\D_k$ related to the space of solutions of an equation on $\Sigma$, 
which is locally 
the dimensional reduction of the Dirac equation in $\R^4$ coupled to
a self-dual Yang-Mills field. 

The
virtual Dirac bundle is  
a priori $$-{\pr_\M}_!(\cE_{\tM}^k)=
-\R^0{\pr_\tM}_*(\cE_{\tM}^k)+
\R^1{\pr_\tM}_*(\cE_{\tM}^k)-
\R^2{\pr_\tM}_*(\cE_{\tM}^k)\in K(\tM)$$ a
formal sum of three coherent sheaves. Corollary~\ref{vanish} shows
that one of these sheaves always vanishes: 
if $k> 0$, then $\R^2=0$, if $k\leq 0$ then $\R^0=0$. From now on $k$ is
assumed to be positive.

We would like to use Porteous' Theorem~\ref{porteous} for $\D_k$. As
we 
explained in Subsection~\ref{porteoussection} for this we need to 
show that we
can think of the virtual Dirac bundle as the virtual 
degeneracy sheaf of a homomorphism
of vector bundles. More precisely we prove:

\begin{theorem}
There exist two vector bundles $V$ and $W$ over $\tM$ together with a 
homomorphism $f:V\rightarrow W$ of vector bundles, 
whose kernel and cokernel are respectively  
$\R^0{\pr_\M}_*(\cE_{\tM}^k)$ and 
$\R^1{\pr_\M}_*(\cE_{\tM}^k)$. In other
words there is an exact sequence of sheaves:
$$0\rightarrow \R^0{\pr_\tM}_*(\cE_{\tM}^k)
\rightarrow V\stackrel{f}{\rightarrow} W \rightarrow
\R^1{\pr_\tM}_*(\cE_{\tM}^k)\rightarrow
0.$$
\label{deg}
\end{theorem}

\noindent {\it Proof\footnote{The idea of the proof was suggested by Manfred Lehn.}.}
First we need a lemma. 
\begin{lemma} Let $X$ be a smooth quasi-projective variety 
and $\Sigma$ a
smooth projective curve. 
If $E$ is a locally free sheaf over $X\times {\Sigma}$ 
then there exists a vector bundle $F$ over $X \times \Sigma$ 
with a surjective vector
bundle homomorphism $g_E:F\rightarrow E$ such that 
$R^0{\pr_{X}}_*(F)=0$. 
We will call $F$ a {\em sectionless resolution} of $E$.
\end{lemma}
\begin{proof} The lemma is a special case of Proposition 2.1.10 of 
\cite{huybrechts-lehn}. We only have to note that 
${\pr_{X}}_*:{X}\times
\Sigma\rightarrow X$ is a smooth projective morphism of relative dimension
$1$ and $E$ being locally free is flat over $X$. 
The proof is rather simple so we sketch it here.

Let us denote by $E_x$ the vector bundle $E\mid_{\{x\}\times \Sigma}$
over $\Sigma$. Fix an ample line bundle $L$ on $\Sigma$. Then it is
well known that for big enough $k$ the vector bundle $E_x\otimes L^k$ 
is generated by its sections and $H^1(\Sigma;E_x\otimes L^k)=0$. Let us denote by $X_k\subset X$ those
points $x$ for which $E_x\otimes L^k$ is generated by its sections and
$H^1(\Sigma;E_x\otimes L^k)=0$. It
is standard that $X_k$ is a Zariski open subset of $X$. Thus we have a
covering $X=\bigcup X_k$ of $X$ by Zariski open subsets. It is again
standard that the Zariski topology of an algebraic variety is
noetherian\footnote{Cf. Example 3.2.1 on p. 84 of \cite{hartshorne}.},
which yields that we have some $k$ such that $X_k=X$. It is
now immediate that $$F=\pr_{\Sigma}^*(L^{-k})\otimes
\pr_X^*\left( (\pr_X)_*(E\otimes \pr_\Sigma^*(L^k))\right)$$ has the required
properties.

The result follows. 
\end{proof} 

\begin{proposition} Let ${\Sigma}$ be a smooth projective curve 
and $X$ be a smooth quasi-projective variety. Let 
$\cE=E\stackrel{f}{\rightarrow} F$ be a complex of vector bundles
on $X\times {\Sigma}$. Let $g_F:A\rightarrow F$ be 
a sectionless resolution of $F$. Let $M$ be the fibred product of 
$f$ and $g_F$. This comes with projection maps $p_F:M\rightarrow F$ and
$p_{A}:M\rightarrow A$. Let $g_M:A_2\rightarrow M$ 
be a sectionless resolution of $M$, and denote $j=g_M \circ p_{A_2}$. Finally, let
$A_1=\ker{g_M}$ and $i:A_1\rightarrow A_2$ the embedding. 
The situation is shown 
in the following diagram:
$$
\begin{array}{cc}
&E\hskip.2cm\stackrel{f}{\longrightarrow}\hskip.2cm F\\
&\hskip.3cm\nwarrow\hskip1cm\\
&\ \hskip.5cm M \hskip.4cm\uparrow\\
&\nearrow\hskip.4cm\searrow\\
0\longrightarrow A_1\stackrel{i}{\longrightarrow}& A_2\hskip.2cm
\stackrel{j}{\longrightarrow} \hskip.2cm A
\end{array}
$$
In this case 
the cohomology of the complex 
$$R^1{\pr_X}_*(A_1)\stackrel{i_*}{\longrightarrow} 
R^1{\pr_X}_*(A_2) \stackrel{j_*}{\longrightarrow} 
R^1{\pr_X}_*(A)$$
calculates the sheaves $\R^0{\pr_X}_*(\cE)$,
 $\R^1{\pr_X}_*(\cE)$ and $\R^2{\pr_X}_*(\cE)$
respectively. 
In other
words 
\begin{eqnarray}\R^0{\pr_X}_*(\cE)&\cong& \ker(i_*)\label{r0}\\  
\R^1{\pr_X}_*(\cE)&\cong &\ker(j_*)/\im(i_*) \label{r1}\\
\R^2{\pr_X}_*(\cE)&\cong& \coker(j_*)\label{r2}.
\end{eqnarray}
\label{diagram}
\end{proposition}
\begin{proof} Let us recall the definition of the fibred product: 
$$M=\ker(f-g_F:E\oplus A\rightarrow F).$$ This comes equipped with
 two obvious projections $p_E:M\rightarrow E$ and
$p_{A}:M\rightarrow A$. Because $g_F$ is surjective,
$f-g_F$ is also surjective. Thus $M$ is a vector bundle. By construction
the kernel
of $p_E$ is isomorphic to the kernel of $g_F$. Denote it by $B$. 
This says that the
following diagram is commutative and has two exact columns: 
$$
\begin{array}{c}0\longrightarrow 0\\
		\uparrow \hskip1cm \uparrow \\
		 E\stackrel{f}{\longrightarrow} F \\
		\mbox{\scriptsize{$p_E$}} \uparrow  \hskip1cm \uparrow
		\mbox{\scriptsize{$g_F$}}\\
		M\stackrel{p_{A}}{\longrightarrow} A\\
		\uparrow\hskip1cm\uparrow\\
		B\stackrel{\cong}{\longrightarrow} B\\
		\uparrow\hskip1cm\uparrow\\
		0\longrightarrow 0
\end{array}		
$$

If ${\cal A}$ denotes
the complex ${\cal A}=M\stackrel{p_{A}}{\rightarrow}A$ and ${\cal
B}$ the complex ${\cal B}=B\stackrel{\cong}{\rightarrow} B$, then the
above diagram is just a short exact sequence of complexes 
$$0\longrightarrow {\cal B}\longrightarrow {\cal A}\longrightarrow
{\cE}\longrightarrow 0.$$

Clearly $\R^i{\pr_X}_*({\cal B})$ vanishes for all $i$. 
(Any hypercohomology of an isomorphism is $0$.) Thus the long exact
sequence of the above short exact sequence gives the isomorphisms 
\begin{eqnarray}
\R^0{\pr_X}_*({\cal E})&\cong & \R^0{\pr_X}_*({\cal A})
\label{1r0}\\
\R^1{\pr_X}_*({\cal E})&\cong & \R^1{\pr_X}_*({\cal A})
\label{1r1}\\
\R^2{\pr_X}_*({\cal E})&\cong & \R^2{\pr_X}_*({\cal A})
\label{1r2}
\end{eqnarray}
 
Because $A$ is a sectionless resolution of $M$, we have 
$R^0{\pr_X}_*(A)=0$
thus the long exact sequence of the push forward of the complex $\cal
A$ breaks up into two exact sequences:
$$0\rightarrow \R^0{\pr_X}_*({\cal A})\rightarrow 
R^0{\pr_X}_*(M)\rightarrow 0,$$
and 
$$
0 \longrightarrow 
\R^1{\pr_X}_*({\cal A})\longrightarrow R^1{\pr_{\cal
W}}_*(M)\stackrel{{p_{A}}_*}\longrightarrow R^1{\pr_X}_*(A)
\longrightarrow \R^2{\pr_X}_*({\cal A})\longrightarrow 0.$$
Thus 
\begin{eqnarray}
\R^0{\pr_X}_*({\cal A})&\cong& R^0{\pr_X}_*(M)
\label{2r0}\\ 
\R^1{\pr_X}_*({\cal A})&\cong& \ker({p_{A}}_*)
\label{2r1}\\ 
\R^2{\pr_X}_*({\cal A})&\cong& \coker({p_{A}}_*). \label{2r2}
\end{eqnarray}

Now consider the short exact sequence:
\begin{eqnarray}0\longrightarrow A_1\stackrel{i}{\longrightarrow} A_2 
\stackrel{g_M}{\longrightarrow} M\longrightarrow 0. \label{exact1}\end{eqnarray}
$R^0{\pr_X}_*(A_2)=0$ because $A_2$ is a 
sectionless resolution of $M$ and hence we get the exact
sequence of sheaves:
\begin{eqnarray}0\longrightarrow R^0{\pr_X}_*(M)\longrightarrow 
R^1 {\pr_X}_*(A_1)\stackrel{i_*}{\longrightarrow} 
R^1{\pr_X}_*(A_2)
\stackrel{{g_M}_*}{\longrightarrow} R^1{\pr_{\cal
W}}_*(M)\longrightarrow 0.
\label{exact}
\end{eqnarray}
Thus $\ker(i_*) \cong R^0{\pr_X}_*(M)$ which by (\ref{2r0}) and
(\ref{1r0}) proves (\ref{r0}). 
Since ${g_M}_*$ is a surjection $\coker(j_*)\cong
\coker({p_{A}}_*)$. This together with (\ref{2r2}) and (\ref{1r2})
give (\ref{r2}).

Finally, consider the commutative diagram:
$$
\begin{array}{c}R^1{\pr_{\cal
W}}_*(M)\stackrel{\cong}{\longrightarrow} 
R^1{\pr_X}_*(M)\\
\mbox{\scriptsize{${g_M}_*$}}\uparrow \hskip1.5cm \downarrow
\mbox{\scriptsize{${p_{A}}_*$}}\\
R^1{\pr_X}_*(A_2)\stackrel{j_*}{\longrightarrow}R^1{\pr_X}_*(A)
\end{array}
$$
Since ${g_M}_*$ surjective by (\ref{exact}) we get that $\ker(j_*)/
\ker({g_M}_*)\cong \ker({p_{A}}_*)$. From (\ref{exact}) clearly 
$\ker({g_M}_*)\cong \im(i_*)$, thus $\ker(j_*)/\im(i_*)\cong 
\ker({p_{A}}_*)$. This together with (\ref{2r1}) and (\ref{1r1})
prove (\ref{r1}).
\end{proof}

\begin{corollary} If $\R^2{\pr_X}_*(\cE)=0$, 
in the situation of Proposition~\ref{diagram},
 then there exist two vector bundles $V$ and
$W$ over $X$ together with a homomorphism $f:V\rightarrow W$,
whose kernel and cokernel are $\R^0{\pr_X}_*(\cE)$ and 
$\R^1{\pr_X}_*(\cE)$ respectively. I.e. the following sequence
is exact:
$$0\rightarrow \R^0{\pr_X}_*(\cE)
\rightarrow V\stackrel{f}{\rightarrow} W \rightarrow
\R^1{\pr_X}_*(\cE)\rightarrow
0.$$
\label{cor}
\end{corollary}

\begin{proof} From the long exact sequence corresponding to
(\ref{exact1}), we have $R^0{\pr_X}_*(A_1)=0$. Let $V$ be the
vector bundle $R^1{\pr_X}_*(A_1)$.

Moreover $R^1{\pr_X}_*(A_2)$ and $R^1{\pr_X}_*(A)$ are
also vector bundles because $A_2$ and $A$ are sectionless
resolutions. Furthermore the assumption $\R^2{\pr_X}_*(\cE)=0$
shows that $j_*$ is surjective. Let $W$ be the vector bundle
$\ker(j_*)$, and $f$ be the map $i_*:V\rightarrow W$. 

The result follows from  Proposition~\ref{diagram}.
\end{proof}

The proof of Theorem~\ref{deg} is completed by 
Corollary~\ref{cor} noting that by
Corollary~\ref{vanish} we have $\R^2{\pr_{\tM}}_*(\cE^k_{\tM})=0$.
$\square$

\newpage
\section{The downward degeneracy locus $DD_k$}
\label{downdegeneracylocus}

\begin{definition} The {\em downward degeneracy locus} $$DD_k:=\{\cE\in \tM: 
\Hy^0(\Sigma,\cE_\tM^k)\neq 0)\}$$ 
is the locus where $\D_k$ fails to be a
vector bundle, i.e. where $f$ of Theorem~\ref{deg} fails to be an 
injection.  
\end{definition}

The aim of this section is to give a description of the degeneracy
locus $DD_k$. For this we need a refinement  of Theorem~\ref{th-nil},
which still follows from the proof of Proposition (19) of \cite{thaddeus1}.

\begin{proposition} The nilpotent cone is a compact union of $3g-3$ 
dimensional manifolds: $$N=\N\cup\bigcup^{g-1}_{k=1}D_k,$$
where each $D_k$ is biholomorphic to the total space of the negative vector
bundle $E^-_k$ over $F_k$, the $k$-th component of the fixed point set of
the $\C^*$-action. 
The component $D_k$ can also be characterised as the locus of those stable Higgs
bundles $\cE=\compE$ which have a unique 
subbundle $L_\cE$ of degree $(1-k)$ killed by the non-zero Higgs
field $\Phi$. 
\label{cone}
\end{proposition} 

\begin{proof}The first part is proved in Theorem~\ref{th-nil}.

For the second part consider a universal Higgs bundle $\cE_\M$
over $\M\times \Sigma$ restricted to $D_k\times \Sigma$. Let us denote
it by $\cE_k=\E_k\stackrel{\uPhi_k}{\rightarrow}\E_k\times
\K$. Consider the kernel of $\uPhi_k$. Because $D_k$ parametrizes
nilpotent stable Higgs bundles with non-zero Higgs field
$\ker(\uPhi_k)$ is a line bundle over $D_k\times \Sigma$. Recall from 
Proposition 7.1 of \cite{hitchin1} that for $\compE\in F_k\subset D_k$ 
we have $\deg(\ker(\Phi))=1-k$. Since $D_k$ is smooth we have that 
$\deg(\ker(\Phi))=1-k$ for every $\compE\in D_k$.

The result follows.    
\end{proof}

\begin{remark} A completely analogous result holds for
$\widetilde{N}$ with $\tN$, $\widetilde{D}_k$ and $\widetilde{F}_k$ instead of
$\N$, $D_k$ and $F_k$.
\end{remark}

\begin{theorem} Let $k=1,..,g-1$. 
The degeneracy locus $DD_k$ has the following
decomposition: $$DD_k=\tN^k\cup\bigcup_{i=1}^k \widetilde{D}^k_i,$$
where $\tN^k=DD_k\cap \tN$, and $\widetilde{D}_i^k\subset \widetilde{D}_i$ are
those nilpotent stable Higgs bundles whose unique line bundle $L_\cE$ of 
Proposition~\ref{cone} has the property that $H^0(\Sigma,L_\cE\otimes
L_p^{k-1})\neq 0$.

Furthermore $$\widetilde{D}^k_k:=\{\cE\in\widetilde{D}_k: L_\cE=L^{1-k}_p\}$$
and hence\footnote{Recall the definition of $\eta^Y_X$ from Notation~\ref{cohclass}.} 
\begin{eqnarray}
\eta^\tM_{\widetilde{D}^k_k}\smallsetminus [\J]= \eta^\M_{D_k}
\in H^{6g-6}(\M), 
\label{splitE}
\end{eqnarray}
where $\eta^\tM_{\widetilde{D}^k_k}\smallsetminus [\J]$ means the
coefficient of $\eta^\J_{pt}$ in the decomposition of (\ref{splitM}).
\label{degeneracy}
\end{theorem} 

\begin{proof} Let $\cE=\compE$ be a stable Higgs bundle with 
$\Phi\neq 0$ and 
$\Hy^0(\Sigma,\cE\otimes L_p^{k-1})\neq 0$. It is easy to see that
this hypercohomology is the vector space of all morphisms from 
$\cE_0\otimes
L_p^{1-k}=L_p^{1-k}\stackrel{0}{\rightarrow}L_p^{1-k}\otimes K$ to $\cE$.
Consider a nonzero such morphism $f$. Consider $L$ the line subbundle
of $E$ generated by the image of $f$ of
Lemma~\ref{narasimhan-seshadri}. 
Clearly $L$ is killed by the Higgs field $\Phi$. This shows that
$\cE\in \widetilde{N}$ and $L=L_\cE$. We also see that 
$\Hy^0(\Sigma,\cE\otimes L_p^{k-1})\cong H^0(L_\cE\otimes L_p^{k-1}$). The
first part of the statement follows. 

By the above argument it follows that 
$\widetilde{D}^k_k=\{\cE\in \widetilde{D}_k: 
H^0(\Sigma, L_\cE \otimes L_p^{k-1})\neq 0\}$, however $L_\cE$ is of
degree $1-k$, thus  
$\widetilde{D}^k_k=\{\cE\in\widetilde{D}_k: L_\cE=L^{1-k}_p\}$, as claimed. 
This means that for every $\cE\in D_k$ there is a unique line bundle 
$L=L^{1-k}_p\otimes L^*_\cE$
such that $\cE\otimes L\in \widetilde{D}^k_k$. This shows (\ref{splitE}).
\end{proof}

\begin{remark} By definition 
$\widetilde{\N}^k=W^0_{2,2k-1}$ are non-Abelian Brill-Noether
loci as defined in  \cite{sundaram}\footnote{Cf. \cite{teixidor}.}.
\end{remark}

\newpage
\section{Proof of Theorem~\ref{main}}
\label{proof}

In this final section we prove Theorem~\ref{main}. 

\paragraph{\it Proof of Theorem~\ref{main}.} The proof proceeds by
first showing that $\ch_0(\D_k)=4g-4$, then $c_{4g-3}(\D_k)=0$ and
concludes 
using Porteous' Theorem~\ref{porteous}  for $\D_k$. 

First we make some calculations.

\begin{lemma} The virtual bundle  $\D_k$ has rank
$4g-4$, i.e. $ch_0(\D_k)=4g-4$. Moreover
\begin{eqnarray}
c(\D_k)=\left(1+\alpha_\M+\frac{\alpha^2_\M-\beta_\M}{4}\right)^{2g-2}
\label{chern}
\end{eqnarray} 
in the decomposition
(\ref{splitM}). 
\label{calc}
\end{lemma}

\begin{proof} It follows from the hypercohomology long exact sequence 
that 
$$\D_k=-{\pr_\tM}_!(\cE^k_\tM)={\pr_\tM}_!(\E^k_\tM\otimes \K)-
{\pr_\tM}_!(\E^k_\tM).$$ We can calculate the Chern character of the
right hand side by the Grothendieck-Riemann-Roch theorem. This gives
$$\ch(\D_k)={\pr_\tM}_*\left(\ch(\E^k_\tM)(\ch
(\K)-1)\td(\Sigma)\right).$$
Now $\td(\Sigma)=1-(g-1)\sigma$ and $\ch(\K)=1+(2g-2)\sigma$.
Moreover ${\pr_\tM}_*$ maps a cohomology class $a\in H^*(\tM)\otimes
H^*(\Sigma)$ of the form $$a=a_0\otimes 1 + \sum_{i=1}^{2g} a_1^i
\otimes e_i +
a_2\otimes \sigma$$ to the class $a_2\in H^*(\tM)$. We will use the
notation $a\smallsetminus \sigma=a_2$ and $a\smallsetminus 1=a_0$. From this it follows that 
$$\ch(\D_k)=\left(\ch(\E^k_\tM)((2g-2)\sigma)(1-(g-1)\sigma)\right)
\smallsetminus\sigma=(2g-2)(\ch(\E^k_\tM)\smallsetminus 1).$$  

Observe that $$\ch(\E^k_\tM)\smallsetminus 1 =\ch((\E^k_\tM)_p)\in
H^*(\tM),$$ where $(\E^k_\tM)_p=\E^k_\tM\mid_{\tM \times \{p\}}$.
It follows from (\ref{c1}) and (\ref{c2}) 
that $c_1((\E^k_\tM)_p)=\alpha_\M$ and 
$c_2((\E^k_\tM)_p)=(\alpha_\M^2-\beta_\M)/4$. Hence the formal Chern
roots of $(\E^k_\tM)_p$ are $(\alpha_\M+\sqrt{\beta_\M})/2$ and 
$(\alpha_\M-\sqrt{\beta_\M})/2$. Thus $$\ch((\E^k_\tM)_p)=
\exp\left(\frac{\alpha_\M+\sqrt{\beta_\M}}{2}\right)
+\exp\left(\frac{\alpha_\M-\sqrt{\beta_\M}}{2}\right)=
2e^{\alpha_\M/2} \cosh\left(\sqrt{\beta_\M}/2\right),$$
and hence 
$$\ch(\D_k)=(4g-4)e^{\alpha_\M/2}\cosh\left(\sqrt{\beta_\M}/2\right).$$
This shows that $\rank(\D_k)=\ch_0(\D_k)=4g-4$ and formal calculation
gives (\ref{chern}).
\end{proof}

(\ref{chern}) has the following immediate corollary:

\begin{corollary} $c_{4g-3}(\D_k)=0$.$\square$
\label{vanishing} 
\end{corollary}
 
To prove Theorem~\ref{main} we exhibit $g$ linearly independent elements 
$r_0,r_1,..,r_{g-1}\in H^{6g-6}_{cpt}(\M)$ for which $j_\M(r_i)=0$. 

To construct $r_k$ for $0<k<g$ consider the Zariski open subvarieties
$$\tD_{\geq k}=\tM\setminus(\tN\bigcup_{i=1}^{k-1}\widetilde{D}_i)$$ and 
$$D_{\geq k}=\M\setminus(\N\bigcup_{i=1}^{k-1} D_i)$$ of $\tM$ and $\M$ respectively. 
Restricting 
the sequence of Theorem~\ref{deg} to $\tD_{\geq k}$ yields:
\begin{eqnarray}0\longrightarrow
\R^0{\pr_\M}_*(\cE_{\tM}^k)\mid_{\tD_{\geq k}}
\longrightarrow V\mid_{\tD_{\geq k}}\stackrel{f\mid_{\tD_{\geq k}}}{\longrightarrow} 
W\mid_{\tD_{\geq k}} \longrightarrow
\R^1{\pr_\M}_*(\cE_{\tM}^k)\mid_{\tD_{\geq k}}\longrightarrow
0.\label{sequ}
\end{eqnarray}

The degeneracy locus of $f\mid_{\tD_{\geq k}}$ (where $f\mid_{\tD_{\geq k}}$ fails
to be an injection) is $DD_k\cap \tD_{\geq k}$ which is $\widetilde{D}^k_k$ from
Theorem~\ref{degeneracy}. This has codimension $4g-3$. Furthermore 
$$\rank(W)-\rank(V)=\rank\left(\R^1{\pr_\M}_*(\cE_{\tM}^k)\right)-
\rank\left(\R^0{\pr_\M}_*(\cE_{\tM}^k)\right)=\rank(\D_k)=4g-4$$
by Lemma~\ref{calc}. Thus the degeneracy locus has the expected
codimension hence we are in the situation of Porteous'
Theorem~\ref{porteous}, which gives: 
$$\eta_{\widetilde{D}^k_k}^{\tD_{\geq k}}=c_{4g-3}(W\mid_{\tD_{\geq k}}-V\mid_{\tD_{\geq k}})
\in H^{8g-6}(\tD_{\geq k}).$$
The right hand side equals $c_{4g-3}(\D_k\mid_{\tD_{\geq k}})$ by
(\ref{sequ}), which vanishes by Corollary~\ref{vanishing}. Also
$$\eta^{\tD_{\geq k}}_{\widetilde{D}^k_k}\smallsetminus [\J]= \eta^{D_{\geq k}}_{D_k}$$
by (\ref{splitE}). It follows that 
\begin{eqnarray}\eta^{D_{\geq k}}_{D_k}=0\in H^{6g-6}(D_{\geq k}).\label{0}\end{eqnarray}

From now on we work over $\M$. We show by induction on $i$ that there is
a formal linear combination $$r^i_k=\sum_{j=k-i}^k 
\lambda_j\cdot\left[\eta^{D_{\geq k-i}}_{D_j}\right]$$ of cohomology classes in 
$H^{6g-6}(D_{\geq k-i})$, such that $\lambda_k=1$ and 
the corresponding cohomology class 
$\sum_{j=k-i}^k\lambda_i\cdot\eta^{D_{k-i}}_{D_j}$ is $0$ in
$H^{6g-6}(D_{\geq k-i})$. 

For $i=0$ the statement is just (\ref{0}). Suppose that there is 
such formal linear combination $r^i_k$. Consider the following bit of the
long exact sequence of the pair $D_{\geq k-i}\subset D_{\geq k-i-1}$:
$$ H^{6g-6}(D_{\geq k-i},D_{\geq k-i-1})\longrightarrow
H^{6g-6}(D_{\geq k-i-1})
\longrightarrow H^{6g-6}(D_{\geq k-i}).$$
Because $D_{\geq k-i-1}\setminus D_{\geq k-i}=D_{k-i-1}$ is of real codimension
$6g-6$, the Thom isomorphism
transforms this sequence to:
\begin{eqnarray}
H^{0}(D_{k-i-1})\stackrel{\tau}{\longrightarrow} H^{6g-6}(D_{\geq k-i-1})
\stackrel{\rho}{\longrightarrow} H^{6g-6}(D_{\geq k-i}),\label{thom}
\end{eqnarray}
where $\tau$ is the Thom map and $\rho$ is restriction. Clearly 
$\rho\left(\eta^{D_{\geq k-i-1}}_{D_j}\right)=\eta^{D_{\geq k-i}}_{D_j}$.
Thus $$\rho\left(\sum_{j=k-i}^k\lambda_j\cdot\eta^{D_{\geq k-i-1}}_{D_j}\right)
=\sum_{j=k-i}^k\lambda_j\cdot\eta^{D_{\geq k-i}}_{D_j}=0.$$  
The exactness of (\ref{thom}) yields that the cohomology class 
$$\sum_{j=k-i}^k \lambda_j\cdot\eta^{D_{\geq k-i-1}}_{D_j}$$ is in the
image of $\tau$. Because $H^0(D_k)\cong \Qu$ there is a rational number 
$-\lambda_{k-i-1}\in \Qu$ such that 
\begin{eqnarray}\tau(-\lambda_{k-i-1})=
\sum_{j=k-i}^k \lambda_j\cdot\eta^{D_{\geq k-i-1}}_{D_j}
\in H^{6g-6}(D_{\geq k-i-1})\label{akar}.
\end{eqnarray} 
However a well known property of the Thom map gives
$\tau(1)=\eta^{D_{\geq k-i-1}}_{D_{k-i-1}}$, thus from (\ref{akar}) the 
formal linear combination 
$$r^{i+1}_k=\sum_{j=k-i-1}^k 
\lambda_j\cdot\left[\eta^{D_{\geq k-i-1}}_{D_j}\right]$$ is
$0$, when considered as a class in $H^{6g-6}(D_{\geq k-i-1})$. This proves the 
existence of formal linear combinations $r^i_k$ for all $0\leq i\leq k-1$.

Using $r^{k-1}_k$ an identical argument gives the formal linear
combination $$r^\prime_k=\lambda\cdot\left[\eta^\M_\N\right]+
\sum_{j=1}^k \lambda_j\cdot\left[\eta^{\M}_{D_j}\right]$$ with the
property that $\lambda_k=1$ and $r^\prime_k$ when considered as an element of
$H^{6g-6}(\M)$ is $0$. Now the compactly supported cohomology class 
$$r_k=\lambda \cdot\overline{\eta}^\M_\N+
\sum_{j=1}^k \lambda_j\cdot \overline{\eta}^{\M}_{D_j}\in H^{6g-6}_{cpt}(\M)$$
has the property that $j_\M(r_k)=r^\prime_k=0$, where by abuse of
notation
$r^\prime_k$ denotes the cohomology class in $H^{6g-6}(\M)$
corresponding to the formal linear combination $r^\prime_k$.

We have found $g-1$ linearly independent compactly supported 
cohomology classes $r_1,..,r_{g-1}\in H^{6g-6}_{cpt}(\M)$. Clearly 
$\overline{\eta}^{\M}_\N$ is not in the span of $r_1,..,r_{g-1}$. Moreover
for each $0<i<g$ we have $\int_\M \overline{\eta}^{\M}_\N \wedge r_i=0$ since
$j_\M(r_i)=0$. Furthermore   
$$\int_\M \overline{\eta}^{\M}_\N\wedge \overline{\eta}^{\M}_\N=\int_\N
c_{3g-3}(T^*_\N)= 0.$$ Thus $\overline{\eta}_\N^\M$ 
is perpendicular to  $r_1,..,r_{g-1}$
and $\overline{\eta}_\N^\M$, which constitutes a basis for $H^{6g-6}_{cpt}(\M)$, and
so  $j_\M(\overline{\eta}^{\M}_\N)=0$.

Putting all this together: we have $g$ linearly independent
middle dimensional compactly supported classes $r_0=\overline{\eta}^{\M}_\N$
and $r_1,..,r_{g-1}$ 
in the kernel of the forgetful map 
$j_\M: H^{6g-6}_{cpt}(\M)\rightarrow H^{6g-6}(\M)$. 

Theorem~\ref{main} is finally proved.  {$\square$\vskip0.4cm}

\newpage
\thispagestyle{empty}

\chapter{Cohomology}
\label{cohomology}

As we already noted in the Introduction, to
understand the physical model of \cite{bershadsky-et-al} one needs to
have a good understanding of the cohomology ring of $\tM$.
The present chapter\footnote{This chapter describes a joint work with Michael
Thaddeus.} attempts to fill the gap in the literature by
providing at least a half-proved complete description of
$H^*(\M)^\Gamma$, and in turn of $H^*(\tM)$, 
which agrees with computer calculations of genus up to $7$. 

We start with constructing equivariant structures on the universal
bundles of Section~\ref{universalbundles}. Then we construct the
equivariant virtual Mumford bundle in Section~\ref{eqmumford}, 
and investigate its degeneracy
locus in Section~\ref{updegeneracylocus}.  Applying the equivariant 
Porteous' theorem  we are able to
deduce a fundamental proposition, which implies, through the general
arguments in Subsection~\ref{equivariantcoh}, that $H_\circ^*(\tM)$
and in turn $H^*(\tM)$ are generated by some universal classes. 

We finish by providing a conjectured complete description of the 
subring $H^*_I(\M)\subset H^*(\M)^\Gamma$, 
generated by the classes $\al_\M$, $\be_\M$ and $\ga_\M$. 
We support it by a few results. We
prove that $H^*_I(\M)$ has the same Poincar\'e polynomial as the
conjectured ring. We also find the first two relations and 
prove here that the second is Newstead's relation $\be^g=0$,
which easily yields that the Chern classes of $\M$ are zero in
degrees at least $2g$.

\newpage
\section{Equivariant universal bundles}

We would like to incorporate into the universal bundle the $\C^*$-action on $\tM$.
Recall that $\C^*$ acts on $\tM$ by scalar multiplication of the
Higgs field. Moreover let $\C^*$ act trivially on
$\Sigma$ so that  we get the diagonal $\C^*$-action on $\tM\times \Sigma$. 
We will need equivariant universal bundles with respect to this action:

\begin{proposition} Let 
$\cE_\tM=\E_\tM\stackrel{\uPhi}{\rightarrow}\E_\tM\otimes
\pr_\Sigma^*(K)$
be a universal Higgs bundle over $\tM\times \Sigma$.
Then there is a $\C^*$-equivariant structure on $\cE_\tM$, 
i.e. an
equivariant bundle 
structure on $\E_\tM$ and $\E_\tM\otimes K_\Sigma$ such that
the universal Higgs field 
${\uPhi}:\E^\circ_\tM\rightarrow \E^\circ_\tM\otimes K_\Sigma^\circ$ is
equivariant. 
\end{proposition}  

\begin{proof} Recall the construction of the universal Higgs bundle
over $\tM\times \Sigma$ from Proposition~\ref{universal}. In the gauge
theory picture $\C^*$ acts on $\calB$ by scalar multiplication of the
Higgs field\footnote{This $\C^*$ does not have anything
to do so far with the constant scalar gauge transformations 
$\C^*\subset \G^c$. However see
Section~\ref{compactthesis}, where we suggest that they are
closely related.}. If we let $\C^*$ act trivially on $\calC$, then we have
that $\pr:\calB\to \calC$ is $\C^*$-equivariant. Also let
$\E_\calC\otimes \bL^{-1}$ be a $\C^*$-equivariant holomorphic bundle
over $\calC\times \Sigma$ with the trivial $\C^*$-action. Thus the
pullback of $\E_\calC\otimes \bL^{-1}$ by the $\C^*$-equivariant map
$\pr$ gives a $\C^*$-equivariant structure on $\E_{\calB}$, we denote
the resulting $\C^*$-equivariant bundle by $\E^\circ_{\calB}$. This
descends in the quotient to an equivariant universal bundle
$\E^\circ_\tM$. 

We need a $\C^*$-equivariant structure on $\calO_\tM$ with homogeneity $1$, in other
words a $\C^*$-equivariant line bundle with equivariant first Chern class
$u$, where $u$ is the integer generator of the $\C^*$-equivariant cohomology
of a point. Such an equivariant line bundle exists by Theorem 1 of 
\cite{edidin-graham2}. In our case we can think of this line bundle as
the invertible subsheaf of the sheaf $\Omega^2(\tM)$ of holomorphic
two-forms 
on $\tM$ generated by the holomorphic symplectic form 
$\omega_h$ on $\tM$, which is homogeneous of degree $1$. 
In any case let us denote it by $\calO^\circ_\tM$. Moreover we denote
by $K^{\circ}_\Sigma$ the equivariant line bundle
\begin{eqnarray} \pr_\tM^*(\calO^\circ_\tM)\otimes
\pr_\Sigma^*(K)\label{marcsakharomnap}\end{eqnarray} on $\tM\times
\Sigma$, 
where by abuse of
notation $K$ stands for the trivial equivariant structure on the
canonical bundle $K$ on $\Sigma$. 

Now it is easy to check that the universal Higgs field
$\uPhi_{\tM}:\E_\tM^{\circ}\rightarrow \E_\tM^{\circ}\otimes K_\Sigma^\circ$ is $\C^*$-equivariant. 
\label{equivariant}
\end{proof}

Having proved the existence of the equivariant universal bundle
$\E_\tM^\circ$, we can restrict it to $\E_\M$ to get a universal bundle
over $\M\times \Sigma$ and consider, analogously to
(\ref{universalclassesM}), 
the K\"unneth decomposition of $\End(\E^\circ_\M)$ to get equivariant
universal classes:
\begin{eqnarray}c_2(\End(\E^\circ_\M))=2\eqal\otimes \sisi+\sum_{i=1}^{2g}4\eqpsi_i 
\otimes \xisi_i -
\eqbe_\M\otimes 1\label{equniversalclassesM} \end{eqnarray}
in $$H^4_\circ(\M\times \Sigma)\cong
\sum_{r=0}^{4} H_\circ^r(\M)\otimes_{\Qu [u]}
H_\circ^{4-r}(\Sigma)\cong \sum_{r=0}^{4} H_\circ^r(\M)\otimes
H^{4-r}(\Sigma)$$
for some {\em equivariant universal classes} 
$\eqal\in H_\circ^2(\M)$, $\eqpsi_i\in H_\circ^3(\M)$ and $\eqbe\in
H_\circ^4(\M)$. 

We will need to know the restriction of the
equivariant universal classes 
to  the fixed point set of the $\C^*$-action
on $\M$. First consider $\N=F_0$. The equivariant universal bundle
$\E_\M^\circ$ as constructed in the  proof of Proposition~\ref{equivariant} 
clearly restricts to
$\N\times \Sigma$ to $\E_\N$ with the trivial $\C^*$-action on it. 
Consequently $$\eqal\mid_\N=\al_\N\in H^2(\N)\subset
H^2_\circ(\N),$$ $$\eqpsi_i\mid_\N=\psi^i_\N\in H^3(\N)\subset
H^3_\circ(\N)$$ and $$\eqbe\mid_\N=\be_\N\in H^4(\N)\subset H^4_\circ(\N).$$ 

Consider now the restriction of the universal classes to $F_d$ for $d>0$.
Because the $\C^*$-action is trivial on $F_d$ we have 
$H^*_{\circ}(F_d)\cong H^*(F_d)\otimes H^*_{\circ}(pt)\cong
H^*(F_d)\otimes \Qu [u]$. 
Note that since all the universal classes are 
invariant under $\Gamma$, their restrictions to $F_d$ 
live\footnote{Recall that $\bd=2g-2d-1$.} in $(H^*(F_d))^\Gamma\cong H^*(\Sigma_\bd)$. 
Recall the ring $H^*(\Sigma_\bd)$ from 
Subsection~\ref{symmetric}. 

\begin{lemma} The equivariant universal classes restrict to $F_d$ as follows:
\begin{eqnarray} \eqal\mid_{F_d} &=& (2d-1)(\eta-u)+\sigma;\nonumber\\
   {\eqpsi_i\mid_{F_d}}&=& \left\{ 
\begin{array}{ccc}\ \ \frac{(\eta-u)}{2}\xi_{i+g} &\mbox{if}&\ i\leq g; \\
                  -\frac{(\eta-u)}{2}\xi_{i-g} &\mbox{if} &\ i>g; 
   	         \end{array} \right. \nonumber\\
   \eqbe\mid_{F_d}&= &(\eta-u)^2. \nonumber
\end{eqnarray}
\label{restrictions}
\end{lemma}
\begin{proof} 
First we construct an equivariant universal Higgs bundle over $F_d\times
\Sigma$. Let $\bP_d$ a normalized Poincar\'e bundle over
$\J_d\times \Sigma$
as in Subsection~\ref{cohomologyofJ}. By abuse of notation we also
denote by $\bP_d$ the pullback\footnote{Recall the projection
$\pr_{\J_d}$ from (\ref{fibredproduct}).} 
of $\pr_{\J_d}^*(\bP_d)$ to $F_d\times
\Sigma$. We also denote by $\Delta_{\bd}$ the pullback
$\pr_{\Sigma_{\bd}}^*(\Delta_{\bd})$
of the universal divisor\footnote{Cf. Subsection~\ref{symmetric}.} 
on $\Sigma_{\bd}\times \Sigma$. Then $\bP_d^2 K_\Sigma^{-1}
\Lambda_{\Sigma}^{-1} (\Delta_{\bd})$ defines a line bundle over $F_d\times
\Sigma$ which is trivial over $x\times \Sigma$ for all $x\in F_d$. By
the push-pull formula this is the pullback of some line bundle
$L_{F_d}$ 
over $F_d$. So \begin{eqnarray*}\calO(\Delta_{\bd})\cong L_{F_d} K_{\Sigma}
\Lambda_{\Sigma} \bP^{-2}_d;\label{delta}\end{eqnarray*} hence the latter has a section, which we
denote by $\phi_d$, vanishing on
$\Delta_{\bd}$. Now we let $$\E_{F_d}=\bP_d\oplus \bP_d^{-1}\Lambda_{\Sigma}
L_{F_d}$$ and $$\uPhi_d=\left(\begin{array}{cc} 0&0\\ \phi_d &
0\end{array}\right),$$ so $\uPhi_d\in H^0(F_d\times \Sigma;
\Hom(\E_{F_d},\E_{F_d}\otimes K_\Sigma))$. By construction
$\E_{F_d}\stackrel{\uPhi_d}{\to} \E_{F_d}\otimes K_\Sigma$ is a
universal Higgs bundle over $F_d\times \Sigma$. 

Then we let
$\calO^\circ$ be the weight $1$ equivariant
structure on the trivial bundle
$\calO_{F_d\times\Sigma}$. Moreover we let 
$\bP_d^\circ=\calO^\circ\bP_d$, $\K^\circ=\calO^\circ \K$ and
$\Lambda_\Sigma^\circ=\calO^\circ\Lambda_\Sigma$. Then it follows that as
equivariant bundles: 
 \begin{eqnarray}\calO(\Delta_{\bd})\cong L_{F_d} K^\circ_{\Sigma}
\Lambda^\circ_{\Sigma} (\bP^\circ_d)^{-2}.\label{eqdelta}\end{eqnarray}

Now if we let $\E^\circ_{F_d}= \bP^\circ_d\oplus 
(\bP^\circ_d)^{-1}\Lambda^\circ_\Sigma
L_{F_d}$, then 
$\E^\circ_{F_d}\stackrel{\uPhi_d}{\to} \E^\circ_{F_d}\otimes
K^\circ_\Sigma$ is an equivariant universal Higgs bundle over
$F_d\times \Sigma$. It
follows from Lemma~\ref{ramanan} that the equivariant universal
classes restricted to $F_d$ will appear as the K\"unneth components of
$c^\circ_2(\End(\E^\circ_{F_d}))$. We can calculate them as follows:
$$c^\circ_2(\End(\E^\circ_{F_d}))=-\left(c^\circ_1((\bP^\circ_d)^2(\Lambda^\circ_\Sigma)^{-1}
L_{F_d}^{-1})\right)^2= 
-\left(c^\circ_1\left(\calO(-\Delta_\bd)(K^\circ_\Sigma)\right)\right)^2$$
from (\ref{eqdelta}). We can calculate this since from
(\ref{cherndelta}) we have  $$c_1(\Delta_\bd)=\eta\otimes 1 +
 \sum_{i=1}^{g}(\xi_i
\otimes \xisi_{i+g}-\xi_{i+g}\otimes \xisi_i)  + 
\bd\otimes \sisi\in H^2(F_d\times \Sigma)\cong
\sum_{r=0}^{2} H^r(F_d)\otimes H^{2-r}(\Sigma)$$ and from (\ref{marcsakharomnap})
$$c^\circ_1(K_\Sigma^\circ)=u\otimes 1 + (2g-2)\otimes \sisi \in H^2(F_d\times \Sigma)\cong
\sum_{r=0}^{2} H^r(F_d)\otimes H^{2-r}(\Sigma).$$ Thus we have 

\begin{eqnarray*} c^\circ_2(\End(\E^\circ_{F_d})) & = & 
-\left(c^\circ_1\left(\calO(-\Delta_\bd)(K^\circ_\Sigma)\right)\right)^2
\\ & = &
-\left((\bd-(2g-2))\otimes \sisi +
 \sum_{i=1}^{g}(\xi_i
\otimes \xisi_{i+g}-\xi_{i+g}\otimes \xisi_i) +  (\eta-u)\otimes 1 \right)^2 \\
& = &
- \left( \sum_{i=1}^{g}(\xi_i \otimes \xisi_{i+g}-\xi_{i+g}\otimes \xisi_i)
\right)^2 - 
2(\eta-u) \sum_{i=1}^{g}(\xi_i
\otimes \xisi_{i+g}-\xi_{i+g}\otimes \xisi_i)\\ & &  - (\eta-u)^2 -
2(\eta-u)(1-2d)\otimes \sisi 
\\  &=& 2\left((2d-1)(\eta-u)+\sisi\right)\otimes \sisi_\Sigma -2 (\eta-u)\sum_{i=1}^{g}(\xi_i
\otimes \xisi_{i+g}-\xi_{i+g}\otimes \xisi_i) \\ & & - (\eta-u)^2\otimes 1.
\end{eqnarray*}
Comparing this with (\ref{equniversalclassesM}) proves the result. 
\end{proof}

\newpage
\section{The equivariant virtual Mumford bundle $\bM_d^\circ$}
\label{eqmumford}

In Section~\ref{generationsection} we will prove the generation theorem by considering the
degeneration locus of the equivariant virtual Mumford bundle. 
 To do this, first we define the equivariant virtual Mumford bundle as:
$$\bM^\circ_d=-{\pr_\tM}_!(\E^{\circ -d}_\tM)=
-R^0{\pr_\tM}_*(\E^{\circ -d}_\tM)+R^1{\pr_\tM}_*
(\E^{\circ -d}_\tM )\in K(\tN),$$
where $\E^{\circ -d}_\tM$ denotes $\E^\circ_\tM\otimes \pr_\Sigma^*(L_p^{-d})$.

The next theorem says that the equivariant virtual Mumford bundle can
be thought of as the equivariant degeneracy sheaf of a homomorphism
of equivariant vector bundles.  

\begin{theorem} There exist two equivariant vector bundles $V^\circ$
and $W^\circ$, together with an equivariant homomorphism
$f^\circ:V^\circ\to W^\circ$   such that the following sequence of coherent sheaves is
exact: $$0\to R^0{\pr_\tM}_*(\E^{\circ-d}_\tM) \to V^\circ \stackrel{f^\circ}{\to} W^\circ
\to 
R^1{\pr_\tM}_*(\E^{\circ-d}_\tM)\to 0.$$
\label{degeneracym}
\end{theorem}
\begin{proof} Choose an effective divisor $D$ on $\Sigma$ 
such that $H^1(\Sigma,E\otimes \calO(D))=0$ for all
stable Higgs bundle $\compE\in \tM$. Such $D$ exists because the
Harder-Narasimhan type of vector bundles occurring in stable Higgs pairs is
bounded\footnote{Cf. Corollary 3.3 of \cite{nitsure}.}. Then tensoring
$$0\to \calO_\Sigma \to \calO(D)\to \calO_D\to 0$$ with
$\E^{\circ -d}_\tM$ and
pushing down by $\pr_\tM$ yields \begin{eqnarray*}
0\to R^0{\pr_\tM}_*(\E^{\circ -d}_\tM)\to R^0{\pr_\tM}_*(\E^{\circ
-d}_\tM\otimes \calO(D))\to \\ \to 
R^0{\pr_\tM}_*(\E^{\circ -d}_\tM\otimes \calO_D)\to
R^1{\pr_\tM}_*(\E^{\circ -d}_\tM)\to 0,\end{eqnarray*}
because of the condition on $D$. However our hypothesis on $D$ also implies
that the second and third terms of the exact sequence are vector
bundles, proving the desired result. 
\end{proof}

\newpage
\section{The upward degeneration locus $UD_d$}
\label{updegeneracylocus}

Here we consider the degeneracy locus $UD_d$ of $f^\circ$ of the above
Theorem~\ref{degeneracym}.
It can be defined as $$UD_d=\left\{(E,\Phi)\in \tM:
H^0(\Sigma;E\otimes L_p^{-d})\not=0\right\},$$
and called the {\em upward degeneracy locus}, in contrast with the
downward degeneracy locus of Subsection~\ref{downdegeneracylocus}.  
If $E$ is stable then we see that $H^0(\Sigma;E\otimes L_p^{-d})=0$,
thus $T^*_\tN\cap UD_d=\emptyset$. Furthermore if the destabilizing
bundle of $E$ is of degree less than $d$, then $(E,\Phi)\not\in
UD_d$. More specifically we have the following 
description of the degeneration locus
$UD_d$:

\begin{theorem} For $1\leq d \leq g-1$
the degeneration locus $UD_d$ has the following decomposition: 
$$UD_d=\bigcup_{k=d}^{g-1} UD_d^k,$$ where $UD^k_d\subset
\widetilde{U}_k$ are those stable Higgs bundles $(E,\phi)$ for which $L$, 
the destabilizing line bundle of $E$, has the property that
$H^0(\Sigma,L\otimes L_p^{-d})\not= 0$.

Finally the real codimension of $UD_d^d$ in $\tM$ is $2(2g+2d-2)$. 
\label{descriptionudd}
\end{theorem} 

\begin{proof} By definition $\compE\in UD_d$ if and only if
$H^0(\Sigma,E\otimes L_p^{-d})\not= 0$. A non-zero section however generates
a line bundle of $E\otimes L_p^{-d}$ of non-negative degree, which in
turn gives a line subbundle $L$ of $E$ of degree $d$ (which, by its
uniqueness, should be the destabilizing line bundle of E), such that
$H^0(\Sigma;L\otimes L_p^{-d})\not= 0$. The first statement follows.

Now consider the map $f_k:U_k\to \J_{k-d}$ sending $\compE$ to $L\otimes
L_p^{-d}$ and consider also the Abel-Jacobi map 
$u_k:\Sigma_{k-d}\to \J_{k-d}$. We have $(E,\Phi)\in
U^k_d$ if and only if $f_k(E,\Phi)\in u_k(\Sigma_{k-d})$. The 
theorem now easily follows.  
\end{proof}

\newpage
\section{Generation theorem}
\label{generationsection}

The aim of this section is to prove that $H^*_\circ(\tM)$ is generated
by the equivariant universal classes. To prove this we apply 
Corollary~\ref{eqporteous} to $UD_d^d$. First we need some notation. 

Consider the equivariant
Chern class $$c^\circ_{2g+2d-2}(\bM^\circ_d)\in
H_\circ^{2(2g+2d-2)}(\tM).$$
Write it down in the K\"unneth decomposition 
$$H^*_\circ(\tM)\cong H^*_\circ(\M)\otimes_{\Qu[u]} 
H_\circ^*(T^*_\J)\cong H^*_\circ(\M)\otimes H^*(\J),$$
to get\footnote{Recall that $\tau_S=\prod_{i\in S} \tau_i \in H^*(\J)$.}
\begin{eqnarray} c^\circ_{2g+2d-2}(\bM^\circ_d)=\sum_{S\subset\{1\dots
2g\}}\eqze_{d,S}\otimes \tau_S.\label{chernmumford} \end{eqnarray}
In this way, for any
$S\subset\{1\dots 2g\}$
we get ${\eqze_{d,S}}\in 
\Qu[{\eqal},\eqbe,\eqpsi_i]$ 
of degree $2g+2d-2-\deg(\tau_S)$. When $d=1$ we also express 
$$c^\circ_{2g+r}(\bM^\circ_1)=\sum_{S\subset\{1\dots
2g\}}\zeta^{\circ r}_{S}\otimes \tau_S,$$
to get \begin{eqnarray}
\zeta^{\circ r}_{S}\in \Qu[\eqal,\eqbe,\eqpsi]\label{eqmumfordrelation}\end{eqnarray} for each
$S\subset \{1\dots 2g\}$, which we call the {\em equivariant
Mumford relations}\footnote{The name is justified by noting that
forgetting the $U(1)$-equivariant structure $\zeta^{\circ r}_{S}$ goes
to the Mumford relation $\zeta^r_S$, defined in (\ref{mumfordrelation}).}. 

Now we prove the following fundamental proposition.

\begin{proposition} 
For any $S\subset \{1\dots 2g\}$ and its complement
$S^\prime=\{1\dots 2g\}\setminus 
S$ we have 
\begin{eqnarray} {\eqze_{d,S^\prime}}\mid_{\widetilde{U}_d}=\teqed\cdot \xi_{S}, \label{restriction}\end{eqnarray} 
where $\teqed\in H^*_\circ(\widetilde{U}_d)$ 
is the equivariant Euler class of the normal bundle of 
$\widetilde{U}_d$ in $\tM$ and $$\xi_S=\prod_{i\in S} \xi_i\in
H^*(\Sigma_\bd)\subset H_\circ^*(\widetilde{F}_d)\cong H^*_\circ(\widetilde{U}_d).$$  
Equivalently \begin{eqnarray} (i_d)_*(\xi_{S})=
{\eqze_{d,S^\prime}}\mid_{\widetilde{U}_{\leq d}},
\label{pushforward}\end{eqnarray}
where $i_d:\widetilde{U}_d\to \widetilde{U}_{\leq d}$ is the embedding. 
\label{fundamentalproposition}
\end{proposition}

\begin{proof} 
Recall the equivariant virtual Mumford bundle $\bM^\circ_d$ from the
previous section, and restrict it to $\widetilde{U}_{\leq d}$. A
simple calculation gives that $\ch_0(\bM^\circ_d)=2g+2d-3$ and in turn
that $$\rank(W^\circ)-\rank(V^\circ)+1=2g+2d-2.$$ According
to Theorem~\ref{descriptionudd} the codimension of the 
degeneracy locus of $f^\circ$ of
Theorem~\ref{degeneracym} is $2g+2d-2$, thus the degeneracy locus
has the expected dimension. 
Corollary~\ref{eqporteous} then yields: 
$$\eta^{\circ \tM}_{UD^d_d}=c^\circ_{2g+2d-2}(W^\circ-V^\circ)=c^\circ_{2g+2d-2}(\bM^\circ_d).$$

Let us denote $Gr_d=UD_d^d\cap \widetilde{F}_d$. Then it is immediate
that 
\begin{eqnarray} 
c^\circ_{2g+2d-2}(\bM^\circ_d)\mid_{\widetilde{F}_d} = \eta^{\circ \tM}_{UD^d_d}\mid_{\widetilde{F}_d}=
\teqed \cdot \eta^{\circ \widetilde{F}_d}_{Gr_d}\in
H^*_\circ (\widetilde{F}_d)\cong H^*_\circ(\widetilde{U}_d).
\label{mar}
\end{eqnarray}

On the other hand  recall that  $F_d\times \J$ is a $2^{2g}$-fold cover of $\widetilde{F}_d$, where
$F_d$ is the moduli space of complexes 
$L\stackrel{s}{\to}L\Lambda K$ with $\deg{L}=d$ and fixed line
bundle $\Lambda$ with $\deg(\Lambda)=1$. We can define a map 
$f:F_d\to \J$ by sending 
$f(L\stackrel{s}{\to}L\Lambda K)=L\otimes L^{1-d}$. 

Let us denote by
abuse of notation 
$Gr_d\subset F_d\times \J$ the pullback of $Gr_d$ to $F_d\times \J$. 
Then Theorem~\ref{descriptionudd} shows that $Gr_d$ is nothing else than the
graph of the map $f$.

Now let $S\subset \{1\dots 2g\}$. Then
clearly $f^*(\tau_S)=\xi_S$ so denoting by 
$\eqed\in H_\circ^*(F_d)\cong H_\circ^*(U_d)$ the equivariant Euler
class of the normal bundle of $U_d$ in $\M$ we get
$$
\begin{array}{cclp{4cm}}
\eqed \cdot \xi_S  &= & \eqed \cdot f^*(\tau_S) & \\ &=&
\eqed \cdot (\pr_{F_d})_*\left( 
\eta^{F_d\times \J}_{Gr_d}\cdot \pr^*_\J(\tau_S)\right) & \mbox{from
(\ref{pullformula})} \\ &=&
(\pr_{F_d})_*\left(\teqed \cdot 
\eta^{F_d\times \J}_{Gr_d}\cdot \pr^*_\J(\tau_S)\right) & \mbox{since
$\pr_{F_d}^*(\eqed)=\teqed$} 
\\ &=&
(\pr_{F_d})_*\left(c^\circ_{2g+2d-2}(\bM^\circ_d)\cdot \pr^*_\J(\tau_S)\right)
& \mbox{from (\ref{mar})} 
 \\ &=&
(\pr_{F_d})_*\left((\sum_{R\subset \{1\dots 2g\}}
{\eqze_{d,R}}\otimes \tau_R )\cdot \pr^*_\J(\tau_S)\right) & \mbox{from (\ref{chernmumford})}\\ & = &
{\eqze_{d,S^\prime}}\mid_{F_d}&
\end{array}
$$
which proves (\ref{restriction}).
\end{proof}

\begin{remark} For $S=\emptyset$, (\ref{pushforward}) 
says that  the cohomology class \begin{eqnarray} \eta^{\circ
\tM}_{\widetilde{U}_d}=(i_d)_*(1)=(i_d)_*(\xi_{\emptyset})=\eqze_{d,\{1\dots 2g\}
}\label{higherstrata} \end{eqnarray} in $$H^{2(g+2d-2)}_\circ(\widetilde{U}_{\leq d})\cong
H^{2(g+2d-2)}_\circ(\tM),$$ since for $k>d$ the stratum $\widetilde{U}_k$ has codimension at
least $2(g+2d)$. In particular for $d=1$ and forgetting the
$U(1)$-equivariant structure, we have that  $$\eta^\tM_{\widetilde{U}_1}=\zeta^0_{\{1\dots
2g\}}\in H^{2g}(\tM),$$ i.e. that the cohomology class of the
first stratum $\widetilde{U}_1$ in $\tM$ agrees with the first Mumford
relation of degree $2g$ ! This is not so surprising if we believe 
the generation theorem --to be proved at the end of the section--
i.e. that the $\eta^\tM_{\widetilde{U}_1}$ cohomology class is expressed
as a (complex) degree $g$ polynomial of the universal classes, since then this
polynomial --$\tN$ and $\widetilde{U}_1$ being disjoint in $\tM$ --
should restrict to $0$ on $\tN$, which therefore should be some
multiple of the first Mumford relation. A similar argument shows that
the cohomology classes (\ref{higherstrata}) of the higher strata are also generated by the
Mumford relations. The next formula gives the exact statement: 
\begin{eqnarray}\sum_{i=0}^\infty
c_i(\bM^\circ_{d})=\left(1+\eqal+\frac{\eqal-\eqbe}{4}\right)^{d-1}
\sum_{i=0}^\infty
c^\circ_i(\bM^\circ_1),\label{exactstatement}\end{eqnarray}
which can be obtained similarly to (\ref{chern}). From this it follows
that the equivariant classes $\eqze_{d,S}$, and in particular the
equivariant cohomology classes of the higher strata, are in the
ideal generated by the equivariant Mumford relations of 
(\ref{eqmumfordrelation}). This fact will 
be used in the proof of Theorem~\ref{mumfordtheorem}. 
\end{remark}

To complete the proof of Proposition~\ref{fundamentalproposition}
we need only to prove the following lemma: 

\begin{lemma} Let $f:X\to Y$ be a map of compact manifolds. Then if
we denote by $Gr\subset X\times Y$ the graph of $f$ and let 
$a\in H^*(Y)$ then we have the formula:
\begin{eqnarray} f^*(a)=(\pr_X)_*\left( \eta^{X\times Y}_{Gr}
\cdot\pr^*_Y(a)\right).\label{pullformula}
\end{eqnarray}
\label{graph}
\end{lemma}

\begin{proof} Let $b\in H^*(X)$. Note that by definition
of the graph $Gr$ \begin{eqnarray}\pr_Y\circ i_{Gr} = 
f\circ \pr_X\circ i_{Gr} : Gr\to Gr. \label{obvious}\end{eqnarray}
and \begin{eqnarray}  \pr_X\circ i_{Gr}:Gr \to X \ \mbox{\rm is
a homeomorphism}, \label{obvious2} \end{eqnarray} where $i_{Gr}:Gr\to
X\times Y$ is the embedding.  
Since $(\pr_X)_*$ is an $H^*(X)$-module homomorphism we have 
$$(\pr_X)_*(\eta^{X\times Y}_{Gr}\cdot \pr^*_Y(a))\cdot b=
(\pr_X)_*(\eta^{X\times Y}_{Gr}\cdot\pr^*_Y(a) \cdot \pr_X^*(b)).$$
Therefore \begin{eqnarray*}
\int_X (\pr_X)_*(\eta^{X\times
Y}_{Gr}\cdot \pr^*_Y(a))\cdot b & = &\int_{X\times Y}\eta^{X\times
Y}_{Gr} \cdot \pr^*_Y(a) \cdot \pr_X^*(b)\\ &=& \int_{Gr}
 i_{Gr}^*(\pr^*_Y(a))\cdot i_{Gr}^*(\pr^*_X(b)) \\ 
& = & \int_{Gr} 
 i_{Gr}^*(\pr_X^*(f^*(a)))\cdot i_{Gr}^*(\pr^*_X(b))\\ & = & \int_X
f^*(a)\cdot b,\end{eqnarray*}
where we used (\ref{obvious}) and (\ref{obvious2}). Now Poincar\'e
duality yields the result.
\end{proof}

Now we have the following important corollary of
Proposition~\ref{fundamentalproposition}: It will yield the
generation theorem at the end of the present section, 
and an obvious generalization of it will provide a
purely geometric proof for the Mumford conjecture in Section~\ref{mumfordproof}.

\begin{corollary}
Set $\calR_0$ to be the $\Qu[u]$-submodule of $H^*_\circ(\tM)$ generated by the
universal classes of Theorem~\ref{eqgeneration}. Furthermore for 
$d=1\dots g-1$ set $$\calR_d=\langle \left\{ {\eqze_{d,S}}:  
 S\subset\{ 1\dots 2g \} \right\} \rangle_{\Qu[u,\eqal,\tau_i]},$$ where
$\langle,\rangle_{\Qu[u,\eqal,\tau_i]}$ stands for the generated $\Qu[u,\eqal,\tau_i]$-module. Then
\begin{eqnarray} 
i^*_{d}(\calR_{d^\prime})=0 \mbox{\ \rm for\ } d<d^\prime \mbox{\
\rm \ and \ }
i_d^*(\calR_d)=\langle \teqed \rangle \subset H^*_\circ(\widetilde{U}_d),\label{generationcondition}\end{eqnarray}  
\label{fundamental}
where $\teqed$ is the equivariant normal bundle to the stratum $\widetilde{U}_d$.
\end{corollary}

\begin{proof} When $d=0$ the statement (\ref{generationcondition}) is
equivalent\footnote{Recall that we set $\widetilde{e}^\circ_0=1$.} 
to the corresponding generation theorem for $\tN=\widetilde{F}_0$, which we
already know. 

For $d>0$ recall that
$H^*_\circ(\widetilde{F}_d)=H^*(\widetilde{F}_d)\otimes \Qu[u]$ is generated by
classes $u\in H^2_\circ(\widetilde{F}_d)$, $\eta\in
H^2_\circ(\widetilde{F}_d)$ and $\tau_i\in H^1_\circ(\widetilde{F}_d)$, $\xi_i\in
H^1_\circ(\widetilde{F}_d)$ for $i=1\dots 2g$. 
Since from Lemma~\ref{restrictions} we have 
$\eqal\mid_{\widetilde{F}_d}=(2d-1)(\eta-u)+\sigma$
it follows easily that  
$$\langle \xi_S: S\subset \{1\dots 2g\}
\rangle_{\Qu[u,\eqal,\tau_i]}=H^*_\circ(\widetilde{F}_d)\cong H^*_\circ(\widetilde{U}_d) .$$
Therefore the statement
(\ref{generationcondition})
follows from Proposition~\ref{fundamentalproposition}.
\end{proof}

The main result of this section is the following corollary:

\begin{theorem} The equivariant cohomology ring $H^*_\circ(\tM)$ is
generated as a $\Qu [u]$-module by the equivariant 
universal classes $\eqal\in H^2_\circ(\tM)$,
$\eqbe\in H^4_\circ(\tM)$ and  $\tau_i\in H^1_\circ(\tM)$,
$\eqpsi_i \in H^3(\tM)$ for $i=1\dots 2g$.
\label{eqgeneration} 
\end{theorem}  

\begin{proof} The proof rests on Proposition~\ref{propositionkirwan}
as explained in Remark 2 after it. Namely we have a strongly
$U(1)$-perfect stratification $\tM=\bigcup_{d=0}^{g-1} \widetilde{U}_d$,
and the sets $\calR_d$ of Corollary~\ref{fundamental}, 
which contain elements generated by the
universal classes, satisfy the conditions of
Proposition~\ref{propositionkirwan}. The result follows. 
\end{proof}

Since the forgetful map $H^*_\circ(\tM)\to H^*(\tM)$ is surjective, we
have the following immediate corollary.

\begin{corollary} The ordinary cohomology ring $H^*(\tM)$ is
generated by the ordinary 
universal classes $\alpha\in H^2(\tM)$,
$\beta\in H^4(\tM)$ and  $\tau_i\in H^1 (\tM)$,
$\psi_i \in H^3(\tM)$ for $i=1\dots 2g$. Consequently the ring $H^*(\M)^\Gamma$
is generated by universal classes:
\label{generation} $\alpha\in H^2(\M)$,
$\beta\in H^4(\M)$ and  
$\psi_i \in H^3(\M)$ for $i=1\dots 2g$.
\end{corollary}

\newpage
\section{A conjectured complete set of relations for $H_I^*(\M)$}
\label{relations}

In the previous chapter we showed that the ring $H^*(\M)^\Gamma$ is
generated by universal classes. To have a complete description of the
ring it is sufficient to give a complete set of relations in these
universal classes. 
The idea is to first determine the relations in $H^*_\circ(\M)^\Gamma$ based
on the injectiveness of (\ref{injective}). Namely we consider
$R_d$ the kernel of the composition of the projection $\Qu
[\eqal,\eqbe,\eqpsi_i]$ to $H^*_\circ(\M)^\Gamma$ with the restriction of  $H_\circ^*(\M)^\Gamma$ to 
$H_\circ^*(F_d)^\Gamma$ for $0\leq d \leq g-1$. Then (\ref{injective})
gives that $\bigcap_{d=0}^{g-1} R_d$ is the ideal of relations of
$H^*_\circ(\M)^\Gamma$, in other words the ring $$H^*_\circ(\M)^\Gamma\cong \Qu
[\eqal,\eqbe,\eqpsi_i]/\bigcap_{d=0}^{g-1} R_d$$ is given by generators
and relations. 
 
Though we could not yet proceed this way, computer calculations with
the software package\footnote{Cf. http://www.math.uiuc.edu/Macaulay2/.}   
Macaulay 2
gave us enough numerical evidence to be able to formulate a conjecture
about a complete
description of the subring of $H^*(\tM)^\Gamma$ 
generated by $\alpha_\M$,$\beta_\M$ and
$\gamma_\M$. 

To explain this conjecture we have to define certain polynomials:

\begin{definition} Define polynomials in $\Qu [\alpha,\beta,
\gamma]$ as follows $$\rho_{r,s,t}=\sum_{i=0}^{\min(r,s)}
\binomial{r}{i}\binomial{g-t-i}{g-t-s} \al^{r-i} \be^{s-i} (2\ga)^{t+i},$$
for $r,s,t\geq 0$. 

Let $R$ denote the graded ring given by generators $\al,\be,\ga$ of
degree $2$, $4$ and $6$ respectively, and
relations:  $$\rho_{r,s,t} \mbox{\ \ \rm for \ } r,s,t\geq 0 \mbox{\rm
\ a \ } r+3s+3t> 3g-3.$$
\end{definition}

\begin{conjecture} The subring $H^*_I(\M)\subset H^*(\M)^\Gamma$, generated
by $\alpha_\M, \beta_\M $ and $\gamma_\M$ is isomorphic to $R$. 
\label{presentation}
\end{conjecture}

\begin{remark}1. Computer calculations with Macaulay 2 show that the
statement is true for $2\leq g \leq 7$.

2. We used the notation $H^*_I(\M)$ for the subring generated by
$\al_\M,\be_\M$ and $\ga_\M$ because it can be thought of as the invariant
subring of the action of $Sp(2g,\Z)$ on $H^*(\M)^\Gamma$, which is
given abstractly by letting $Sp(2g,\Z)$ act on $H^3(\M)^\Gamma$ in the
usual symplectic manner, which induces an action\footnote{As a matter of fact this action 
is induced from the action of the
mapping class group of $\Sigma$ on the cohomology of the
representation space of $\pi(\Sigma)$ to $SL(2,\C)$, but we do not
need this fact, so we simply think of this action as given abstractly.}  on the whole of
$H^*(\M)^\Gamma$, because of Corollary~\ref{generation}. 
\end{remark}

We give some further evidence supporting Conjecture~\ref{presentation}
in the rest of the section. First we show that $H^*_I(\M)$ and the
conjectured ring have the same Poincar\'e polynomial.

\begin{lemma} An additive basis for the ring $R$ is given by
$\al^r\be^s\ga^t$ for $r,s,t\geq 0$ and $r+3s+3t\leq 3g-3$.
Consequently its Poincar\'e polynomial is of the form 
\begin{eqnarray}
 \sum_\stack{r,s,t \geq 0}{r + 3s + 3t \leq 3g-3}
T^{r+2s+3t}.
\label{poincareofR}
\end{eqnarray}
\end{lemma}

\begin{proof} The result easily follows by noting that 
the highest order term of $\rho_{r,s,t}$ in the
lexicographical ordering is clearly $\al^r\be^s\ga^t$.  
\end{proof}

\begin{theorem} The Poincar\'e polynomial
$P_T^I(\M)$ of $H^*_I(\M)$ equals (\ref{poincareofR}). 
\end{theorem} 

\begin{proof} First we define an $Sp(2g,\Z)$-action on $H^*(\Sigma_n)$
induced by the usual symplectic action on $H^1(\Sigma_n)$. This gives
$H^*(F_d)^\Gamma$ an $Sp(2g,\Z)$-module structure. Observe that
$H^*(\M_{\leq d})^\Gamma$ being generated\footnote{Just like
Corollary~\ref{generation} it follows from Corollary~\ref{fundamental}.} by tautological classes
admits a natural $Sp(2g,\Z)$-module structure such that the surjective
map $$i_{\leq d}^*:H^*(\M)^\Gamma\to H^*(\M_{\leq d})^\Gamma$$ is an
$Sp(2g,\Z)$-module homomorphism. Moreover from (\ref{pushforward}) it
easily follows that $$(i_d)_*: H^*(U_d)^\Gamma\to H^*(\M_{\leq d})^\Gamma$$ is an
$Sp(2g,\Z)$-module homomorphism. It follows that the short exact
sequence (\ref{gammashortexact}) is a sequence of $Sp(2g,\Z)$-modules
yielding the short exact sequence: 
$$0\to H^*_I(U_d)\to H_I^*(\M_{\leq d})\to H_I^*(\M_{<d})\to 0,$$
where $H^*_I$ denotes the $Sp(2g,\Z)$-invariant part of
$H^*()^\Gamma$. 
Recall $P_T^I(U_d)=P_T^I(\Sigma_\bd)$ from (\ref{invpoinsigman}) and 
$P_T^I(\N)$ from (\ref{invpoincareofN}). We have 
\begin{eqnarray*} P_T^I(\M) & = & 
P_T^I(\N)+\sum_{d=1}^{g-1} T^{g+2d-2} P^I_T(U_d) \\ & = &
 \sum_\stack{r,s,t\geq 0}{r+s+t\leq g-1} T^{r+2s+3t} 
 + \sum_{d=1}^{g-1} T^{g+2d-2}\sum_\stack{q,s\geq 0}{q+2s\leq \bd}
 T^{q+s} 
\\ & = &
 \sum_\stack{r,s,t\geq 0}{r+s+t\leq g-1} T^{r+2s+3t} +  
 \sum_{t=0}^{g-2}T^{3t}\sum_\stack{q,s\geq 0}{q+2s\leq \bd}T^{g-t+q+s}
 \\ &=& 
 \sum_\stack{r,s,t\geq 0}{r+s+t\leq g-1} T^{r+2s+3t}  + 
\sum_{t=0}^{g-2}T^{3t}\sum_\stack{r,s\geq 0, r+s\geq g-t}{r+3s\leq
 3g-3-3t} 
T^{r+2s} \\ & = & 
 \sum_\stack{r,s,t \geq 0}{ r + 3s + 3t \leq 3g-3}
T^{r+2s+3t}.
\end{eqnarray*}

We first introduced $t=d-1$ and then $r=q-s+g-t$. The result
follows.  
\end{proof}

\begin{theorem} The polynomial 
$\rho_{0,g,0}=\beta^g$ of complex degree $2g$ is zero in $H^*(\M)^\Gamma$. 
\end{theorem}

\begin{proof} In order to prove such a statement we want to extend it to an
equivariant relation. Namely we prove that $$\sum_{r=0}^g \eqze_{r,g-r} u^r=0$$
on $\M$. By (\ref{injective}) it is sufficient to check this relation 
on $F_d$ for every $0\leq d \leq g-1$. For $d=0$ the vanishing is
automatic since Zagier's relations (\ref{additivebasis}) hold on
$\N=F_0$. 
For $d>0$ we use
Lemma~\ref{restrictions} and substitute $x=u$ and $y=1$ in Zagier's
generating function (\ref{generatingxirs}) to get: 
$$ \sum_{r=0}^g \eqze_{r,g-r} u^r\mid_{F_d}= 
\left(e^{\si u}
\frac{(1-(\eta-u)(\eta-2u))^{d-1}}{(1-\eta(\eta-u))^d}\right)_{2g}.$$
Recall from Subsection~\ref{symmetric} that $(\dots)_m$ means the
parts of degree $2m$. 
To prove that it  vanishes in $H^*(F_d)^\Gamma \cong
H^*(\Sigma_{\bd})$, express it as 
\begin{eqnarray*}
e^{\si u} \frac{(1-(\eta-u)(\eta-2u))^{d-1}}
{(1+\eta u)^d \left( 1-\frac{\eta^2}{1+\eta u} \right)^d }
\end{eqnarray*}
which equals
\begin{eqnarray}
\sum_{i=1}^\infty \binomial{d+i}{i} 
\frac{\eta^{2i} e^{\si u}}{(1 + \eta u)^{d+i}} 
(1-(\eta-u)(\eta-2u))^{d-1}.
\label{sum}
\end{eqnarray}               

It immediately follows from Lemma~\ref{symmetricvanish} that  
$$\left( \frac{e^{\si u} (\eta u)^{2k}}{(1+(\eta u))^{d+k}} 
\right)_{2(g-d+k+1)+n} = 0$$ for $n \geq 0$
and hence that
$$\left( \frac{e^{\si u} \eta^{2k}}{(1+(\eta u))^{d+k}} 
\right)_{2g-2d+2+n} = 0$$ for $n \geq 0$.  
Consequently in (\ref{sum}) each term vanishes at total degree $2g$. 

The result follows. 
\end{proof}

\begin{remark} 1. A similar argument shows the vanishing of the first relation 
$$\rho_{1,g-1,0}=g\eqal(\eqbe)^{g-1}+(g-1)(\eqbe)^{g-2}(2\eqga)=\eqze_{1,g-1}$$ of complex degree
$2g-1$, by showing the vanishing of the equivariant class 
$$\sum_{r=1}^g r \eqze_{r,g-r} u^{r-1} = 0$$ on $\M$. This goes
similarly to the argument above, though the calculation is more
tedious.

2. To settle Conjecture~\ref{presentation} 'all' we have to do is to
extend these polynomials to equivariant relations as we did
above in a special case. 
However it seems that our polynomials are simpler than the ones
which may occur naturally from equivariant relations. The reason for
this may be that $\rho_{r,s,t}$ can be defined in a much 
simpler fashion, than Zagier's 
polynomials $\zeta_{r,s,t}$, which do have some geometric origin. 

To picture this difference 
consider a conjectured relation $\rho_{r,s,t}$ of 
Conjecture~\ref{presentation}. Since it should also be a relation on
$\N$ it has to be written as a linear combination of Zagier's
relations (\ref{additivebasis}). Computer calculations show that this
is indeed the case, though the linear combinations tend to be fairly
complicated. We do not even have a general formula for this linear
combination! 
\end{remark}

\chapter{Resolution tower}
\label{minfty}

By considering the spaces $\tM_k$ as defined in
Definition~\ref{poles}, we complete here\footnote{ This chapter is
based on a
joint work with Michael Thaddeus.} the picture of the cohomology ring
$H^*(\tM)$ in the framework described at the beginning of
Section~\ref{cohomologyofmoduli}. 
Namely we show that the tower $$\tM\cong \tM_0\subset
\tM_1 \subset \dots \subset \tM_k \subset \dots$$ gives a resolution
of the cohomology ring $H^*(\tM)$, in the sense that
$$i_0^*:H^*(\tM_\infty)\to H^*(\tM)$$ is a free graded commutative
resolution of $H^*(\tM)$, where $$\tM_\infty=\lim_{\to} \tM_k$$ is defined
as the direct limit of the above tower. 

We also show that the space $\tM_\infty$ is important on its own right. To
prove this, we use it to give a simple geometric proof of the
Mumford conjecture in Section~\ref{mumfordproof}. 
In the last section we explain why $\tM_\infty$
is so useful by showing that it is a model for the
classifying space of $\bG$, the gauge group modulo constant scalars. 
Also we prove that this homotopy
equivalence is even preserved on the level of the strata. We conclude by
showing that the above tower has the property that its
homotopy groups are stabilizing, thus getting a picture similar to the
Atiyah-Jones theorem for the moduli space of instantons on $S^4$.

Finally we mention that $\tM_k$ and even $\tM_\infty$ appeared already
in the
work of Donagi and Markman \cite{donagi-markman}, where they showed that
$\tM_k$ is a complex Poisson manifold, and its Hitchin map is an algebraic
completely integrable Hamiltonian system.

\newpage
\section{The moduli space of Higgs $k$-bundles $\tM_k$}

In this section we list some basic properties of the spaces
$\tM_k$. They are completely analogous to the properties of $\tM$. As 
the proofs are also following the same lines we do not spell out the
details here, but hope, that in case of any doubt, 
the reader can complete the arguments.

\subsubsection{The $\C^*$-action on $\M_k$}
Recall from Subsection~\ref{nagysrac} that $\M_k$ is a  smooth
quasi-projective varieties of dimension $6g-6+3k$. Moreover $\C^*$
acts on $\M_k$ by multiplication of the Higgs $k$-field. Completely
analogously\footnote{This follows from the gauge
theory construction of $\M_k$ of Section~\ref{gaugemk}.} to $\M$, 
$\M_k$ is a K\"ahler manifold, and 
$U(1)\subset \C^*$ acts on it in a Hamiltonian way, with proper moment
map $\mu_k$ of absolute minimum $0$. Thus as in Subsection~\ref{rizsa1}
we have a stratification $\M_k=\bigcup_{d=0}^{n} U_d^k$ with
upward flows, where $n$ is the number of components of the fixed point
set, to be determined later. We call this the {\em Hitchin
stratification}. Just as in Subsection~\ref{stratificationM} we have
the Shatz stratification $\M_k=\bigcup_{d=0}^{n} U_d^{\prime
k}$, defined by $U_d^{\prime k}=\{(E,\Phi_k): E\in \calC_d\}$, 
which coincides with the Hitchin stratification. 

From Subsection~\ref{rizsa1} we know
that $U_d^k$ retracts to $F_d^k$, the $d$-th component of the fixed
point set of the $U(1)$-action. A stable Higgs $k$-pair $(E,\Phi_k)$ 
is fixed by the circle action either if $E$ is stable and $\Phi_k=0$
or if $$E=L\oplus L^{-1}\Lambda$$ and
$$\Phi_k=\left(\begin{array}{cc} 0 & 0\\ \phi_k & 0
\end{array}\right),$$
where $0\not=\phi_k\in H^*(\Sigma; L^{-2}\Lambda K L_p^k)$. From the
stability of the pair we have that $\deg(L)>0$, and from the
assumption $\phi_k\not=0$ that $\deg(L)\leq g-1+k$. It follows that
$n=g-1+k$ and 
the components of the fixed point set of the $U(1)$-action are $F_0^k\cong \N$
and $F_d^k$ for $0<d\leq g-1+k$ are $2^{2g}$-fold covers of the
symmetric product $\Sigma_{\bd+k}$, with covering group $\Gamma$. Now
the tangent space of $\M_k$ at a point $(E,\Phi_k)\in F_d^k$ is naturally
$$\Hy^1(\Sigma; \End_0(E)\stackrel{[\Phi_k,\cdot]}{\longrightarrow}
\End_0(E)\otimes K\otimes L_p^k).$$ Tracing back the action of $\C^*$
on it one gets that the only negative weight appearing is $-1$ and the
corresponding weight space is $\Hy^1(\Sigma; L^{-2}\Lambda\to 0)\cong
H^1(\Sigma; L^{-2}\Lambda)$. By Riemann-Roch it has dimension $g+2d-2$. Thus
the real codimension of $U_d^k$ in $\M_k$ which is the same as the
index of the critical submanifold $F_d^k$ is $2(g+2d-2)$. 

\subsubsection{Poincar\'e polynomial of $\tM_k$}
The fact that the indices are even implies that 
the stratification is perfect, thus we have the following formula for
the $\Gamma$-invariant Poincar\'e polynomial of $\M_k$:
$$P_t(\M_k)^\Gamma=P_t(\N)+\sum_{d=1}^{g-1+k} t^{2(g+2d-2)}P_t(\Sigma_{\bd+k}),$$
from which we get the following formula for the Poincar\'e polynomial
of $\tM_k$:
\begin{eqnarray}
P_t(\tM_k)=P_t(\J)P_t(\M_k)^\Gamma=P_t(\tN)+\sum_{d=1}^{g-1+k}
t^{2(g+2d-2)}P_t(\Sigma_{\bd+k})P_t(\J).
\label{poincaremk}\end{eqnarray}

\subsubsection{Generators for $H^*(\tM_k)$} We have an equivariant universal
bundle $\E^\circ_{\M_k}$, which gives equivariant universal classes
$\eqal$, $\eqbe$ and $\eqpsi_i$ in $H^*_\circ(\M_k)^\Gamma$. We
have now the analogue of
Corollary~\ref{fundamental}:
\begin{proposition}
Set $\calR_0$ to be the $\Qu[u,\tau_i]$-submodule of $H^*_\circ(\tM_k)$ generated by the
above equivariant 
universal classes. Furthermore for 
$d=1\dots g-1+k$ set $$\calR_d=\langle  {\eqze_{d,S}}: S\subset\{
1\dots 2g \} \rangle_{\Qu[u,\eqal,\tau_i]},$$ 
where
$\langle,\rangle_{\Qu[u,\eqal,\tau_i]}$ stands for the generated
${\Qu[u,\eqal,\tau_i]}$-module. 
Then $i^*_{d}(\calR_{d^\prime})=0$ for $d<d^\prime$ and
\begin{eqnarray}
i_d^*(\calR_d)=\langle e_d \rangle \subset H^*_\circ(\widetilde{U}^k_d).\end{eqnarray}  
\label{fundamentalk}
\end{proposition} 
It follows from Proposition~\ref{propositionkirwan} that the
equivariant cohomology ring of $H^*(\tM_k)$ 
is generated as an algebra by $u$ and universal classes
$\eqal$,$\eqbe$, $\eqpsi_i$ and $\tau_i$.

\newpage
\section{The moduli space of Higgs $\infty$-bundles $\tM_\infty$}
\label{moduliminfty}

Let us fix $s_p$, a non-zero section of $L_p$. This induces embeddings
$i_k:\tM_k\to \tM_{k+1}$ given by $i_k(E,\Phi_k)=(E,\Phi_k\otimes
s_p)$. It clearly respects the $\C^*$-action and 
$i_k(\widetilde{N}_d^k)\subset \widetilde{N}_d^{k+1}$ which for $d>0$ is induced from the
map $\Sigma_{\bd+k}\to \Sigma_{\bd+k+1}$ given by $D\mapsto D+p$. It
follows from (\ref{surjectsymmetric}) and Corollary~\ref{surjection} 
that \begin{eqnarray} i_k^*:H^*(\tM_k)\to
H^*(\tM_{k+1}) \mbox{\rm\ is a surjection.}\label{surjectmk}\end{eqnarray} 

Now consider the direct limit of the embeddings $i_k$, and denote it
by $$\tM_\infty=\lim_{\longrightarrow} \tM_k. \label{Minfty}$$ 
Then we have the inverse limit $$H^*(\tM_\infty)=\lim_{\longleftarrow}
H^*(\tM_k),$$
since $H^*$ is a contravariant functor. Recall $\bG$ and $P_t(B\bG)$ 
from (\ref{poincarebg}). From (\ref{surjectmk}) we have
that 
\begin{eqnarray} P_t(\tM_\infty) & = & \lim_{k\to \infty} P_t(\tM_k)= 
\lim_{k\to \infty} \left( P_t(\tN)+\sum_{d=1}^{g-1+k}
t^{2(g+2d-2)}P_t(\Sigma_{\bd+k})P_t(\J)\right) \nonumber\\  
& = & 
P_t(\tN)+\sum_{d=1}^{\infty}
t^{2(g+2d-2)}P_t(\J) \lim_{k\to
\infty}\left(P_t(\Sigma_{\bd+k})\right)\nonumber\\ 
& = & P_t(\tN)+\sum_{d=1}^{\infty}
t^{2(g+2d-2)}P_t(\J) P_t(\Sigma_\infty) \nonumber\\ & = & 
(1+t)^{2g}\left(\frac{(1+t^3)^{2g}-t^{2g}(1+t)^{2g}}{(1-t^2)(1-t^4)}\right)+
\sum_{d=1}^{\infty}
t^{2(g+2d-2)}(1+t)^{2g}\frac{(1+t)^{2g}}{(1-t^2)}\nonumber\\ & = & 
\frac{\left\{(1+t)(1+t^3)\right\}^{2g}}{(1-t^2)(1-t^4)} =  P_t(B\bG) 
\label{poincareminfty}
\end{eqnarray}
On the other hand $H^*(\tM_\infty)$ is generated by universal classes,
because the same is true for $H^*(\tM_k)$. It follows that
$H^*(\tM_\infty)$ is a free graded commutative algebra, and thus
$$H^*(\tM_\infty)\to H^*(\tM)$$ is a resolution of the cohomology ring
$H^*(\tM)$. It shows that $H^*(\tM)$ can be understood in the
framework described at the beginning of
Section~\ref{cohomologyofmoduli}.

In the rest of the chapter we give some applications of $\tM_\infty$ to
emphasize its significance.

\newpage
\section{Geometric proof of the Mumford conjecture} 
\label{mumfordproof}

As the embeddings $i_k:\tM_k\to \tM_{k+1}$ are respecting the
$U(1)$-action, 
we have a  $U(1)$-action on $\tM_\infty$. Just like above
$H_\circ^*(\tM_\infty)$ is also generated over $\Qu [u]$ 
by the universal equivariant
classes. Furthermore a calculation, completely analogous to
(\ref{poincareminfty}), shows that
$$P_t^\circ(\tM_\infty)=\frac{P_t(B\bG)}{1-t^2}$$ from which we see that
$H^*_\circ(\tM_\infty)$ is a free graded commutative algebra on the
equivariant universal classes and $u$. Observe also that the
stratification $$\tM_\infty=\bigcup_{d=0}^\infty \widetilde{U}_d^\infty$$
is $U(1)$-perfect, so we are in a position to apply
Proposition~\ref{propositionkirwan} as explained in Remark 1 after
it. Namely Proposition~\ref{fundamentalk} in the direct limit yields:
\begin{proposition}
Set $\calR_0$ to be the subring of $H^*_\circ(\tM_\infty)$ generated by the
universal classes. Furthermore for 
$d\geq 1$ set $$\calR_d=\langle \left\{ {\eqze_{d,S}}: 
S\subset\{ 1\dots 2g \} \right\} \rangle_{\Qu[u,\eqal,\tau_i]},$$ where
$\langle,\rangle_{\Qu[u,\eqal,\tau_i]}$ stands for the generated
$\Qu[u,\eqal,\tau_i]$-module. 
Then $i^*_{d}(\calR_{d^\prime})=0$ for $d<d^\prime$ and
\begin{eqnarray}
i_d^*(\calR_d)=\langle e_d \rangle \subset H^*_\circ(\widetilde{U}_d^\infty) \end{eqnarray}  
\end{proposition} 

Now as explained in Remark 1 after
Proposition~\ref{propositionkirwan}, 
the above proposition yields that
$\bigcup_{d=1}^\infty\calR_d$ additively generates the kernel of
$$H_\circ^*(\tM_\infty)\to H_\circ^*(\tN).$$ However $H^*(\tM_\infty)$
is a free graded commutative algebra, thus
$\bigcup_{d=1}^\infty\calR_d$ is a complete set of relations for
$\tN$. Moreover combining with (\ref{exactstatement}), 
we have the following theorem:
\begin{theorem} The Mumford relations $\ze^r_S$ for each
$S\subset \{1\dots 2g\}$ and $r\geq 0$ generate 
the relation
ideal of $H^*(\tN)$.
\label{mumfordtheorem}
\end{theorem}

\newpage
\section{Gauge theoretic construction of $\tM_\infty$}
\label{gaugeminfty}

To construct $\tM_\infty$ gauge theoretically, first 
recall the gauge theoretic construction of $\tM_k$ from Subsection~\ref{gaugemk}. 
Recall also that $s_p\in H^0(\Sigma;L_p)$ is a fixed holomorphic
section of $L_p$. It follows that
there are embeddings $\Omega^{1,0}_k\subset \Omega^{1,0}_{k+1}$ and
$\Omega^{1,1}_k\subset \Omega^{1,1}_{k+1}$, given by tensoring with $s_p$.
 Since $s_p$ is holomorphic
$\dbar_{k+1}\mid_{\calC\times \Omega_k}=\dbar_k$, consequently
\begin{eqnarray}\calB_k\subset
\calB_{k+1}.\label{embeddingofbk}\end{eqnarray} 
Thus if we define
the direct limit $$\Omega^{1,i}_\infty=\lim_{k\to \infty} \Omega^{1,i}_k,$$ then the direct limit
of the maps $\dbar_k$ will be $$\dbar_\infty:\calC\times
\Omega^{1,0}_\infty\to \Omega^{1,1}_\infty,$$ and the direct limit 
$$\calB_\infty=\lim_{k\to \infty} \calB_k$$ coincides with $\dbar_\infty^{-1}(0)$,
which is the space of pairs $(E,\Phi)$ where $\Phi$ is a holomorphic
Higgs  $\infty$-field.  
We denote by $\pr_\infty:\calB_\infty\to
\calC$, the direct limit of $\pr_k$. We also have $i_k((\calB_k)^s)\subset
(\calB_{k+1})^s$, and we let $$(\calB_\infty)^s=\lim_{\longrightarrow}(\calB_k)^s$$ denote
the space of stable Higgs $\infty$-bundles. 

It follows from the foregoing that we can think of
$\tM_\infty$ as the quotient $(\calB_\infty)^s/\G^c$.  
In order to apply this construction to obtain topological results
about $\tM_\infty$, we make a detailed study of the spaces occurring:

\paragraph{A condition for $\dbar_k$ to be a submersion.}
The map $\dbar_k$ of (\ref{dbark}) is a smooth  
map of Banach manifolds  and the following theorem gives a sufficient
condition for the derivative $T_{\dbar_k}$ to be surjective:

\begin{theorem} The derivative
$T_{\dbar_k}$ is surjective at the point $(E,\Phi)\in
\calC\times\Omega_k^{1,0}$ if and only if 
\begin{eqnarray}\Hy^0\left(\Sigma;\End(E)\otimes
L^{-k}_p\stackrel{[\Phi,\cdot]}{\longrightarrow}\End(E)\otimes
L_p^{-k} \otimes K\right)\label{hypervanish} \end{eqnarray} is trivial. 
\label{derivative}
\end{theorem}  

\begin{proof} At the point $(E,\Phi)$ the derivative of $\dbar_k$ 
$$T_{\dbar_k}:\Omega^{0,1}_k\times 
\Omega^{1,0}_k\to
\Omega^{1,1}_k$$
is given by $$T_{\dbar_k}(\alpha,\beta)=\dbar^E_k\beta +[\alpha,\Phi],$$
where $$\alpha\in \Omega^{0,1}_k \mbox{\rm \ and \ } \beta\in
\Omega^{1,0}_k.$$ There is a natural non-degenerate pairing between $\Omega_k^{1,1}$
and $\Omega^{0,0}_{-k}$ given by integrating
over $\Sigma$ the trace of the tensor product. Suppose now that
$\psi\in  \Omega^{0,0}_{-k}$ is perpendicular
to the image of $T_{\dbar_k}$,
i.e. \begin{eqnarray}\int_{\Sigma}\trace\left(T_{\dbar_k}(\alpha,\beta)\otimes
\psi\right)=0,\label{zero}\end{eqnarray}
for all $\alpha$ and $\beta$. Then for all $\beta\in \Omega_k^{1,0}$ 
\begin{eqnarray*}
\int_\Sigma \trace\left(\beta\otimes\dbar^E_{-k}\psi\right)&=&
\int_\Sigma
\trace\left(\dbar_0^E(\beta\otimes\psi)-\dbar^E_k(\beta)\otimes\psi\right)\\
&=&\int_\Sigma 	d_0^E\trace(\beta\otimes\psi)-\int_\Sigma \trace
\left( \dbar^E_k(\beta)\otimes\psi\right)\\ & = & 0,
\end{eqnarray*}
the first term vanishes because of Stokes' theorem, the second
because of (\ref{zero}) for the choice of $\alpha=0$. However the
pairing between $\Omega_{k}^{1,0}$ and $\Omega_{-k}^{1,0}$ is
non-degenerate, which gives that $\dbar^E_{-k}(\psi)= 0$. 
On the other hand we have for all $\alpha\in \Omega_{k}^{1,0}$ that $$\int_{\Sigma} [\al,\Phi]\otimes
\psi=0$$ from (\ref{zero}) for the choice of $\beta=0$. It follows
that $[\psi,\Phi]=0$.

Putting everything together we have that $T_{\dbar_k}$ is surjective
at $(E,\Phi)$ if and only if $\dbar_{-k}^E(\psi)=0$ and
$[\psi,\Phi]=0$ imply $\psi=0$. 
However this is exactly the Dolbeault description of
the hypercohomology vector space (\ref{hypervanish}). The result follows.
\end{proof}

The following lemma will be useful later:

\begin{lemma} If $k>0$ and $(E,\Phi)$ is a stable Higgs $k$-bundle, 
then the hypercohomology (\ref{hypervanish}) vanishes and thus $T_{\dbar_k}$ is surjective.

If $k>0$, $0\leq 2d\leq k$ and $E\in \calC_d$, then
$H^0\left(\Sigma; \End(E)\otimes L^{-k}_p\right)=0$. Consequently the
hypercohomology (\ref{hypervanish}) vanishes, and thus
$T_{\dbar_k}$ is surjective at $(E,\Phi)$ for any $\Phi$.
\label{2d}
\end{lemma}
\begin{proof} The first statement follows since $(E\otimes L_p^{-k},\Phi)$
is also stable, thus a
result analogous to Theorem~\ref{van} for $k$-Higgs bundles for $k>0$
gives the vanishing of the hypercohomology in question. 

For the second part consider $$0\to L \to E\to V \to 0$$ the Harder-Narasimhan
filtration of $E$. Recall from p. 566 of \cite{atiyah-bott} that
$\End^\prime(E)$ denotes the bundle of those endomorphisms which
preserve this filtration. Any filtration-preserving endomorphism of
$E$ gives an element in $\Hom(L,L)\cong \calO_\Sigma$, thus we have a bundle homomorphism
$\End^\prime(E)\to \Hom(L,L)$, whose kernel consists of endomorphisms which
kill $L$ i.e. $V^*\otimes E\subset E^*\otimes E= \End(E)$. Thus we
have the short exact sequence \begin{eqnarray}0\to V^*\otimes E\to
\End^{\prime}(E)\to \calO_\Sigma \to 0.
\label{endprime}\end{eqnarray} Therefore
we have the following filtration of $\End(E)$:
$$0\subset V^*\otimes L\subset V^*\otimes E\subset
\End^\prime(E)\subset \End(E).$$ This is not yet the Harder-Narasimhan
filtration of $\End(E)$ since 
$$\deg\left( (V^*\otimes E)/ (V^*\otimes
L)\right)=\deg
\left(\left(\End^\prime(E)\right)/(V^*\otimes E)\right)=0.$$ However
this means that $\left(\End^\prime(E)\right)/(V^*\otimes
L)$ is semistable, thus the Harder-Narasimhan filtration of $\End(E)$
is: $$0\subset  V^*\otimes L\subset
\End^\prime(E)\subset \End(E),$$ consequently the highest degree line
subbundle of $\End(E)$ is $V^*\otimes L$ of degree
$2d-1$. Therefore $\End(E)\otimes L_p^{-k}$ has highest degree line
subbundle $V^*\otimes L\otimes L_p^{-k}$ of degree $2d-1-k<0$. Such a
bundle cannot have a section. The result follows.
\end{proof}

\paragraph{Stratifications on $\calB_k$.}
We define two stratifications on the spaces $\calB_k$. The first is the
preimage of the Shatz stratification:
$(\calB_k)_d=\pr_k^{-1}(\calC_d)$, i.e. $(\calB_k)_d$ contains pairs
$(E,\Phi)$ with $E\in \calC_d$. The other is given by the
Harder-Narasimhan filtration of Corollary~\ref{higgsfiltrationtheorem}
for Higgs $k$-bundles, namely we define
$(\calB_k)^0\subset \calB_k$ to be the subspace of stable pairs
$(E,\Phi)$, and $(\calB_k)^d$ to be the subspace of pairs $(E,\Phi)$
with destabilizing Higgs $k$-subbundle of degree $d>0$. 

We also let
$(\calB_k)^d_l=(\calB_k)^d\cap (\calB_k)_l$. It is clear that
$(\calB_k)^d\subset (\calB_k)_d$ for $d>0$, because the line bundle
of the destabilizing Higgs $k$-subbundle of a Higgs $k$-bundle $(E,\Phi)$
will be the destabilizing line
bundle of $E$. Thus for $d>0$ we have \begin{eqnarray} {\rm either} \  
(\calB_k)^d_l=\emptyset \ {\rm for} \  d\neq l \ {\rm or} \
(\calB_k)^d_d=(\calB_k)^d.
\label{double}\end{eqnarray}

Since these stratifications are compatible with the embeddings
(\ref{embeddingofbk}) we get stratifications $(\calB_\infty)_d$ and
$(\calB_\infty)^d$ in the direct limit. 

Now we have, analogously to (7.8) of \cite{atiyah-bott}:

\begin{theorem} The decomposition $\calB_\infty=\bigcup_{d=0}^\infty
(\calB_\infty)^d$ has the property:
$$\overline{(\calB_\infty)^d}\subset \bigcup_{i=d}^\infty(\calB_\infty)^i.$$
\label{closed} 
\end{theorem}

\begin{proof} First we show that $(\calB_\infty)^0\subset \calB_\infty$ is
open. For $k>0$ Lemma~\ref{2d} shows that $(\calB_k)^0$ is a Banach
submanifold of $\calB_k$, moreover the tangent space of $(\calB_k)^0$ 
is naturally isomorphic to the tangent space of $\calB_k$, which
proves that $(\calB_k)^0\subset \calB_k$ is open indeed. 

Now if $x\in \overline{(\calB_\infty)^d}$ for $d>0$, then $x\not\in
(\calB_\infty)^0$, since $(\calB_\infty)^0$ is open. However
$$\pr_\infty(x)\in \overline{\calC_d} \subset \bigcup_{i\geq d} \calC_i$$ 
from (7.8) of
\cite{atiyah-bott}, thus (\ref{double}) proves the result. 
\end{proof}

\begin{theorem} For $k>0$ and $0\leq 2l\leq k$ 
the space  $(\calB_k)_{\leq l}$ is naturally a Banach
manifold: it is the total space of a rank $4g-4+4k$ 
smooth, complex vector bundle over
the Banach manifold $\calC_{\leq l}$. 
\label{banach}
\end{theorem}

\begin{proof} For $2l \leq k$ the derivative of
$\dbar_k:\calC_{\leq l}\times \Omega^{1,0}_k\to \Omega^{1,1}_k$ is
surjective by Theorem~\ref{derivative} and Lemma~\ref{2d}. 
Thus the inverse function theorem gives that $
(\calB_k)_{\leq l}=\dbar_k^{-1}(0)$ is a Banach manifold indeed. 
Moreover the fibre of the map $$(\pr_k)_{\leq l}:(\calB_k)_{\leq l}\to \calC_{\leq
l}$$ over the point $E\in \calC_{\leq l}$ is $$H^0(\Sigma;\End(E)\otimes
K\otimes L_p^k)\subset \Omega^{1,0}_k$$ of dimension $4g-4+4k$ since
$$H^1(\Sigma;\End(E)\otimes K \otimes L_p^k)\cong
(H^0(\Sigma;\End(E)\otimes L_p^{-k}))^*=0,$$ by Lemma~\ref{2d}. 
Consequently the map $(\pr_k)_{\leq l}$ is a locally trivial
fibration with fibres $\C^{4g-4+4k}$.
\end{proof}

\begin{corollary} The projection $$\pr_\infty:\calB_\infty\to \calC$$ is a
locally trivial fibration with fibres homeomorphic to $\C^\infty$. 
\label{loctriv}
\end{corollary}

\begin{proof} Since $4g-4+4k\to \infty$ as $k\to \infty $
Theorem~\ref{banach} 
gives that 
$$(\pr_\infty)_{\leq l}:(\calB_\infty)_{\leq l}\to (\calC)_{\leq l}$$ 
is a locally trivial fibration with fibres $\C^\infty$. But clearly
$$\lim_{l\to \infty}(\pr_\infty)_{\leq l}=\pr_\infty,$$ which gives the
desired result. 
\end{proof}

\begin{theorem} If $k>0$, $0\leq 2l \leq k$ and  $d\leq l$ the stratum
$(\calB_k)^d\subset (\calB_k)_{\leq l}$ is 
a  Banach submanifold of $(\calB_k)_{\leq l}$ of complex 
codimension $2g-2+k$.
\label{conormal}
\end{theorem}

\begin{proof} We proceed similarly to the discussion on p. 566 of
\cite{atiyah-bott}. We have that the $\G^c$-orbit of a Higgs $k$-bundle
$\cE=\compE\otimes L_p^k$ in 
$(\calB_k)_{\leq l}$ is, locally, a manifold of finite codimension and
its normal bundle can be identified with the hypercohomology vector space 
$\Hy^1(\Sigma;\End(\cE))$, where
$$\End(\cE)=\End(E)\stackrel{[\Phi,\cdot]}{\longrightarrow}
\End(E)\otimes K \otimes L_p^k$$ is the complex of Higgs $k$-endomorphisms
of $\E$.

In the same way we can identify the normal
bundle to $(\calB_k)^d$. Let $\End^\prime(\cE)$ denote the complex of
Higgs $k$-endomorphisms which respects the
Harder-Narasimhan\footnote{Cf.
Corollary~\ref{higgsfiltrationtheorem}.}
filtration of 
$\cE$ and define the complex $\End^{\prime\prime}(\cE)$  by the exact
sequence \begin{eqnarray}0\to \End^\prime(\cE)\to \End(\cE)\to
\End^{\prime\prime}(\cE)\to 0\label{hypshort}.\end{eqnarray}

Alternatively, one defines $\End^{\prime}(\cE)$ to be the complex
$$\End^{\prime}(\cE)=\End^\prime(E)\stackrel{[\Phi,\cdot]}{\longrightarrow} \End^\prime(E)\otimes K \otimes
L_p^k$$ and $$\End^{\prime\prime}(\cE)=\End^{\prime\prime}(E)\stackrel{[\Phi,\cdot]}{\longrightarrow}
\End^{\prime\prime}(E)\otimes K \otimes L_p^k,$$ using the notation
of p. 566 of \cite{atiyah-bott}. From this definition and 7.4 of \cite{atiyah-bott}
it follows that $$\Hy^0\left(\Sigma;\End^{\prime\prime}(\cE)\right)=0.$$ 
Because of this vanishing, the hypercohomology long exact sequence of the short exact
sequence (\ref{hypshort}) gives the exact sequence:
$$0\to \Hy^1\left(\Sigma;\End^{\prime}(\cE)\right)\to \Hy^1\left(\Sigma;\End(\cE)\right)\to
\Hy^1\left(\Sigma;(\End^{\prime\prime}(\cE)\right)\stackrel{\delta}{\to}
\Hy^2\left(\Sigma;\End^{\prime}(\cE)\right)\to\dots$$

Clearly the conormal to $(\calB_k)^d$ is the factor of
$\Hy^1\left(\Sigma;\End(\cE)\right)$ by  $\Hy^1\left(\Sigma;\End^{\prime}(\cE)\right)$, which by the above
exact sequence is isomorphic to $\ker(\delta)$. 

Now we need the following lemma:

\begin{lemma}
For $k>0$ the vector space
$$\Hy^2\left(\Sigma;\End^{\prime}(\cE)\right)=0$$ is trivial.  
\end{lemma}
\begin{proof}
It is sufficient to show that 
$H^1\left(\Sigma;\End^{\prime}(E)\otimes K \otimes L_p^k\right)=0$. For this we
need by Serre duality that $H^0\left(\Sigma; \left(\End^{\prime}(E)\right)^*\otimes L_p^{-k}\right)=0$.
Taking the dual of (\ref{endprime}) and tensoring by $L_p^{-k}$ we get the short exact
sequence $$0\to L_p^{-k}\to \left(\End^{ \prime}(\cE)\right)^*\otimes
L_p^{-k} \to  V\otimes E^*\otimes
L_p^{-k}\to 0.$$ Since $V\otimes E^*\otimes L_p^{-k}$ has Harder-Narasimhan
filtration: $$0\to L_p^{-k}\to V\otimes E^*\otimes L_p^{-k}\to
V\otimes L^*\otimes
L_p^{-k}\to 0,$$
it follows that $H^0\left(\Sigma;V\otimes E^*\otimes L_p^{-k}\right)=0$, and in turn that 
$H^0\left( \Sigma; \left(\End^{\prime}(E)\right)^*\otimes L_p^{-k}\right)=0.$ The
result follows. \end{proof}
 
The above lemma yields that the conormal to $(\calB_k)^d$ is
isomorphic to $\Hy^1\left(\Sigma;\End^{\prime\prime}(\cE)\right).$ 
Finally we need  the following result:

\begin{lemma} For $k>0$ the dimension of
$\Hy^1\left(\Sigma;\End^{\prime\prime}(\cE)\right)$ depends only on $k$;
it is given by: $$\dim\left(\Hy^1(\Sigma;\End^{\prime\prime}(\cE)\right)=2g-2+k.$$
\end{lemma}
\begin{proof} First we show that
$\Hy^2\left(\Sigma;\End^{\prime\prime}(\cE)\right)=0$. This follows from
$$H^1\left( \Sigma;\End^{\prime\prime}(E)\otimes K \otimes L_p^k\right)=0.$$ Observe that
$\End^{\prime\prime}(E)\cong L^*\otimes V$ thus
$\deg\left(\End^{\prime\prime}(E)\right)=1-2d$, consequently
$$H^0\left( \Sigma;\left(\End^{\prime\prime}(E)\right)^*\otimes
L_p^{-k}\right)=0.$$ Now Riemann-Roch proves  the lemma. 
\end{proof}
Theorem~\ref{conormal} follows. 
\end{proof}

\begin{corollary}
We have $$H_q( (\calB_k)^{\leq d+1},(\calB_k)^{\leq d};\Z)=0$$ for
$2d+2\leq k$ 
and $q < 2(2g-2+k).$  
\label{relativevanish} 
\end{corollary}

\begin{proof}
This is a consequence of Theorem~\ref{closed}, Theorem~\ref{conormal}
and the Thom isomorphism in homology: $$H_q( (\calB_k)^{\leq d+1},(\calB_k)^{\leq
d};\Z)\cong H_{q-2(2g-2+k)}((\calB_k)^{d+1};\Z).$$
\end{proof}

\newpage
\section{Homotopy types}
\label{homotopy}

We learned in Section~\ref{moduliminfty} that the cohomology rings of
$\tM_\infty$ and $B\bG$ are isomorphic. We show in the next subsection
that this is because they are in fact {\em homotopy equivalent}. 

\subsection{Homotopy type of $\tM_\infty$}






We start with a result in the gauge theory setting of the previous section.

\begin{proposition} The space $(\calB_\infty)^0$ is contractible.
\label{weak}
\end{proposition}

\begin{proof} 
First we prove that \begin{eqnarray} \pi_i\left(
(\calB_\infty)^0\right)\ \mbox{is trivial for} \ i\geq 0 
\label{homotopyvanish} \end{eqnarray}
For this we show that 
\begin{eqnarray}
H_i((\calB_\infty)^0;\Z)=\left\{ \begin{array}{ccc} 0 & \mbox{for} &
i>0\\ \Z & \mbox{for} & i=0  \end{array}
\right. \label{homologyvanish}.\end{eqnarray} 
To see that
$(\calB_\infty)^0$ is connected note that the map
$$(\pr_\infty)^0:(\calB_\infty)^0\to \calC$$ has connected image and
fibres. 
Note also 
that $\calB_\infty$ is contractible from Corollary~\ref{loctriv}. 

Thus
(\ref{homologyvanish}) follows from
$$H_*(\calB_\infty,(\calB_\infty)^0;\Z)=0.$$ Taking direct limits it follows from
$$H_*((\calB_\infty)^{\leq d},(\calB_\infty)^0;\Z)=0 \ {\rm for\ each\
}d.$$ We prove this by
induction on $d$. For $d=0$ it is trivial. Suppose we proved it for
$d$ and consider the homology long exact sequence of the triple 
$(\calB_\infty)^0\subset (\calB_\infty)^{\leq
d}\subset(\calB_\infty)^{\leq d+1}$: 
$$\rightarrow H_q\left((\calB_\infty)^{\leq d},(\calB_\infty)^{0};\Z\right)
\rightarrow H_q\left((\calB_\infty)^{\leq d+1},(\calB_\infty)^{0};\Z\right)
\rightarrow H_q\left( (\calB_\infty)^{\leq
d+1},(\calB_\infty)^{\leq d};\Z\right)\rightarrow.$$
By induction $H_*((\calB_\infty)^{\leq d},(\calB_\infty)^0;\Z)=0$ thus we
need only to prove  $$H_*\left( (\calB_\infty)^{\leq
d+1},(\calB_\infty)^{\leq d};\Z\right)=0.$$ Taking direct limits it follows from
Corollary~\ref{relativevanish} since $2g-2+k \to \infty$ as $k \to \infty$. 
Thus $$H_*\left(\calB_\infty,(\calB_\infty)^0;\Z\right)=0$$ indeed,
proving (\ref{homologyvanish}).

We also show that \begin{eqnarray}
\pi_1((\calB_\infty)^0)\ \rm{is \
Abelian}.\label{abelian}\end{eqnarray} 
Consider the homotopy long exact
sequence of the fibration $(\calB_\infty)^0\to \left((\calB_\infty)^0\right)_\bG\to
B\bG$: $$\pi_2(B\bG)\to \pi_1\left((\calB_\infty)^0\right)\to 
\pi_1\left(\left((\calB_\infty)^0\right)_\bG\right)\stackrel{\pr_*}{\to} \pi_1(B\bG),$$
Because the indices  of the Bott-Morse function $\mu_k$ 
--the moment map of the $U(1)$-action on $\tM_k$-- are all even it
follows from Bott-Morse theory that
$$\pi_1\left(\left((\calB_k)^0\right)_\bG\right)\cong \pi_1(\tM_k)\cong
\pi_1(\tN).$$ According to p. 581 of \cite{atiyah-bott}
$$\pi_1(\tN)\cong \pi_1(B\G)\cong \pi_1(B\bG).$$ Thus $\pr_*$ is an
isomorphism. Thus $\pi_1\left((\calB_\infty)^0\right)$ is a factor group of
the Abelian group $\pi_2(B\bG)$, proving (\ref{abelian}). 

Now (\ref{homologyvanish}) together
with the Hurewitz theorem\footnote{As in Theorem 2.1.1 of Section 2 of
Chapter  13 of \cite{james}.} imply that the abelianization
of $\pi_1\left((\calB_\infty)^0\right)$ is $0$ thus from
(\ref{abelian}) it is $0$, and in turn we get (\ref{homotopyvanish}).

The next step is to show that $(\calB_\infty)^0$ is a CW-space\footnote{Which
means that it has the homotopy type of a CW-complex}. 
Consider the fibration \begin{eqnarray}(\calB_\infty)^0\to
((\calB_\infty)^0)_{\bG^c} \to B\bG \label{fibration1}\end{eqnarray} 
from (\ref{eqfibration1}). We show that its total space and base space are CW-spaces. It will then
follow from Corollary (13) of \cite{stasheff}, that the fibre
$(\calB_\infty)^0$ is a CW-space. 

Note that  $\bG^c$ acts freely on
$(\calB_\infty)^0$ and the quotient is $\tM_\infty$. 
Thus we have the fibration $$E\bG^c\to ((\calB_\infty)^0)_{\bG^c}\to
\tM_\infty.$$
In this fibration the base space, being the direct limit of finite
dimensional manifolds, is a CW-complex, and the fibre, being
contractible, is a CW-space. Then Proposition (0) 
of \cite{stasheff} yields that the total space  
$((\calB_\infty)^0)_{\bG^c}$ is a CW-space. 

Furthermore $B\G$ is a CW-space 
because according to Proposition 2.4 of \cite{atiyah-bott}
it is a component of a mapping space from a compact Hausdorff space $\Sigma$ to
a CW-complex $BU(2)$, and a
theorem of Milnor\footnote{Cf. \cite{milnor}} says that such a
space is a CW-space. Recall now
(\ref{product}), i.e. that $B\G\sim BU(1)\times B\bG$. Thus we have a fibration
$B\G\to BU(1)$ from a CW-space to a CW-complex, according to
Corollary (13) of \cite{stasheff} it follows that the fibre $B\bG$ is itself a CW-space.  

Putting everything together we have a connected CW-space $(\calB_\infty)^0$ with
trivial homotopy groups. Whitehead's theorem\footnote{Cf. e.g. Theorem
2.1.3 of Section 2 of Chapter 13 of \cite{james}.} concludes the proof.
\end{proof}

Thus we have $\bG^c$ acting freely on the contractible space
$(\calB_\infty)^0$, with quotient $\tM_\infty$, which gives the
following immediate

\begin{corollary}
The space $\tM_\infty$ is
homotopy equivalent to $B\bG$.
\label{homotopyminfty}
\end{corollary}

\subsection{Homotopy type of the strata}
\label{homotopystrat}

In this subsection we prove that not only the whole spaces $\tM_\infty$ and $B\bG$ are
homotopy equivalent, but even as stratified spaces. This explains why
the calculation (\ref{poincareminfty}) was the same as the Atiyah-Bott
calculation (\ref{poincaretn}) of $P_t(\tN)$.

\begin{proposition} The map $$(\pr_\infty)^0_d:(\calB_\infty)^0_d\to
\calC_d$$ is a homotopy equivalence. 
\label{weakstratum}
\end{proposition}
\begin{proof} First we note that
$(\calB_\infty)^0_d=(\calB_\infty)_d\setminus (\calB_\infty)^d$. Now
both are locally trivial fibrations with fibre $\C^\infty$ over $\calC_d$ and the codimension
of $(\calB_\infty)^d$ in $(\calB_\infty)_d$ is infinite. It follows
that  $(\calB_\infty)^0_d$ is a locally trivial fibration over
$\calC_d$ with fibre retracting to $S^\infty$, which is
contractible. Now $\calC_d$ is paracompact (and consequently
numerable), 
since it is  a metric subspace
of the metric space $\calC$, and also locally contractible since it
is a Banach manifold. Thus Theorem 6.3 of \cite{dold}\footnote{This
reference was suggested by Ioan James.} yields that
$(\pr_\infty)^0_d$ is a homotopy equivalence.  
\end{proof}

\begin{corollary} We have the homotopy equivalence
\begin{eqnarray} \widetilde{U}^\infty_d\sim (\calC_d)_{\bG^c}.\label{marcsakketnap}\end{eqnarray}
Consequently  $(\calC_d)_{\bG^c}$ is homotopy equivalent to $\J_d\times \Sigma_\infty$.
\end{corollary}
\begin{proof} The map $(\pr_\infty)^0_d$ induces a map of
$B\bG$-spaces: 
$$\begin{array}{ccc} (\calB_\infty)^0_d &\to &\calC_d\\
				\downarrow &  & \downarrow \\
		\left((\calB_\infty)^0_d\right)_{\bG} &\to & (\calC_d)_{\bG}\\
                               \downarrow &     & \downarrow \\
				B\bG &   \cong  &    B\bG
                  \end{array}. $$
As one can choose a CW-complex model for $B\bG$, e.g. from 
Corollary~\ref{homotopyminfty}, the previous theorem and Theorem 6.3
of \cite{dold} gives (\ref{marcsakketnap}).

For the last statement note that $\widetilde{U}^k_d$ is the moduli
space of complexes $L\stackrel{\phi}{\to} VKL_p^k$, which is
uniquely 
determined by $L$ and $\phi$, thus $\widetilde{U}^k_d\cong \J_d\times
\Sigma_{\bd+k}$.
It follows that \begin{eqnarray}\widetilde{U}^\infty_d\cong\J_d\times
\Sigma_\infty.\label{infty}\end{eqnarray}
The result follows.    
\end{proof}

\begin{corollary} We have $$(\calC_d)_\G\sim \Sigma_\infty\times \Sigma_\infty.$$
\label{homotopycdg}
\end{corollary}
\begin{proof} Since (\ref{product}) is a product we get $$(\calC_d)_\G\sim
BU(1)\times (\calC_d)_\bG\cong BU(1)\times \J_d \times \Sigma_\infty.$$
The result follows from (\ref{doldthom}).
\end{proof}

\subsection{Stabilization of homotopy groups}

Finally we have two results about homotopy groups stabilizing in the
resolution tower.  The second of which is reminiscent of the
the Atiyah-Jones conjecture about the stabilization of the homotopy
groups of the moduli space of instantons on $S^4$.

\begin{theorem} For $k\geq 0$ we have $$(\pr_d^k)_*:\pi_i(\widetilde{U}^k_d)\to
\pi_i(\widetilde{U}^\infty_d) \mbox{\rm\ \ for \ } 0\leq i \leq \bd+k-1.$$ 
\end{theorem} 

\begin{proof} Because of (\ref{infty}), it is sufficient to show the stabilization of homotopy
groups for the resolution tower (\ref{resolutionsymmetric}) of $\Sigma$. 
By (12.2) of \cite{macdonald} 
we have an isomorphism \begin{eqnarray}(i_n)_*:H_i(\Sigma_n;\Z)\to H_i(\Sigma_\infty;\Z)
\mbox{\rm \ \ 
for \ } 0\leq i \leq n-1.\label{isosymm}\end{eqnarray} 
The isomorphism for the fundamental groups
is clear. Thus
$\pi_2(\Sigma_\infty,\Sigma_n)$, being the factor group of
the Abelian group $\pi_2(\Sigma_\infty)$, 
is also Abelian. Now the relative
Hurewitz theorem gives the result. 
\end{proof}
\begin{remark} 1. By a theorem of  Dold and Thom we have a complete
description
 of the homotopy type of $\Sigma_\infty$, namely
$\pi_k(\Sigma_\infty)\cong H_k(\Sigma)$, and 
\begin{eqnarray}\Sigma_\infty\sim \prod_{i>0}
K(H_i(\Sigma;\Z),i)\label{doldthom},\end{eqnarray}
which combined with (\ref{infty}) gives an explicit description of the homotopy type of
$\widetilde{U}_d^\infty$. 

2. It is also interesting to note that (\ref{doldthom}) together
with Proposition 2.4 and (2.6) of \cite{atiyah-bott} show that
$B\G_1\sim \Sigma_\infty,$ where $\G_1$ is the group of gauge
transformations on a principal $U(1)$-bundle on $\Sigma$.
\end{remark}

Our final result is the following

\begin{theorem} For $k\geq 0$ we have
$$\pi_{i}(\tM_k)\stackrel{i_k^*}{\cong}\pi_{i}(\tM_\infty)\cong
\pi_{i}(B\bG) \mbox{\rm \ \ for \ } 0\leq i \leq 4g-8+k.$$
\end{theorem}

\begin{proof} 
First  we show that $$H_i(\tM_{\infty},\tM_{k};\Z)=0 \mbox{\rm \ \
for \ } 0\leq i \leq 4g-7+k.$$ 
This follows from \begin{eqnarray} (i_{\tM_k})_*:H_i(\tM_{k};\Z)\to
H_i(\tM_{\infty};\Z) \mbox{\ is an isomorphism
for\ } 0\leq i \leq 4g-7+k,\nonumber \end{eqnarray} which is a
consequence of the five lemma applied to the diagram 
$$\begin{array}{ccccccccc} \dots & \to & H_q(\widetilde{U}^\infty_{<d};\Z) & \to &
H_q(\widetilde{U}^\infty_{\leq d};\Z) & \to &
H_{q-2(g+2d-2)}(\widetilde{U}_d^\infty;\Z) & \to & \dots \\
&&\uparrow &&\uparrow&& \uparrow && \\
 \dots & \to & H_q(\widetilde{U}^k_{<d};\Z) & \to &
H_q(\widetilde{U}^k_{\leq d};\Z) & \to &
H_{q-2(g+2d-2)}(\widetilde{U}_d^k;\Z) & \to & \dots
\end{array}$$
and (\ref{isosymm}).   
We also need that $\pi_1(\tM_{\infty},\tM_{k})=0$, 
this follows
from the fact that $(i_0)_*:\pi_1(\tN)\to \pi_1(\tM_{k})$ is an
isomorphism for each $k$
from standard Bott-Morse theory, since each index is even. Finally we
have 
that $\pi_2(\tM_{\infty},\tM_k)$ is
Abelian, because it is a factor of
the Abelian group $\pi_2(\tM_{\infty})$. 
Now the relative Hurewitz 
theorem\footnote{E.g. Theorem 2.1.2 of Section 2 of Chapter 13 of
\cite{james}.} 
gives that $$\pi_i(\tM_{\infty},\tM_{k};\Z)=0 \mbox{\rm \ \
for \ } 0\leq i \leq 4g-7+k,$$ which in turn proves the result. 
\end{proof}

\backmatter 

\setcounter{chapter}{8}
\chapter{Conclusion}
\label{Conclusion}

In this thesis we have  attempted to give a general picture of the
geometrical and topological properties of $\M$, the moduli space of
rank $2$ Higgs bundles with fixed determinant of degree $1$ over $\Sigma$.
Examining the symplectic geometry of $\M$ we found two Morse stratifications on
it and a natural compactification of it. The downward flows were found
to be 
responsible for the intersection numbers, and the upward flows for the
cohomology ring. Investigating the latter we constructed a resolution
tower for $\tM$ and found that its direct limit was a model for the
classifying space of the gauge group modulo constant scalars. 

However we have not yet explained the relation between the
compactification and the rest of the thesis. In the next and final section we
intend to fill this gap by providing a heuristic and at some places
conjectural summary of the thesis
from the point of view of the compactification.

\setcounter{section}{0}
\section{Compactification of the thesis}
\label{compactthesis}

Let us go back to the end of Chapter~\ref{compact} and recall
Theorem~\ref{fura}, where we transformed the problem of
intersection numbers from $\M$ to a problem concerning the cohomology
ring of $Z$. Without using this correspondence we were able to
calculate the intersection numbers in Chapter~\ref{intersection}. 
Here however we are  focusing on $Z$ and
explain the cohomological calculations of the later chapters from this angle. 

Since $Z$ is a symplectic quotient of $\M$ we have the Kirwan map
$$r:H^*_\circ(\M)\to H^*(Z),$$ which has the fundamental property that
it is surjective. Consequently generators for $H^*_\circ(\M)$ give
generators for $H^*(Z)$. We denote by $\al_Z$, $\be_Z$ and $\psi^i_Z$
the images of the corresponding equivariant universal classes by the
Kirwan map. We have also that $r(u)=c_1(L_Z)$ the first Chern class of
the contact line bundle on $Z$. Moreover these generators can be obtained
from the universal bundle $\E_Z$, which is the restriction of
$\E_\M^\circ$ in the quotient.

In order to be able to go on we have to consider $H_{\circ cpt}^*(\M)$ the
compactly supported equivariant cohomology of $\M$. It is trivial
below the middle dimension and is $g$-dimensional at the middle
dimension. Indeed $H^{3g-3}_{\circ cpt}(\M)$ is generated by the
equivariant compactly supported cohomology classes of the components
of the nilpotent cone. An analogue of Theorem~\ref{eqgeneration} for
compactly supported cohomology says that 
$H_{\circ cpt}^*(\M)^\Gamma$   is generated\footnote{Recall
that $H^*_{cpt}$ is an $H^*$-module.} as a $\Qu[u,\al,\be,\psi_i]$-module
from $H^{3g-3}_{\circ cpt}(\M)$. 

Going back to $Z$ we have, as the $\C^*$-equivariant analogue of 
Theorem~\ref{fura}, the isomorphism\footnote{This statement
seems to be true in the general setting of Subsection~\ref{rizsa2}. As
such, it may be interesting in its own right, especially in relation
with the recent paper \cite{tolman-weitsman}.} of rings: 
$$H^*(Z)\cong H^*_{\circ}(\M)/H^*_{\circ cpt}(\M),$$
where on the right hand side we have the quotient ring of the ordinary
equivariant cohomology by the image of the
compactly supported equivariant cohomology, which is an ideal.
Considering the obvious equivariant structure on the virtual Dirac
bundle $\D_k$, we can work out the equivariant cohomology classes of
the components of the nilpotent cone. The corresponding $\D_k$ over
$Z$ however, will be an honest vector bundle, thus its Chern classes
vanish in degrees beyond the rank, giving relations in $H^*(Z)^\Gamma$. We know that each
component of the nilpotent cone has trivial cohomology class, which
over $Z$ says, that all the above relations will be some multiple of $c_1(L_Z)$.

The picture which emerges is that the relations for
$H^*(Z)^\Gamma$ are of two types: those which are multiples of
$c_1(L_Z)$, these correspond to the intersection numbers on $H^*(\M)$,
and the rest, which correspond to relations in $H^*(\M)^\Gamma$. Thus
from this perspective
Chapter~\ref{intersection} and Chapter~\ref{cohomology} attempt to
give a complete description of the cohomology ring of the projective
variety $Z$ and in turn for the compactification $\cM$!

We can also look at Chapter~\ref{minfty} from the point of view of the
compactification. Namely we can form $\widetilde{Z}_k$ the highest level K\"ahler
quotients of each $\tM_k$, since their moment maps $\mu_k$ are proper,
by taking the quotient of $\tM_k\setminus \tilde{N}^k$ by the
$\C^*$-action, 
where $\tilde{N}^k$ denotes the downward Morse flow, or
equivalently the nilpotent cone in $\tM_k$. 
It can be seen that the inclusions of the spaces $\tM_k$ induces
inclusions for $\widetilde{Z}_k$. Thus we can form the direct limit
$\widetilde{Z}_\infty$. Similarly we have $\cM_k$ the compactification
of $\tM_k$ and their direct limit $\cM_\infty$. 

The rational cohomology of $\widetilde{Z}_\infty$ is generated by the
universal classes and an extra degree $2$ class:
$c_1(L_{\widetilde{Z}_\infty})$, the first Chern class of the contact
line bundle on $\widetilde{Z}_\infty$. Its Poincar\'e polynomial can
be shown to be equal to the Poincar\'e polynomial of the free graded
commutative algebra on these generators, showing that
$H^*(\widetilde{Z}_\infty)$ is a free graded commutative algebra. It
follows that it is isomorphic to $H^*(B\G)$ the rational cohomology of
the classifying space of the {\em whole} gauge group. We suggest that
this is because they are both homotopy equivalent to $\cM_\infty$:
\begin{eqnarray}\cM_\infty \sim\widetilde{Z}_\infty\sim
B\G.\label{conj}\end{eqnarray}
The main technical difficulty which arises in attempting to prove
this is the fact that the spaces $\widetilde{Z}_k$ and $\cM_k$ 
are not smooth, they have $\Z_2$-orbifold singularities, thus the calculation
of their homology with integer coefficients is a bit subtle. 

To avoid this problem one may want to consider the highest level
K\"ahler quotient and the compactification only
homotopically. In other words let us take $\widetilde{Z}_k^\prime$ to
be $(\tM_k\setminus \tilde{N}^k)_{U(1)}$, i.e. take the homotopy
quotient instead of the singular topological quotient. In an analogue
way we form the homotopy-simplectic cut $\cM_k^\prime$. 
Then it can be
shown, without problems coming from torsion, that we have
$$\widetilde{Z}_\infty^\prime\sim \cM_\infty^\prime \sim (\tM_\infty)_{U(1)}.$$
Now (\ref{eqfibration1}) gives the fibration $$\tM_\infty \to
(\tM_\infty)_{U(1)} \to BU(1),$$
which we propose to be homotopy equivalent to the fibration
(\ref{groupfibration}), in particular, that it is a product. 

Since the codimension of the singular locus in
$\widetilde{Z}_\infty$ is $\infty$, one hopes to conclude that
$\widetilde{Z}_\infty\sim \widetilde{Z}_\infty^\prime$ and similarly
$\cM_\infty\sim \cM_\infty^\prime$ yielding (\ref{conj}).

\newpage
\thispagestyle{empty}



\begin{thebibliography}{ZZZZZ}

\label{bib}

\bibitem[ACGH]{arbarello-et-al} E. Arbarello, M. Cornalba, P.A. Griffiths,
J. Harris, {\em Geometry of algebraic curves I} Springer, Berlin, 1985

\bibitem[Ati1]{atiyah} M. Atiyah, {\em Geometry of Yang-Mills fields.}
Fermi lectures, Scuola Normale 
Pisa, 1979

\bibitem[Ati2]{atiyah2} M. Atiyah, Vector bundles over an elliptic
curve. {\em Proc. Lond. Math. Soc} {\bf 7} (1957) 414-452 

\bibitem[At,Bo]{atiyah-bott} M. Atiyah, R. Bott. 
The Yang-Mills equations over Riemann surfaces.
{\em Philos. Trans. Roy. Soc. London} Series A, {\bf 308} (1982) 523-615 

\bibitem[At,Hi]{atiyah-hitchin} M. Atiyah, N. Hitchin. {\em The
geometry and dynamics of magnetic monopoles}, Princeton University
Press, 1987



\bibitem[ADHM]{atiyah-etal-1} M. Atiyah, V. Drinfeld, N. Hitchin,
Y. Manin, Construction of instantons, {\em Phys. Lett.} {\bf 65A} (1978) 185-187

\bibitem[A,H,S]{atiyah-etal-2} M. Atiyah, N. Hitchin, I. Singer, 
Self-duality in four dimensional Riemannian geometry, {Proc. R. Soc.
London Ser. A} {\bf 362} (1978) 425-461

\bibitem[Bai]{bailey} W.L. Bailey, Jr. On the embedding of V-manifolds in
projective space. {\em Americal Journal of Mathematics} {\bf 79} 
(1957) 403-430 


\bibitem[Bar]{baranovsky} V. Baranovsky, Cohomology ring of the moduli space of
stable vector bundles with odd determinant, {\em
Izv. Russ. Acad. Nauk} 58/4 (1994) 204-210

\bibitem[BJSV]{bershadsky-et-al} M. Bershadsky, A. Johansen, V. Shadov,
C. Vafa. Topological reduction of $4$D SYM to $2$D
$\sigma$-models. {\em Nuclear Phys. B}  {\bf 448} No.1-2 (1995) 166-186

\bibitem[B,P,V]{barth-peters-van} W. Barth, C. Peters, A. Van de Van. 
{\em Compact complex surfaces}, Springer Verlag, 1984 

\bibitem[Bia]{bialynicki} A. Bialynicki-Birula. Some theorems on 
actions of algebraic groups, {\em Ann. of Math.} {\bf 98} (1973) 480-497

\bibitem[Bi,Da]{bielawski-dancer} R. Bielawski, A.S. Dancer, The geometry
and topology of toric hyperk\"ahler manifolds. Preprint, 1996

\bibitem[Bi,Ra]{biswas-ramanan}
I. Biswas, S. Ramanan. An infinitesimal study of the
moduli space of Hitchin pairs, {\em J. London Math. Soc.} {\bf 49}
(1994) 219-231

\bibitem[Br,Pr]{brion-procesi} M. Brion, C. Procesi, 
Action d'un tore dans une vari\'et\'e 
projective, {\em Operator algebras, unitary representations, enveloping 
algebras, and invariant theory} (A.Connes, M.Duflo, A. Joseph, and 
R.Rentschler, eds.), Birkhauser, (1990) 509-539

\bibitem[Das]{daskalopoulos} 
G.D. Daskalopoulos. The topology of the space of stable bundles on a
compact Riemann surface. {\em J. Diff. Geo.} {\bf 36} No.3 (1992) 699--746 

\bibitem[DeRh]{derham} De Rham. {\em Differential manifolds}. English
edition, Berlin, Heidelberg, New York: Springer, 1984

\bibitem[Dod]{dodziuk} J. Dodziuk. Vanishing theorems for square integrable
harmonic forms, {\em Proc. Indian. Acad. Sci. Math. Sci.} {\bf 90} 
(1981) 21-27

\bibitem[Dold]{dold} A. Dold. Partitions of unity in the theory of
fibrations. {\em Annals of Math.} {\bf 79} No.2 1963

\bibitem[Do,Th]{dold-thom} A. Dold, R. Thom. Quasifaserungen und
unendliche
symmetrische Produkte, {\em Ann. of Math.} {\bf 67} (1958) 230-281

\bibitem[Do,Ma]{donagi-markman} R. Donagi, E. Markman. Spectral
covers, algebraically completely integrable, Hamiltonian systems, and
moduli of bundles, 1-119 in {\em Integrable Systems and Quantum
Groups},
Springer, 1996 

\bibitem[Do,Kr]{donaldson-kronheimer} S.K. Donaldson,
P.B. Kronheimer. {\em The geometry of four manifolds}, Oxford
University Press, 1990

\bibitem[Dr,Na]{drezet-narasimhan} J.-M. Drezet, M.S. Narasimhan,
 Groupe Picard des vari\'et\'es
de modules de fibr\'es semistables sur les courbes alg\'ebraiques, 
{\em Invent. Math.} {\bf 97} 53-94 (1989)

\bibitem[Du,He]{duistermaat-heckman} J.J. Duistermaat, 
G.J. Heckman. On the variation in the cohomology of the symplectic form of
the reduced phase space. {\em Invent. Math.} {\bf 69} (1982) 259-268 

\bibitem[Earl]{earl} R. Earl, {\em The Mumford relations and the
moduli of low rank stable bundles}, D. Phil thesis, Oxford, 1995

\bibitem[EdGr1]{edidin-graham} D. Edidin, W. Graham. Algebraic cuts, Preprint, 
alg-geom/9608028 

\bibitem[EdGr2]{edidin-graham2} D. Edidin, W. Graham. Equivariant
intersection theory. {\em Invent. Math.} {\bf 131} (1998) 595-634

\bibitem[Goth1]{gothen1} P.B. Gothen, The Betti numbers of the moduli space
of rank $3$ Higgs bundles, {\em Internat. J. Math.} {\bf 5} (1994) 861-875 

\bibitem[Goth2]{gothen2} P.B. Gothen, Private communication

\bibitem[Goto]{goto} R. Goto. On Toric hyperk\"ahler manifolds given by
hyperk\"ahler quotient method, in `Infinite Analysis', {\em Advanced
Series in Mathematical Physics} {\bf 16} (1992) 317-388 

\bibitem[Gr,Ha]{griffiths-harris} P. Griffiths, J. Harris. {\em Principles
of algebraic geometry}. New York, Wiley 1978

\bibitem[Gro]{grothendieck} A. Grothendieck. Sur la classification des
fibr\'es holomorphes sur la sph\'ere de Riemann. {\em Am. J. Math.}
{\bf 79} (1957) 121-138

\bibitem[Gu,St]{guillemin-sternberg} V. Guillemin, S. Sternberg, Birational
equivalences in the symplectic category. {\em Invent. Math} {\bf 97} (1989)
485 

\bibitem[Gun]{gunning} R.C. Gunning, {Lectures on vector bundles over
Riemann surfaces} Princeton University Press, 1967

\bibitem[Ha,Na]{harder-narasimhan} G. Harder, M.S. Narasimhan. 
On the cohomology groups of moduli spaces of vector bundles over
curves. {\em Math. Ann.} {\bf 212} (1975) 215-248

\bibitem[Har]{hartshorne} R. Hartshorne, {\em Algebraic geometry} 
Springer-Verlag, New York, 1977

\bibitem[Hau1]{hausel1} T. Hausel. Compactification of moduli
of Higgs bundles,   to appear in 
{\em J. reine angew. Math.}, Preprint math.AG/9804083,

\bibitem[Hau2]{hausel2} T. Hausel. Vanishing of intersection 
numbers on the 
moduli space of Higgs bundles, Preprint math.AG/9805071

\bibitem[Ha,Th]{hausel-thaddeus} T. Hausel, M. Thaddeus, Cohomology
ring of the moduli space of Higgs bundles, under preparation

\bibitem[Hit1]{hitchin1} N. Hitchin. The self-duality equations on a Riemann 
surface, {\em Proc. London Math. Soc.} (3) {\bf 55} (1987) 59-126

\bibitem[Hit2]{hitchin2} N. Hitchin. Stable vector bundles and integrable systems, {\em Duke Mathematical Journal} (1) {\bf 54} (1987) 91-114

\bibitem[Hit3]{hitchin3} N. Hitchin. Flat connections and geometric 
quantization, {\em Comm. Math. Phys.} {\bf 131} (1990) 347-380 

\bibitem[Hit4]{hitchin4} N. Hitchin. Private communication

\bibitem[HKLR]{hitchin-et-al} N. Hitchin, A. Karlhede, U. Lindstrom,
M. Ro\v cek, Hyper-K\"ahler metrics and supersymmetry,
Comm. Math. Phys. {\bf 108} (1987) No.4 535-589

\bibitem[Hu,Le]{huybrechts-lehn} 
D. Huybrechts, M. Lehn. {\em The geometry of moduli 
spaces of sheaves}, Friedr. Vieweg und Sohn, Braunschweig, 1997

\bibitem[Jam]{james} I.M. James (editor). {\em Handbook of Algebraic
Topology}, North-Holland, 1995 

\bibitem[Ki,Ne]{king-newstead} A.D. King, P.E. Newstead, On the cohomology ring of
the moduli space of rank $2$ vector bundles on a curve, {\em Topology}
{\bf 37} No.2 (1998) 407-418

\bibitem[Kir1]{kirwan1} F.C. Kirwan. {\em Cohomology of quotients in
symplectic and algebraic geometry}, Mathematical Notes 31, Princeton University
Press, 1984

\bibitem[Kir2]{kirwan2} F. Kirwan. The cohomology ring of moduli space
of bundles over Riemann surfaces. {\em J. Amer. Math. Soc.} 5 (1992) 
No.4 853-906

\bibitem[Lan]{lang} S. Lang. {\em Algebra} Addison-Wesley, Reading, 
Massachusetts, 1965

\bibitem[Ler]{larman} E. Lerman. Symplectic Cuts, {\em Math. Res. Lett.} 
{\bf 2} (1995) 247-258

\bibitem[Lau]{laumon} G. Laumon. Un analogue du cone nilpotent.
{\em Duke Math. J.} (2) {\bf 57} (1988) 647-671 

\bibitem[Leb]{lebrun} C. Lebrun. Fano manifolds, contact structures and 
quaternionic geometry, {\em International Journal of Mathematics} (3) {\bf 6}
(1995) 419-437 

\bibitem[Macd]{macdonald} I.G. Macdonald, Symmetric products of an algebraic
curve, {\em Topology} {\bf 1} (1962) 319-343

\bibitem[Miln]{milnor} J.W. Milnor. Spaces having the homotopy type of
a CW complex, {\em Trans. Am. math. Soc.} {\bf 90} (1959) 272-280

\bibitem[M,F,K]{fogarty-etal} D. Mumford, J. Fogarty, F. Kirwan,
{\em Geometric invariant theory}, third edition, Springer-Verlag, Berlin, 1994

\bibitem[Nak]{nakajima} H. Nakajima. {\em Lectures on Hilbert schemes of points
on surfaces.} to appear

\bibitem[Na,Se]{narasimhan-seshadri} M.S. Narasimhan,
C.S. Seshadri. Stable and unitary vector bundles on a compact Riemann
surface, {\em Ann. Math.} (1965) 540-567

\bibitem[New1]{newstead} P.E. Newstead. {\em Introduction to moduli 
problems and orbit spaces}, Tata Inst. Bombay, 1978

\bibitem[New2]{newstead2} P.E. Newstead. Characteristic classes of
stable bundles of rank $2$ over an algebraic curve, {\em
Trans. Amer. Math. Soc.} {\bf 169} July (1972) 337-345

\bibitem[Nit]{nitsure} N. Nitsure. 
Moduli space of semistable pairs on a curve. {\em Proc. London Math. Soc.}
(3) {\bf 62} (1991) 275-300 

\bibitem[Oxb]{oxbury} W. Oxbury. {\em Stable bundles and branched coverings
over Riemann surfaces.} D.Phil thesis, Oxford, 1987

\bibitem[Sche]{scheja} G. Scheja, Riemannsche Hebbarkeitsatze f\"ur
cohomologieklassen, {\em Math. Ann.} {\bf 144} (1961) 345-360 

\bibitem[Schm]{schmitt} A. Schmitt. Projective moduli of Hitchin pairs, 
{\em Int. J. of Math.} {\bf 9}, No.1 (1998) 107-118 (available at alg-geom/9611008).

\bibitem[Se,Se]{segal-selby} G. Segal, A. Selby. Cohomology of the
space of magnetic monopoles, {\em Comm. Math. Phys.} {\bf
177} (1996) 775-787

\bibitem[Sen]{sen} A. Sen. Dyon-monopole bound states, self-dual
harmonic forms on the multi-monopole moduli space, and $SL(2,\Z)$
invariance in string theory. {\em Phys.-Lett.-B} {\bf 329}  No.2-3
(1994) 217-221

\bibitem[Si,Ti]{siebert-tian} B.Siebert, G.Tian: Recursive relations
for the cohomology ring of moduli space of stable bundles, {\em
Turkish J. Math.} {\bf 19} (1996) 131-144

\bibitem[Sim1]{simpson} C. Simpson. Higgs bundles and local systems. {\em 
Publ. Math. I.H.E.S.} {\bf 75} (1992) 5-95

\bibitem[Sim2]{simpson2} C. Simpson. The Hodge filtration on
nonabelian cohomology, alg-geom/9604005

\bibitem[Sim3]{simpson3} C. Simpson. Nonabelian Hodge Theory. {\em Proceedings
of the International Mathematical Congress, Kyoto}, 1990, 747-756

\bibitem[Sim4]{simpson4} C. Simpson. 
The ubiquity of variations of Hodge structure. 
{\em Complex geometry and Lie theory (Sundance, UT, 1989)}, 
{\em Proc. Sympos. Pure Math.}, {\bf 53}, Amer. Math. Soc.,
Providence, 
RI, 1991, 329-348

\bibitem[Sta]{stasheff} J.D. Stasheff. A classification theorem for
fibre spaces. {\em Topology} {\bf 2} (1963) 239-246

\bibitem[Sun]{sundaram} N. Sundaram. Special divisors and vector
bundles. {\em Tohoku Math. J.} (2) {\bf 39} No.2 (1987) 175-213

\bibitem[Tei]{teixidor} M. Teixidor i Bigas, Brill-Noether theory for 
stable vector bundles, {\em Duke Math. J.} {\bf 62} No.2 (1991) 385-400

\bibitem[Tha1]{thaddeus1} M. Thaddeus. Topology of the moduli space of stable bundles
over a compact Riemann surface, 1990, unpublished


\bibitem[Tha2]{thaddeus2} M. Thaddeus. Stable pairs, linear systems and the
Verlinde formula, {\em Invent. Math.} {\bf 117} (1994) 317-353 

\bibitem[Tha3]{thaddeus3} M. Thaddeus. Geometric invariant theory and flips,
{\em Journal of the American Mathematical Society} (3) {\bf 9} (1996) 691-723

\bibitem[Tha4]{thaddeus4} M. Thaddeus. An introduction to the topology of the
moduli space of stable bundles on a Riemann surface, Collection:
Geometry and Physics (Aarhus, 1995), 71-99

\bibitem[Ti,Ya]{yau-tian} 
G. Tian, S.H. Yau. Complete K\"ahler manifolds with zero 
Ricci curvature. I. {\em J.Amer.Math.Soc.} (3) {\bf 3} (1990) 579-609 

\bibitem[To,We]{tolman-weitsman} S. Tolman, J. Weitsman.
The cohomology rings of abelian symplectic quotients. Preprint
math.DG/9807173

\bibitem[Yoko]{yokogawa}
K. Yokogawa. Compactification of moduli of parabolic sheaves and
moduli of parabolic Higgs sheaves; {\em J. Math. Kyoto Univ.} {\bf 33}
(1993) 451-504

\bibitem[Zag]{zagier} D. Zagier. On the cohomology of moduli spaces of rank
$2$ vector bundles over curves, {\em The moduli space of curves},
ed. R. Dijkgraaf, C. Faber and G. van der Geer, Progress in
Mathematics 129, Birkh\"auser, 1995

\end{thebibliography}
\end{document}